%% file: Main.tex
\documentclass[11pt]{article}
\usepackage{fullpage,latexsym,amsmath,geompsfi,amssymb,subfigure}

\newtheorem{theorem}{Theorem}[section]
\newtheorem{prop}[theorem]{Proposition}
\newtheorem{defn}[theorem]{Definition}
\newtheorem{cor}[theorem]{Corollary}
\newtheorem{lem}[theorem]{Lemma}

\renewcommand{\bar}{\overline}
\newcommand{\dis}{\displaystyle}
\newcommand{\al}{\alpha}
\newcommand{\mc}{\mathcal}
\newcommand{\es}{\emptyset}
\newcommand{\ep}{\epsilon}
\newcommand{\ra}{\rightarrow}
\newcommand{\var}{\varepsilon}
\newcommand{\ov}{\overline}
\newcommand{\un}{\underline}
\newcommand{\p}{\partial}
\newcommand{\noi}{\noindent}

\newcommand{\sr}{\stackrel}
\newcommand{\sm}{\setminus}

\input{imsmark.tex}
\begin{document}
\input{title}

\SBIMSMark{1997/4}{March 1997}{}

\input{intro}

\input{prelim}

\input{planetop}

\input{cedisks}

\input{isotopies}

\input{theorem}

\input{examples}

\input{biblio}

\end{document}

%% file: imsmark.tex
\def\SBIMSMark#1#2#3{
 \font\SBF=cmss10 at 10 true pt
 \font\SBI=cmssi10 at 10 true pt
 \setbox0=\hbox{\SBF Stony Brook IMS Preprint \##1}
 \setbox2=\hbox to \wd0{\hfil \SBI #2}
 \setbox4=\hbox to \wd0{\hfil \SBI #3}
 \setbox6=\hbox to \wd0{\hss
             \vbox{\hsize=\wd0 \parskip=0pt \baselineskip=10 true pt
                   \copy0 \break%
                   \copy2 \break%
                   \copy4 \break}}
 \dimen0=\ht6   \advance\dimen0 by \vsize \advance\dimen0 by 8 true pt
                \advance\dimen0 by -\pagetotal
 \dimen2=\hsize \advance\dimen2 by .25 true in
  \openin2=publishd.tex
  \ifeof2\setbox0=\hbox to 0pt{}
  \else 
     \setbox0=\hbox to 3.1 true in{
                \vbox to \ht6{\hsize=3 true in \parskip=0pt  \noindent  
                \input publishd.tex 
                \vfill}}
  \fi
  \closein2
  \ht0=0pt \dp0=0pt
 \ht6=0pt \dp6=0pt
 \setbox8=\vbox to \dimen0{\vfill \hbox to \dimen2{\copy0 \hss \copy6}}
 \ht8=0pt \dp8=0pt \wd8=0pt
 \copy8
 \message{*** Stony Brook IMS Preprint #1, #2 ***}
}

%% file: title.tex
\setlength{\topmargin}{1.5in}
\begin{center}
{\Large \bf Pruning fronts and the formation of horseshoes}

\bigskip
\bigskip

{\large \bf by}

\bigskip
\bigskip

{\Large \bf Andr\'e de Carvalho}

\bigskip

Department of Mathematics\\
University of California\\
Berkeley, CA 94720-3840\\
E-mail: andre@@math.berkeley.edu
\end{center}

\vskip .5 in

\begin{abstract}

\setlength{\baselineskip}{18pt}
Let $f: \pi \ra \pi$ be a homeomorphism of the plane $\pi$.  We define 
open sets $P$, called {\em pruning fronts} after the work of Cvitanovi\'c~
\cite{C}, for which it is possible to construct an isotopy $H: \pi \times 
[0,1] \ra \pi$ with open support contained in $\dis\bigcup _{n \in 
{\Bbb{Z}} } f^{n} (P)$ such that $H(\cdot, 0 ) = f(\cdot)$ and $H(\cdot, 
1) = f_{P} (\cdot)$, where $f_P$ is a homeomorphism under which every 
point of $P$ is wandering.  Applying this construction with $f$ being 
Smale's horseshoe, it is possible to obtain an uncountable family of 
homeomorphisms, depending on infinitely many parameters, going from 
trivial to chaotic dynamic behaviour.  This family is a 2-dimensional 
analog of a 1-dimensional universal family.

\end{abstract}


%% file: intro.tex
\setcounter{section}{-1}

\section{Introduction}

\setlength{\baselineskip}{20pt}

\indent
One of the main concerns in the study of dynamical systems is to 
understand how a family of maps passes from simple to complicated dynamic 
behaviour as we vary parameters.  When the dynamical systems under 
consideration are 1-dimensional, the kneading theory of Milnor and 
Thurston provides a full topological understanding of the transition from 
simple to chaotic behaviour.  In dimension 2, no such theory exists.  In 
fact , it is not clear what restrictions should be imposed on the 
families under consideration in order that understanding them is not too 
hopeless a task.

Families like the H\'enon and the Lozi ones are interesting examples but 
they lack a defining topological characteristic analogous, for example, 
to saying that a 1-dimensional map is {\em unimodal} (i.e., is piecewise 
monotone with exactly one turning point.) 

In this work, we present a method of isotoping away dynamics from a 
homeomorphism of the plane in a controlled fashion.  More precisely, if 
$f: \pi \ra \pi$ is a homeomorphism of the plane $\pi$, we define open 
sets $P$ for which there exists an isotopy $H: \pi \times [0,1] \ra \pi$ 
with (open) support contained in $\dis\bigcup_{n \in {\Bbb{Z}} }f^{n} 
(P)$, such that $H( \cdot, 0 ) = f$ and $H ( \cdot, 1) = f_{P}$, where 
$f_{P}$ is a homeomorphism under which every point of $P$ is wandering.  
Using this construction, with $f$ being Smale's horseshoe, for example, 
it is possible to produce an uncountable family of homeomorphisms of the 
plane, depending on infinitely many parameters, going from trivial 
dynamics (say, only two nonwandering points, one attracting and one 
repelling fixed points) to a full horseshoe.

We call the sets $P$ mentioned above {\em pruning fronts}, after the work 
of P. Cvitanovi\'c \cite{C}.  In \cite{C} they propose sets of symbol 
space for Smale's horseshoe which get ``pruned away'' as we vary 
parameters in a family like the H\'enon one.  Here we give a precise 
definition of pruning fronts and construct the isotopies which ``prune 
away'' the dynamics in $P$.

In forthcoming papers we intend to do two things.  First, for each map 
$f_P$, where $P$ is a pruning front as defined herein, there exists a 
collapsing procedure which produces a ``tight'' map $\varphi_P$ isotopic 
to $f_P$ and with essentially the same dynamics.  More precisely, there 
exists an $f_P$-invariant upper semi-continuous decomposition $G_P$ of the 
sphere $S^2$ (we can extend $f$ to $S^2$ setting $f(\infty )= \infty$ ), 
such that, for every element $g$ of $G_P$, $g$ contains at least one 
element of the nonwandering set of $f_P$ and $h(f_{P};g) = 0$, where 
$h(f_{P};g)$ is the topological entropy of $f_P$ in $g$ as defined by 
Bowen.  $f_P$ then projects to a homeomorphism $\varphi_{P}: K_{P} \ra 
K_{P}$ of the {\em cactoid} $K_{P} = S^{2}/G_{P}$, such that no point of 
$K_P$ is wandering under $\varphi_{P}$ and $h(f_{P}) =h( \varphi_{P} )$.  
Second, we intend to show that the family $\varphi_{P}$ contains the 
Thurston minimal reresentatives in the isotopy classes of $f$ relative to 
periodic orbit collections of $f$.  In other words, we would like to 
show that given a periodic orbit collection $\mc{O}$ of periodic orbits 
of $f$, there exists a pruning front $P=P({\mc{O}})$, such that 
$\varphi_P$ is 
the Thurston minimal representative in the isotopy class of $f$ rel 
$\mc O$.  This last statement should have an algorithmic proof, providing 
another algorithmic proof of Thurston's classification therorem for 
homeomorphisms of surfaces.

The techniques used in the present work are those of point set topology 
of the plane.  In Section 1 we state without proof the 
main background 
results we will need, the most important of which being the Jordan Curve 
Therorem (Theorem \ref{1}) and Whyburn's Separation Theorem (Theorem 
\ref{2}).  In Section 2 we develop the plane toplogy 
tools we will use in 
the remainder of the paper.  In Section 3 we introduce 
the concept of $(c,e)$-{\em disks}, define pruning fronts and prove
some propositions which will be used in Section 5.  In Section 4 we state 
and 
prove some results 
about isotopies of homeomorphisms of the plane, which will also be needed 
in Section 5.  Although these results are 
folkloric, we decided to present them for completeness; the proofs given 
are rather elementary.  Section 5 contains the 
proof of the main theorem, as its title suggests.  Within the first few 
pages we get to define an isotopy which is almost all we need (Proposition 
\ref{43}) and the remainder of the section is devoted to showing how this 
isotopy works and how we fix it in order to get the final isotopy $H$ 
(which depends, of course, on $P$.) It is only in Section 6 that we get to the
second part of the title --- the formation of horseshoes. We present three 
examples of pruning fronts for Smale's horseshoe map. The first of which is,
in fact, a family of such examples and produces, via the main theorem, a 
family of homeomorphisms of the plane whose dynamics mimics that of a  
full unimodal family of endomorphisms of the interval. The second example 
gives rise to a `renormalizable' map, that is, a homeomorphism which 
interchanges two closed disks. The second iterate restricted to each one 
of these disks is again a full horseshoe. Finally, in the third example
we present a pruning front which gives rise to a `lax pseudo-Anosov' 
homeomorphism. Together these examples should suggest different ways 
in which a horseshoe can be formed.

A word about the figures is in order.  One of the hardest things for me 
during the preparation of this work was to translate into precise 
mathematical statements the pictures I had in my mind.  I decided, 
therefore, to add to the text all those pictures I had to draw over and 
over for myself before I understood what were the right mathematical 
statements that described them.  I hope they will be helpful to the 
reader, for as the saying goes, ``a picture is worth a thousand 
words.''

\noindent
{\sc Acknowledgements:} The research presented herein comprises my Ph.D.
dissertation done at the Graduate Center of The City University of New 
York (CUNY) under the supervision of Professor Dennis Sullivan. I would like 
to thank Professor Sullivan for his guidance during the preparation of this
work and the Graduate Center of CUNY for providing a friendly and helpful 
research atmosphere. I had several discussions with Alberto Baider,
Pregrag Cvitanovi\'c, Fred Gardiner, Toby Hall, Michael Handel, Ronnie 
Mainieri, Charles Tresser and Nick Tufillaro and I would like to thank them 
for their help.

%% file: prelim.tex
\section{Preliminaries}

We will denote the 2-dimensional plane
$\pi$ or $\Bbb{R}^2$.  A {\em Jordan curve}\/ $J$ is the
homeomorphic
image of the circle $S^1= \{ (x,y) \in \Bbb{R}^{2}; \ x^{2} + y^{2}=1 \}$ and
a {\em closed arc} $L$ is the homeomorphic image of the closed interval
$[0,1]$, the images of $\{0\}$ and $\{1\}$ being its endpoints.
By an {\em open arc}\/ we will mean the set obtained by taking the
endpoints away from a closed arc.  If $L$ is a closed arc $\sr{\circ}{L}$ 
will denote the corresponding open arc.

The theorems that follow can be found in the books of Newman \cite{Ne},
Moise \cite{M} and Whyburn \cite{Wh}.  Moore's book \cite{Mo} is also a
good reference although a little less palatable. \bigskip

\begin{theorem}[Jordan Curve Theorem] \label{1}
Every Jordan curve
separates the plane into two regions $I$ and $O$ and is the boundary of
each.
\end{theorem}

\begin{defn}
Let $J$ be a Jordan curve and $I$ the bounded region of $\pi \sm J$.  
We call $I$ a {\em Jordan domain} and sometimes refer to it as the {\em 
inner domain} determined by $J$. 
\end{defn}

\begin{theorem}[Separation Theorem (Whyburn)] \label{2}
Let $A$ be compact and $B$ closed
subsets of the plane such that $A \cap B$ is totally disconnected, $a \in
A \sm (A \cap B), \ b \in B \sm (A \cap B)$ and
$\varepsilon$ a positive
number.  Then there exists a Jordan curve $J$ which separates $a$ and
$b$ and is such that $J \cap (A \cup B ) \subset A \cap B$ and every
point of $J$ is at distance less than $\varepsilon$ from some point of $A$.
\end{theorem}

\begin{defn}
Let $U$ be a domain in the plane and $\alpha$
an open (closed) arc whose endpoints lie on $\partial U$ and all others
lie in $U$.  Such an $\alpha$ is called an {\em open (closed) cross-cut}.
\end{defn}

\begin{theorem}
If both endpoints of a cross-cut $\alpha$ in a
domain $U
\subset \pi$ are on the same component of ${\mathcal{C}} U$,
the complement of $U, \ U \backslash \alpha$ has two components and is
contained in the frontiers of both.
\end{theorem}

\begin{cor}\label{3}  Let $J$ be a Jordan
curve, $I$ its inner domain and $\alpha \subset I$ a cross-cut.  Then
$\alpha$ separates $I$ into two Jordan domains $I_{1}$ and $I_2$ whose
boundaries are $L_{1} \cup \alpha$ and $L_{2} \cup \alpha$, where $L_1$
and $L_2$ are the arcs into which the endpoints of $\alpha$ separate $J$.
\end{cor}

\begin{theorem}
\label{4}  Let $f: J_{1} \rightarrow J_2$ be a
homeomorphism between the Jordan curves $J_1$ and $J_2$.  Then it is
possible to extend $f$ to a homeomophism $\tilde{f} : D_{1} \rightarrow
D_{2}$ between the closed disks $D_{1}=J_{1} \cup I_1, \ D_{2}=J_{2} \cup
I_2$ bounded by $J_1$ and $J_2$.
\end{theorem}

\begin{theorem} [Alexander] \label{5}
In $\Bbb{R}^n$, let $B^{n}= \{
x; ||x|| \leq 1 \}$ and $S^{n-1}=\partial B^{n-1} = \{x; ||x||=1 \}$ and
$f: B^{n} \rightarrow B^n$ a homeomorphism such that $f|_{S^{n-1}}
\equiv$ identity.  Then $f$ is isotopic to the identity
through an isotopy that fixes the boundary pointwise.
\end{theorem}
\bigskip

%% file: planetop.tex
\section{Plane Topology}

In this section we will develop some plane topology preliminaries we will
need later on.
\bigskip

\noindent{\sc Notation}:  Unless stated explicitly otherwise, we will use
the following notations:  $J$ will stand for a Jordan curve, $I$ and 
$O$ for its
inner and outer domains respectively, and $D$ for the closed disk $I \cup 
J$.  If $D$ is a closed
disk we will sometimes use $I(D)$ to denote its inner domain.  Subscripts
will match in the obvious way, so that the inner domain determined by the
Jordan curve $J_1$ is $I_1$ and $D_{1}= I_{1} \cup J_1$, etc.

If $k$ is a positive integer $\un{k}$ will stand for the set $\{1,2,
\ldots , k \}$.
\bigskip

\begin{figure}
\begin{center}~
\psfig{file=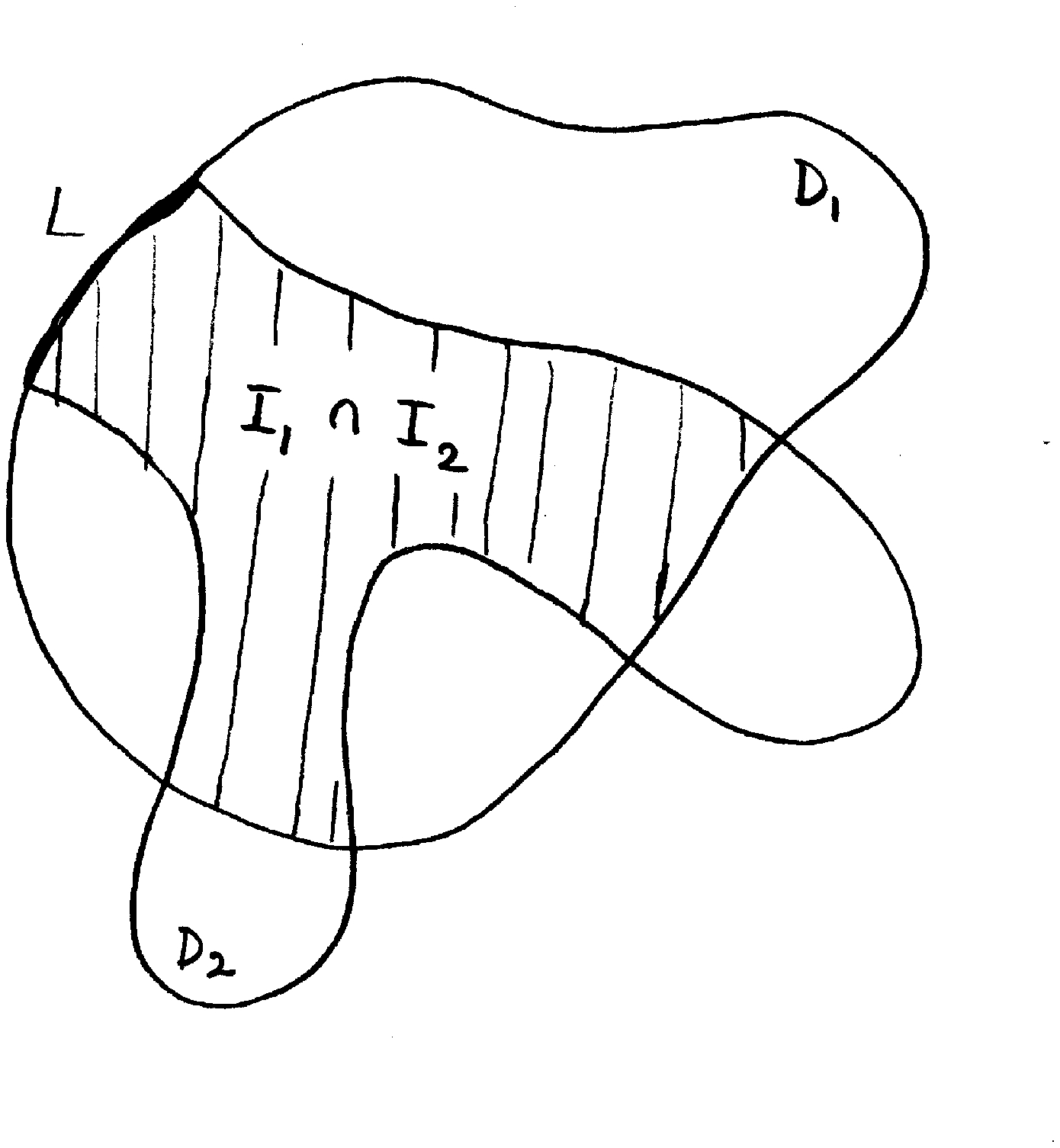,height=2.5in}
\end{center}
\caption{Two disks on the same side of the arc $L$.} \label{f1}
\end{figure}

\begin{defn}
Let $J_{1}, \dots, J_{n}$ be Jordan
curves and $L \subset J_{1} \cap \ldots \cap J_n$ an arc.
We say the closed disks $D_{1}, \ldots , D_n$ {\em lie on the same side
of}\/ $L$, denoted $D_{1}, \ldots , D_{n} |_{L}$, if $L \subset
\overline{I_{1} \cap \ldots \cap I_{n} }$ (see figure~\ref{f1}.)
\end{defn}

\begin{prop} \label{6}
In the plane
$\pi$, let $A$ be a closed arc and $B$ a closed set such
that $A \cap B \subset \{ {\mbox{\rm endpoints of } A} \}$ and there exists
$\varepsilon > 0$ such that every component of $B \sm (A \cap B)$
contains a
point at distance greater than $\varepsilon$ from $A$.  Then there exists
a Jordan curve $J$ such that $A \sm (A \cap B) \subset I$ and $B
\sm (A \cap B) \subset O$
where $I$ and $O$ are the bounded and unbounded components
of ${\mathcal{C}} J$ (the complement of $J$ in $\pi$)
respectively, and $J \cap (A \cup B ) \subset A \cap B \subset \{
\mbox{\rm endpoints of }A \}$.
\end{prop}

\noindent{\sc Proof}: Let $a \in A \sm (A \cap B)$ and $b \in B
\sm (A \cap B)$ such
that $d(b,A) > \varepsilon$.  By Theorem~\ref{2}, there exits a Jordan
curve $J$
separating $a$ from $b$, such that $J \subset V_{\varepsilon}(A)$ (the
$\varepsilon$-neighborhood about $A$) and $J \cap (A \cup B) \subset A
\cap B$.

First notice that $I \subset V_{\varepsilon}(A)$.  This is so because $D
= J \cup I$ is compact and since $A$ is also compact, there exist $x \in
A$ and $y \in D$ which realize $\sup \{d(x,y); \ x \in A, \ y \in D
\}$.  We claim $y \in J$ for if $y \in I$ there would exist $\delta >0$
such that $V_{\delta}(y) \subset I$ and in $V_{\delta}(y)$ there must be
a point whose distance to $x$ is greater than $d(x,y)$.  This shows that
if $J$ is contained in $V_{\varepsilon}(A)$ then so is $D=J \cup I$.

Since $b \notin V_{\varepsilon}(A), \ b \in O$ and since $J$ separates
$a$ from $b, \ a \in I$.  But $A \sm (A \cap B)$ is a connected point
set disjoint from $J$ and $a \in A \sm (A \cap B)$ so that $A
\sm (A \cap B)
\subset I$.  Also, we assumed that each connected component of $B
\sm (A \cap
B)$ had a point outside of $V_{\varepsilon}(A)$, and therefore in $O$.
Since $B \sm (A \cap B)$ is disjoint from $J, \ B \sm (A
\cap B) \subset O$, as we wanted. $\Box$
\bigskip

In the proofs of the statements that follow, indexed unions and 
intersections will be assumed to range from $i=1$ to $i=n$.
\bigskip

\begin{cor} \label{8}
Let $J_{1}, \ldots , J_{n}$ be Jordan curves and $L \subset 
\dis\bigcap^{n}_{i=1}
 J_i$ a closed arc.  Then there exists a Jordan curve $J$ such that
${\stackrel{\circ}{L}} \subset I$ and such that $\left( \dis\bigcup^{n}_{i=1}
J_{i} \right) \sm L \subset O$.
\end{cor}

\begin{figure}
\begin{center}~
\psfig{file=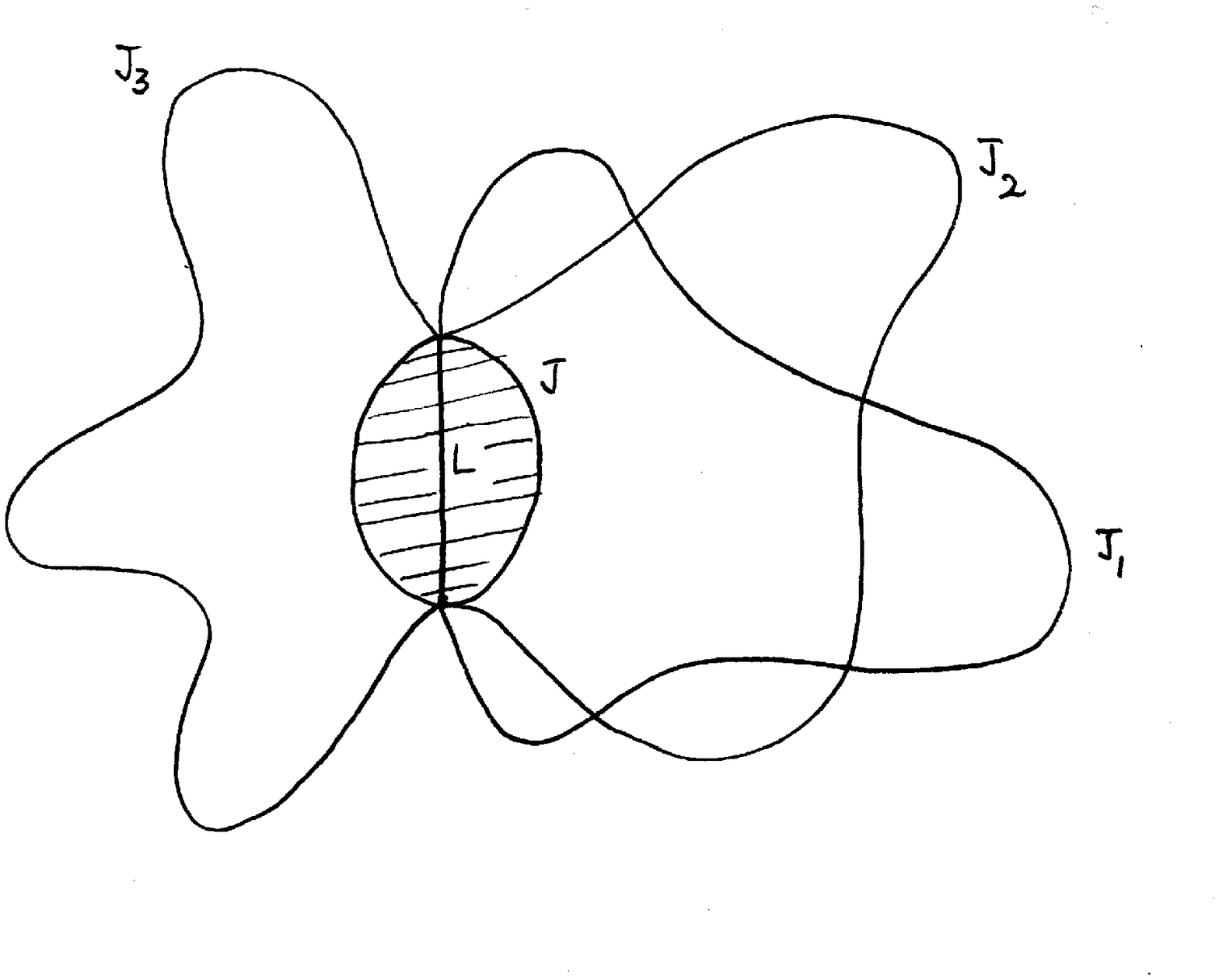,height=2.5in}
\end{center}
\caption{A Jordan neighborhood of a common arc $L$.} \label{f2}
\end{figure}

\noindent {\sc Proof}:  Let $\varepsilon_{i} = \sup \{ d(x,L); \ x \in
J_{i} \backslash L \}$.  Since $L$ is a closed arc, $J_{i} \sm \ov{L}
\neq \emptyset$ and thus $\varepsilon_{i} > 0$.  Let $A = L, \ B
=\overline{ \left( \bigcup J_{i} \right) \backslash L } \
= \ \bigcup \overline{ J_{i} \backslash L}$ and
$\varepsilon = \frac{1}{2} \min \varepsilon_i$.  Then $A \cap B = \{
\mbox{\rm endpoints of } A \}$ and $B \backslash (A \cap B) =
\bigcup (J_{i} \backslash L)$, every component of
which has a point at distance greater than $\varepsilon$ from $A$.  We
can then apply Proposition~\ref{6} in
order to find the desired Jordan curve $J$ (see figure~\ref{f2}.)  $\Box$

\begin{cor} \label{9}
With the notation of Corollary~\ref{8}, $J \cap L=$
\{endpoints of $L \}$ and thus $L$ is a cross-cut in $I$.
\end{cor}

\noindent{\sc Proof}:  Since $\stackrel{\circ}{L} \subset I, \ L \subset
\overline {I} = I \cup J$ so that \{endpoints of $L\} \subset I \cup J$.
On the other hand both endpoints of $L$ are accumulation points of each
$J_{i} \backslash L$ so that \{endpoints of $L$\} $\subset \left( \ov{
\bigcup  J_{i} } \right) \backslash L  \subset
\overline{O} =
O \cup J$.  Therefore \{endpoints of $L$\} $\subset J$. $\Box$
\bigskip

\begin{cor} \label{7}
Let $J$ be a Jordan curve, and $L \subset J$ a closed arc.  Then for any 
$\var >
0$ there exists an open cross-cut $\al \subset I \cap V_{\var}(L)$ with the
same endpoints as $L$. $\Box$
\end{cor}

\begin{prop} \label{10} The
closed disks $D_{1}, \ldots, D_n$ lie on the same side of a closed arc
$L$ if and only if there exists an open arc $\alpha \subset
\dis{\bigcap^{n}_{i=1}}I_i$ with the same endpoints as $L$.  As a
consequence, if $U$ is the Jordan domain bounded by $\alpha \cup L, \ U
\subset \dis{\bigcap^{n}_{i=1}}I_i$.
\end{prop}

\noindent{\sc Proof}:  If there exists such an arc, and $U$ is the Jordan
domain bounded by $\alpha \cup L$, by Theorem~\ref{1}, $\alpha \cup L
\subset \overline{U} \subset \left( \overline{ \bigcap
I_{i} } \right)$. Therefore $D_{1}, \ldots , D_{n} |_{L}$.

If $D_{1}, \ldots , D_{n}|_{L}$,
then $L \subset  \bigcap J_{i}$ and we can use Proposition~\ref{8}
 to find a Jordan curve $J$ satisfying the conclusions of that
proposition.  By Corollary~\ref{9}
and Corollary~
\ref{3}, $L$ separates $I$ into two Jordan
domains $U$ and $V$.  Notice that since $U \cup V = I
\sm L,  \ U \cup V$ does not intersect $L$ or
$\left(\bigcup J_{i} \right) \sm L$, that is, $U \cup
V \subset {\mc{C}} \left(\bigcup J_{i} \right)$.

Since $\stackrel{\circ}{L} \subset \bigcap I_{i}\, , \ \stackrel{\circ}{L}
\subset I$ and $L \cap ( \bigcap I_{i}) = \emptyset, \ (\bigcap I_{i}) 
\cap (U \cup V) \neq \es$.  Assume $U \cap ( \bigcap I_{i})
= \emptyset$.  Since $U \subset {\mc{C}} ( \bigcup J_{i}
)$ and $U$ is connected, $U \subset ( \bigcap  I _{i} )$.
Now, $U$ is bounded by $\al \cup L$, where $\alpha$ is one of the open
arcs into which the endpoints of $L$ separate $J$ (see figure 3.)  
Since $J \cap (\bigcup  J_{i})=$ \{endpoints of $L$\}, $\alpha \cap
(\bigcup J_{i}) = \es$ and since $\alpha
\subset \overline{U} \subset \overline{\bigcap
I_{i} }\, , \ \alpha \subset \bigcap I_{i}$.

Therefore $\alpha$ is the arc we were after.  $\Box$

\begin{figure}
\begin{center}~
\psfig{file=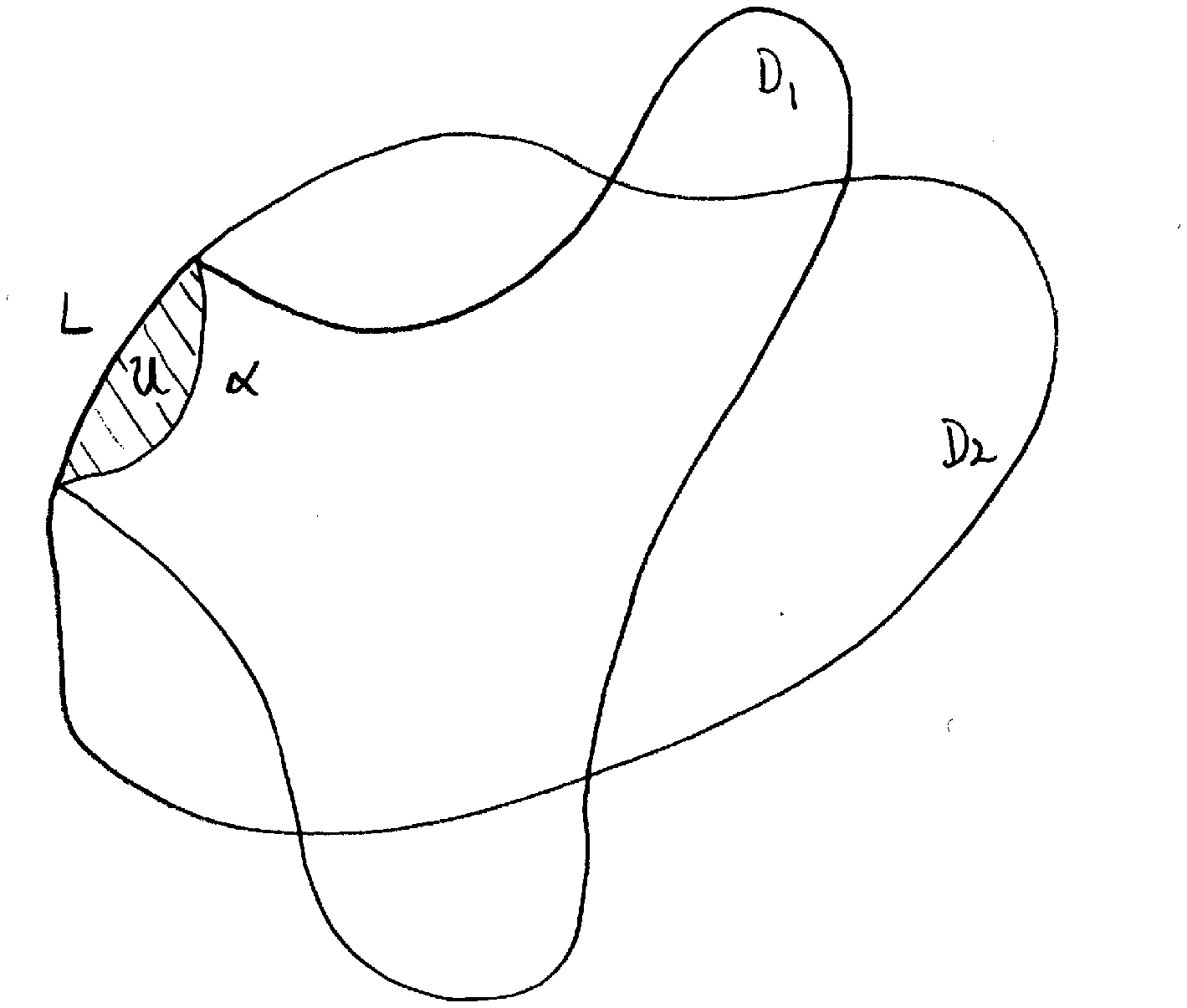,height=2.5in}
\end{center}
\caption{$D_{1}$ and $D_{2}$  are on the same side of $L$ and $\al$ is a 
cross-cut in both $D_1$ and $D_2$.}
\label{f3}
\end{figure}

\begin{cor}  [of the proof]  \label{11}  In Proposition~\ref{10},
$\al$ may be taken to lie in a $\var$-neighborhood of $L$, for any
$\var$ chosen in advance.  $\Box$
\end{cor}

\noi{\sc Remark}:  The arc $\al$ of Proposition~\ref{10}
 and Corollary~\ref{11} is
clearly a cross-cut in each of the domains $I_{i}$ for each $i \in 
\un{n} $. 
\bigskip

\begin{prop} \label{12}
If $D_{1}, \ldots, D_{n}|_{L'}$ and $L$ is the connected
component of $\dis\bigcap^{n}_{i=1} J_{i}$ containing $L'$, then $D_{1},
\ldots, D_{n}|_{L}$.
\end{prop}

\noindent {\sc Proof}:  Let $J$ be a Jordan curve as in
Corollary~\ref{8}
and $U$ and $V$ the components
of $J \sm L$.  By Corollary~\ref{9},
 $U \cup V \subset {\mc{C}}( \bigcup  J_{i})$.
Since $\sr{\circ}{L'} \subset \stackrel{\circ}{L} \subset I$ and $D_{1},
\ldots ,
D_{n}|_{L'}$, by the same reasoning as in the proof of
Proposition~\ref{10},
$( \bigcap  I_{i}) \cap (U \cup
V ) \neq \emptyset$, say, $( \bigcap  I_{i} ) \cap U
\neq \emptyset$.  Since $U \subset {\mathcal{C}} ( \bigcup 
J_{i}), \ U \subset  \bigcap  I_{i}$.  Thus, if $\partial U
= L \cup \alpha$, $\alpha$ satisfies the conditions of
Proposition~\ref{10},
which shows that $D_{1}, \dots, D_{n}|_{L}$ as
we wanted. $\Box$
\bigskip

\begin{prop} \label{13}  If $D_{1},
D_{2}|_{L}, \, D_{2}, \, D_{3}|_{L'}$ and $L'' \subset L \cap L'$ then
$D_{1}, D_{2}, D_{3}|_{L''}$.
\end{prop}

\noindent {\sc Proof}:  The proof is similar to the previous ones and is
left to the reader.  $\Box$
\bigskip

\begin{prop} \label{14}
Let $J_{0},
J_{1}, \ldots, J_{n}$ be Jordan curves, $L \subset J_0$ an open arc and
for $i \in \un{n}, \ L \cap \overline{I_{0} \cap I_{i} }= \emptyset$.
Then given $\varepsilon > 0$ there exists an open cross-cut $\alpha$ in
$I_0$ joining the endpoints of $L$ such that $\alpha \subset
V_{\varepsilon}(L)$ and if $U$ is the Jordan domain bounded by $\alpha
\cup L$, then $(U \cup \alpha) \cap D_{i} = \emptyset$ for each
$i \in \un{n}$.
\end{prop}

\noindent {\sc Proof}:  Consider the set $B=(J_{0} \backslash L) \cup
\left[ \overline{ I_{0}  \cap  ( \bigcup D_{i} ) } \right] $.
 $B$ is clearly closed, since $L$ is an open arc, and we claim
that $B \cap L = \emptyset$.  Since $I_0$ is open, it is an exercise to
show that $\overline{ I_{0} \cap I_{i} }= \overline{ I_{0} \cap D_{i}
}$.  Thus our assumption that $L_{0} \cap \overline{I_{0} \cap I_{i} } =
\emptyset$ is equivalent to $L \cap \overline{ I_{0} \cap D_{i} } =
\emptyset $ for each $i \in \underline{n}$.  Since $\overline{ I_{0} \cap
 \bigcup D_{i} } =
\bigcup \overline{ I_{0} \cap D_{i} }, \ L \cap
( \overline{ I_{0} \cap D_{i} } ) = \emptyset$ and clearly $L \cap (J_{0}
\backslash L ) = \emptyset$, so that $B \cap L = \emptyset$.

Now let $C$ be a component of $\overline{I_{0} \cap D_{i} }$ for some $i 
\in \un{n}$
and assume $C \cap (J_{0} \backslash L)= \emptyset$.  Since $C \cap L =
\es, \ C \cap J_{0} = \es$ and it follows that $C \subset I_{0}$.  But
$D_i$ is
connected so that $C = D_i$.  This shows that if a component of $B$ is
not that which contains $J_{0} \sm L$, it must consist of the union of
one or more of the closed disks $D_i$.  From this it is not hard to see
that there exists $\var > 0$ such that every component of $B \sm \{
\mbox{\rm endpoints of } L \}$ contains a point at distance greater
than $\var$ from $L$.  Let $A = \ov{L}$ and apply
Proposition~\ref{6}
 to $A, B$ and $\var$ as above to find a Jordan curve
$J$ such that $A \sm (A \cap B) = \ov{L} \sm \{ \mbox{\rm endpoints of
}L\}=L \subset I, \ B \sm (A \cap B ) = B \sm \{\mbox{\rm endpoints of }L\}
\subset O$ and $J \cap (A \cup B ) \subset A \cap B = \{\mbox{\rm
endpoints of
}L\}$.  Since $L \subset I$ and $J_{0} \sm L \subset \ov{O}$,  $L$ is a
cross-cut in
$I$ and $I \sm L = U \cup V$, where $U$ and $V$ are disjoint
Jordan domains.  Since $I \cap (J_{0} \sm L)= \es, \ I \sm L = I \sm [ L
\cup
(J_{0} \sm L )]= I \sm J_{0}$ and it follows that $I \sm L = (I \cap I
_{0}) \cup (I \cap O_{0})$ so that either $U = I \cap I_{0}$ and
$V=I \cap O_{0}$ or vice versa.  Assume $U=I \cap I_{0}$ (see figure 
4.)  Then $U
\cap D_{i}= \es$ for every $i \in \un{n}$ since $U = I \cap I_{0}$ and $I
\cap \ov{I_{0} \cap \bigcup D_{i} }= \es$.  Also, if $\al$ is the arc of $J
\sm (A \cup B )$ for which $\p U = \al \cup L$, it is clear that $\al
\subset I
_{0}$ and since $\al \cap \ov{I _{0} \cap \bigcup D_{i} } = \es, \ \al \cap
\bigcup D_{i} = \es$.  Therefore, $\al$ is the arc we were after.  $\Box$
\bigskip

\begin{figure}
\begin{center}~
\psfig{file=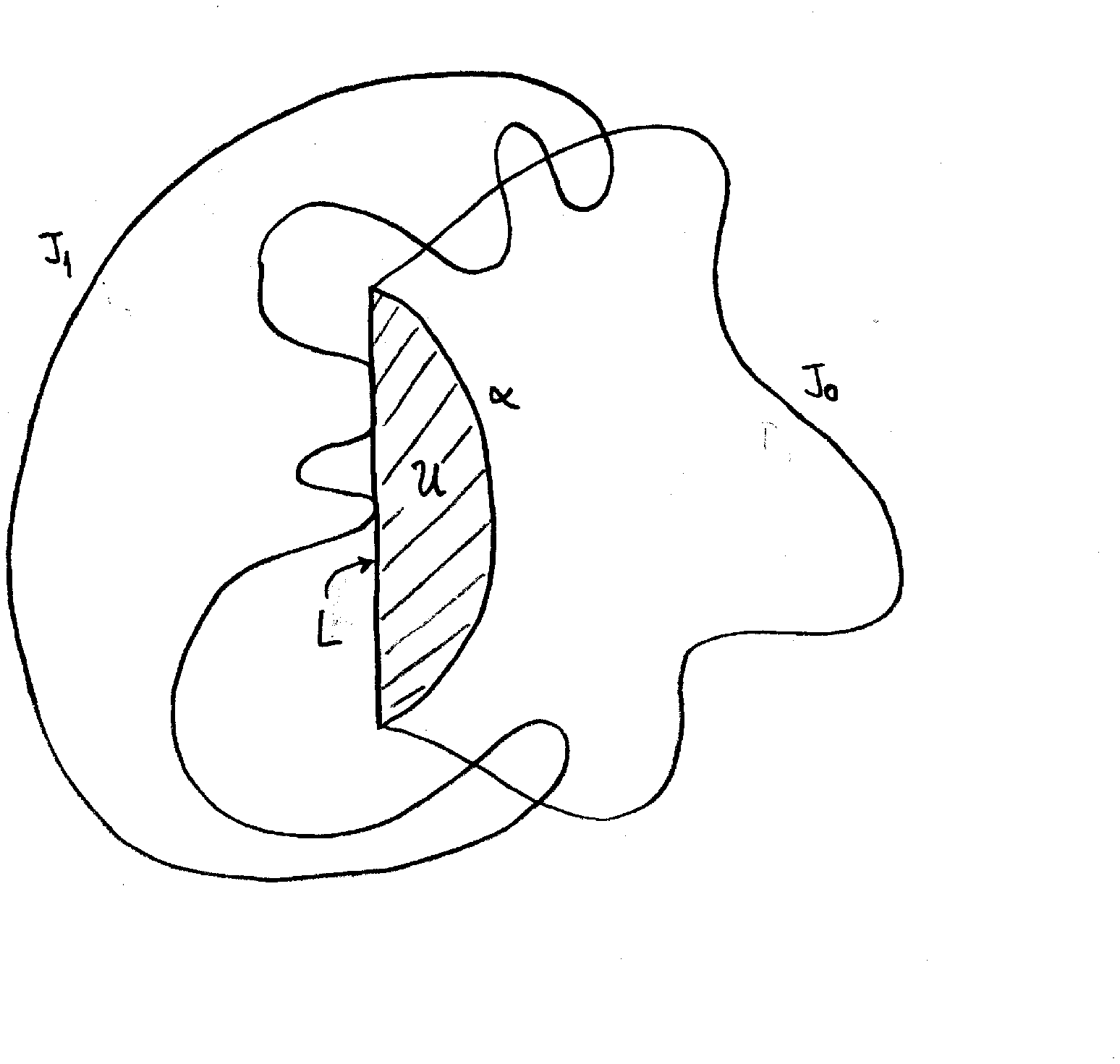,height=2.75in}
\end{center}
\caption{The curve $J_ 1$ only touches the arc $L$ from the outer domain 
determined by $J_0$.} \label{f4} 
\end{figure}

\begin{defn}
Let $A$ be a Jordan curve or an arc and $L, L'
\subset A$ closed arcs.  We say that $L$ and $L'$ are {\em unlinked}\/ if
either $L \subset L'$ or $L' \subset L$ or $L$ and $L'$ intersect at most at
endpoints.
\end{defn}

\noi {\sc Remark}:  Notice that saying that $L$ and $L'$ are unlinked in
a Jordan curve is more than the usual definition of their endpoints being
unlinked.
\bigskip

\begin{prop} \label{15}
Let $J$
be a Jordan curve and $L_{1}, \dots , L_{n} \subset J$ be pairwise unlinked
closed arcs.  Then for every $\var > 0$ there exist disjoint open
cross-cuts $\al_{i} \subset I \cap V _{\var} (L_{i})$ joining the
endpoints of $L_i$, for each $i \in \un{n}$.
\end{prop}

\noi {\sc Proof}:  We will use induction on the number $n$ of arcs.  For
$n=1$, the statement is true by Corollary~\ref{7}.
 Assume we have proven the statement for collections of arcs
with up to $n-1$ elements and $L_{1}, \ldots ,L_{n}$ are unlinked.  Use
Corollary~\ref{7} to find an open cross-cut
$\al_{1} \subset I \cap V_{\var} (L_{1})$ joining the endpoints of $L_1$.
Then for $i > 1$, since $L_{i}, L_{1}$ are unlinked, either $L_{i} \subset
L_{1}$ or $L_{i} \subset \ov{J \sm L_{1} }$.  Let $L_{i_{1}}, \ldots ,
L_{i_{k}}
\subset L_{1}$ and $L_{j_{1}}, \ldots , L_{j_{m}} \subset \ov{ J \sm L_{1}
}$.
These are collections of unlinked arcs with fewer than $n$ elements and
since $L_{i_{1}}, \ldots , L_{i_{k}} \subset L_{1} \cup \al_{1}$ and 
$L_{j_{1}},
\ldots , L_{j_m} \subset \ov{J \sm L_{1} } \cup \al_{1}$, by the inductive
hypothesis it is possible to find collections of cross-cuts $\al_{i_{1}},
\ldots , \al_{i_{k}}$ and $\al_{j_{1}}, \ldots , \al_{j_{m}}$
satisfying the conclusion of the proposition.  Clearly $\al_{1},
\al_{i_{1}}, \ldots , \al_{i_{k}}, \al_{j_{1}}, \ldots , \al_{j_{m}}$
is the desired collection for $L_{1}, \ldots, L_{n}$.  $\Box$
\bigskip

\begin{prop} \label{16}
Let
$J_{0}, \ldots , J_n$ be Jordan curves, and $L_{i} \subset J_{i} \cap
J_{0}, \ i \in \un{n}$,
closed arcs, pairwise unlinked in $J_0$, no two of which
are indentical.  Assume that $D_{0}, D_{i} |_{L_{i}}$ for $i \in
\un{n}$.  Then for each $\var > 0$ there exist disjoint open cross-cuts
$\al_{i} \subset I_{0}$ joining the endpoints of $L_i$ such that $\al_{i}
\subset V_{\var} (L_{i}) \cap I_{i}$ for $i \in \un{n}$.
\end{prop}

\noi {\sc Proof}:  The proof is by induction on the number $n$ of
curves.  If $n=1$, the statement is true by Corollary~\ref{7}.
Assume we have proven the statement for
collections with fewer than $n$ curves and $J_{i}, L_{i}, \ i \in
\un{n}$ satisfy the hypotheses above.  Among $L_{1}, \ldots , L_{n}$
choose all the ones which are not contained in any other (see figure 
5.)  We may assume without loss of generality that they are the first $k$ 
arcs
$L_{1}, \ldots , L_k$.  Since $L_{1}, \ldots , L_k$ are pairwise unlinked
and are not contained in one another, they are pairwise disjoint except
possibly at endpoints.  By Proposition~\ref{15}
there exist disjoint open cross-cuts $\gamma_{i} \subset I_{0} \cap
V_{\var}(L_{i})$
joining the endpoints of $L_i$ for $i \in \un{k}$.  Notice that since the
arcs $L_{1}, \ldots , L_k$ are disjoint except possibly at endpoints, the
interior $U_i$ of the disks bounded by $\gamma_{i} \cup L_{i}, \ i \in
\un{k}$ are pairwise disjoint.  Moreover by
Proposition~\ref{13}
the closed disk bounded by $\gamma_{i} \cup L_{i}$ is on the same side
of $L_i$ as $D_0$ (the disk bounded by $J_0$) for each $i \in \un{k}$.
>From Proposition \ref{10} it follows
that there exist
arcs $\al_{i} \subset U _{i} \cap I_{i}$ joining the endpoints of $L_i$.
Since $U_{i} \subset I_{0} \cap V _{\var} (L_{i})$ it is clear that $\al_i$
is a
cross-cut in $I_0$ and $\al_{i} \subset V_{\var}(L_{i})\cap I_{i}$.  Now,
the
remaining arcs are contained in $L_{1}, \ldots , L_k$, since we chose
all the arcs which were not contained in any other.  For each $L_{i}, \ i
\in \un{k}$, the arcs inside it form an unlinked collection with fewer
than $n$ elements satisfying the hypotheses of the proposition.
Therefore by the inductive assumption we are done.  $\Box$
\bigskip

\begin{figure}
\begin{center}~
\psfig{file=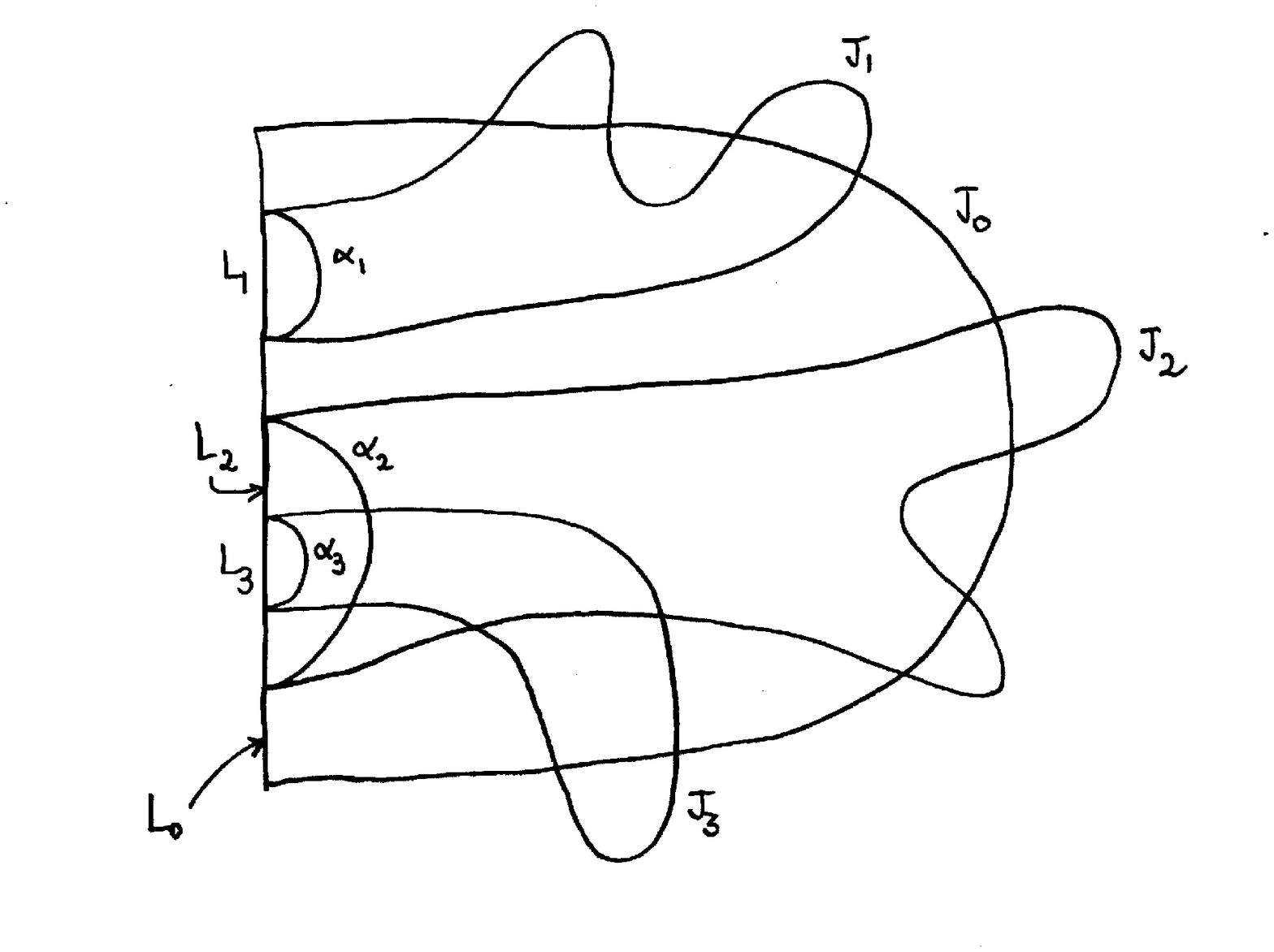,height=2.5in}
\end{center}
\caption{Disjoint cross-cuts joining the endpoints of unlinked arcs.} 
\label{f5} 
\end{figure}

%% file: cedisks.tex
\section{(\boldmath{$c,e$})-Disks and Pruning Fronts}
\label{cedisks}

We will now give a preliminary definition of what we call $(c,e)$-disks.
Later we will add a dyamical hypothesis which is not necessary at present.
\bigskip

\begin{defn}
A closed disk $D$ is called a $(c,e)$-{\em disk}\/
if there are closed arcs $C, E \subset \p D$ specified such that $\p D=C
\cup
E$ and $C$ and $E$ only intersect at endpoints.  In other words, for now,
a $(c,e)$-disk is just a {\em bigon}\/ with sides $C$ and $E$.  We call
the common endpoints of $C$ and $E$ the {\em vertices}\/ of $D$.
\end{defn}

\begin{defn}
\label{longer}
Let $D_{1}, D_{2}$ be $(c,e)$-disks such that
$I_{1} \cap I_{2} \neq \es$.  We say $D_1$ is {\em e-longer}\/ or simply
{\em longer}\/ than $D_2$, denoted $D_{1} \succ D_2$, if (i), (ii) and (iii)
hold (see figure 6):

\begin{description}
 \item [(i)]  $C_{1} \cap I_{2} = \es$ and $E_{2} \cap I_{1} = \es$;

 \item [(ii)] if $C_{1} \cap C_{2} \neq \es$ then $C_{1} \cup C_{2}$ is an
arc and if $E_{1} \cap E _{2} \neq \es$ then $E_{1} \cup E_{2}$ is an arc;

  \item [(iii)] if ${\sr{\circ}{C_{1}} } \cap \ov{ I_{1} \cap I_{2} } \neq
\es$ then $C_{1} \subset C_2$ and if ${\sr{\circ} {E_{2}} }\cap \ov{ I_{1}
\cap I_{2} } \neq \es$ then $E_{2} \subset E_1$.

\end{description}
\end{defn}

\begin{figure}
\begin{center}~
\psfig{file=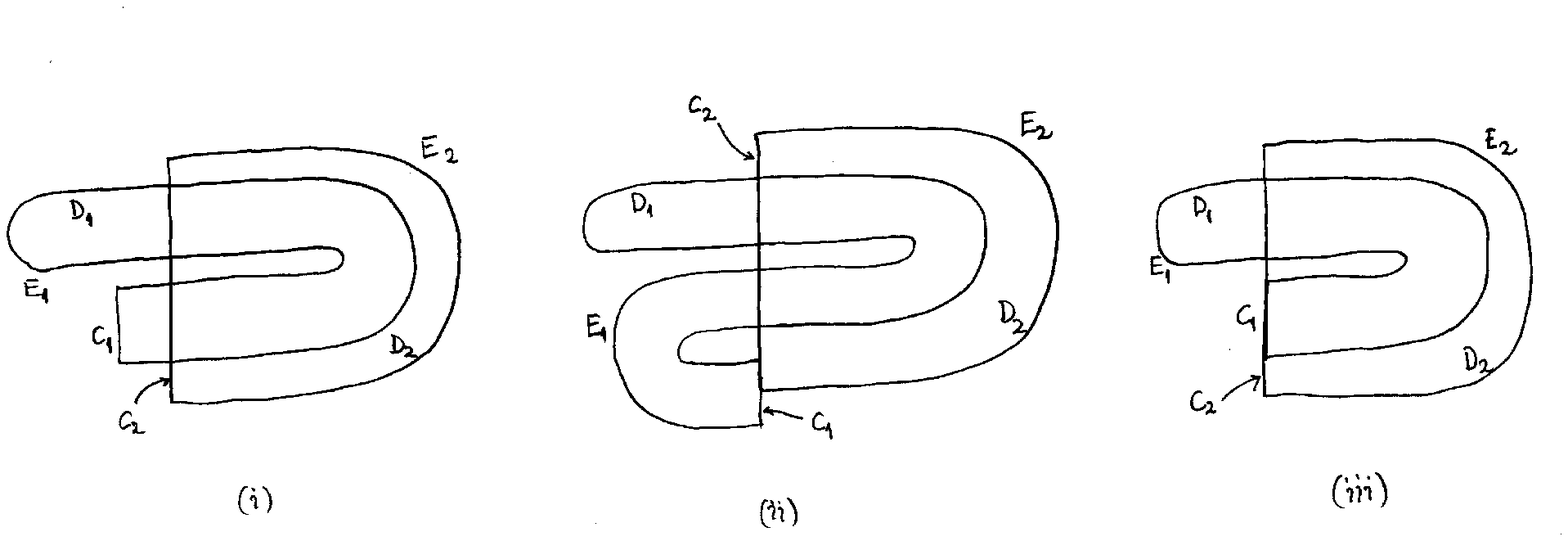,height=2.5in}
\end{center}
\caption{The relation $\succ$.} \label{f6}
\end{figure}

\noi{\sc Notation}:  Let $D$ be a $(c,e)$-disk and $\al$ a cross-cut
joining the vertices of $D$.  We have seen that $\al$ separates the
interior $I$ of $D$ into two Jordan domains whose boundaries are $C \cup
\al$ and $E \cup \al$.  We denote them by $I^{c}(\al)$ and $I^{e}(\al)$,
respectively, and their closures by $D^{c}(\al), \ D^{e}(\al)$ (see figure 7.)  
Moreover, when the disks are indexed and so are the cross-cuts we will
only use the index inside the parentheses so that $D^{c}(\al_{i})$
will denote the disk bounded by $C_{i} \cup \al_{i}$.
\bigskip

\begin{figure}
\begin{center}~
\psfig{file=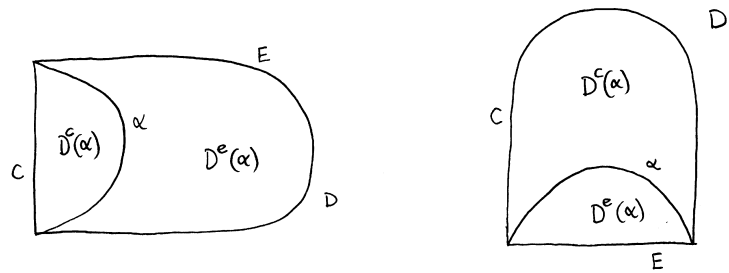,height=1.75in}
\end{center}
\caption{Cut $(c,e)$-disks.} \label{f7}
\end{figure}

\noi{\sc Conventions}:  If $D_{1}, \ldots , D_{L}$ is a collection of
$(c,e)$-disks, we say they are {\em related}\/ by $\succ$ if for any $i,
\ j
\in \un{L}$ either $I_{i} \cap I_{j} = \es$ or $D_{i} \succ D_{j}$ or
$D_{j} \succ
D_{i}$.  When nothing is mentioned about a colection of $(c,e)$-disks it
is assumed they are related by $\succ$.
Cross-cuts in $(c,e)$-disks, when nothing is mentioned to the contrary,
are assumed to be open and to join the vertices of the disk wherein they lie.
\bigskip

The following propositions are easy consequences of what we have
developed so far and we omit the proofs.
\bigskip

\begin{prop}  \label{17}
If $D$
is a $(c,e)$-disk and $\al, \beta, \gamma \subset D$ are open cross-cuts
joining vertices such that $\beta \subset I^{c} (\al)$ and $\gamma \subset
I^{e}(\al)$ then $I^{c}(\beta) \subset I^{c}(\al) \subset I^{c}(\gamma)$ and
$I^{e}(\beta) \supset I^{e}(\al) \supset I^{e}(\gamma)$. $\Box$
\end{prop}

\begin{prop} \label{18}
Let $D_{1}$ and
$D_2$ be $(c,e)$-disks and $D_{1} \succ D_{2}$.  Then

\begin{description}
\item[(i)]  if ${\sr{\circ}{C_{1}} }\cap \ov{ I_{1} \cap I_{2} } \neq \es$,
then $D_{1}, D_{2}|_{C_{1}}$ and

\item[(ii)] if ${\sr{\circ}{E_{2}} } \cap \ov{ I_{1} \cap I_{2} } \neq \es$,
then $D_{1}, D_{2}|_{E_{2} }$.  $\Box$
\end{description}
\end{prop}

\begin{defn}
A collection of pairs $\{(D_{i}, \beta_{i})\}^{L}_{i=1}$,
where $\{D_{i} \}^{L} _{i=1}$ is a collection of
$(c,e)$-disks related by $\succ$ and $\{ \beta_{i} \subset D_{i} 
\}^{L}_{i=1}$ 
is a collection of open cross-cuts joining vertices, will be called a {\em
cut collection}.
\end{defn}

\begin{prop} \label{19}
Let $D_{1}$ and $D_{2}$ be $(c,e)$-disks and $D_{1} \succ D_{2}$.  If
$C_{1}=C_{2}$ or $E_{1}=E_{2}$ then $D_{1}=D_2$.
\end{prop}

\noi{\sc Proof}:  Assume $C_{1}=C_{2}$.  Then the endpoints of $E_{1}$ and
$E_{2}$ coincide (since they are the same as those of $C_{1}$ and $C_2$)
and, by (ii) in the definition of $\succ, \ E_{1} \cup E_2$ is an arc.  But
this can only happen if $E_{1}=E_2$.  $\Box$
\bigskip

\begin{prop} \label{20}
If $D_{1},
D_{2}$ are $(c,e)$-disks and $D_{1} \succ D_2$ and $D_{2} \succ D_1$ then
$D_{1}=D_2$.
\end{prop}

\noi{\sc Proof}:  The proof is easy and is left to the reader.  $\Box$
\bigskip

\begin{prop} \label{21}
Let
$\{(D_{i}, \beta_{i} ) \}^{L}_{i=0}$ be a cut collection and $\var$ a
positive number.
 If $D_{0} \not\prec D_i$ (i.e., either $I_{0} \cap
I_{i} = \es$ or $D_{0} \succ D_{i}$ and $D_{0} \neq D_{i}$) for every $i
\in \un{L}$ then there exists an open cross-cut $\al_{0} \subset
I^{c}(\beta_{0}) \cap V_{\var} (C_{0})$ joining vertices such that for
each $i \in \un{L}$ either (i) or (ii) holds:

\begin{description}
\item[(i)] if ${\sr {\circ}{C_{0}} } \cap \ov{ I_{0} \cap I_{i} } \neq \es$
then $[ I^{c}( \al_{0}) \cup \al_{0} ] \subset I^{c}(\beta_{i})$;

\item[(ii)]  otherwise $[I^{c}(\al_{0}) \cup \al_{0} ] \cap D_{i} = \es$.

\end{description}

If, on the other hand, $D_{0} \not\succ D_i$ for every $i \in \un{L}$,
then there exists an open cross-cut $\al_{0} \subset I^{e}(\beta_{0}) \cap
V_{\var} (E_{0})$ such that for each $i \in \un{L}$ either (iii) or (iv)
holds:

\begin{description}
\item[(iii)]  if ${\sr {\circ}{E_{0}} } \cap \ov{ I_{0} \cap I_{i} }\neq
\es$ then $[ I^{e} (\al_{0}) \cup \al_{0} ] \subset I^{e}(\beta_{i})$;

\item[(iv)]  otherwise $[I^{e}(\al_{0}) \cup \al_{0}] \cap D_{i}= \es$.
\end{description}
\end{prop}

\noi{\sc Proof}:
We will prove (i) and (ii), the proof of (iii) and (iv) being analogous.
Divide the disks $D_i$ into two groups:  (i) those for which $C_{0}
\cap \ov{I_{0} \cap I_{i} } \neq \es$ and (ii) those for which $C_{0} \cap
\ov{I_{0} \cap I_{i} } = \es$.  If $D_i$ is in group (i), $I_{0} \cap
I_{i} \neq \es$, so that by our assumption $D_{0} \succ D_i$ and, by
Proposition~\ref{18}, $D_{0}, D_{i}|_{C_{0}}$.  Clearly
$D^{c}(\beta_{0}), D_{0}|_{C_{0} }$ and $D^{c}(\beta_{i}),
D_{i}|_{C_{i}}$ and, since $C_{0} \subset C_i$, we see that
$D^{c}(\beta_{0}), D_{0} |_{C_{0}}$, $D_{0}, D_{i}|_{C_{0}}$ and
$D_{i},D^{c}(\beta_{i})|_{C_{i}}$, by Proposition~
\ref{13}, imply that $D^{c}(\beta_{0}),
D^{c}(\beta_{i})|_{C_{0}}$ for every $D_i$ in group (i).  It now follows
>from Proposition \ref{10} and
Corollary~\ref{11}
that there exists an open
cross-cut $\al \subset I^{c}(\beta_{0}) \cap V_{\var} (C_{0})$ such that
$[I^{c} (\al ) \cup \al ] \subset I^{c} (\beta_{i})$.

On the other hand, for the disks $D_j$ in group (ii), $C_{0} \cap \ov{ I_{0}
\cap I_{j} } = \es$ and since $I^{c} (\al ) \subset I_0$, it is also the
case that $C_{0} \cap \ov{I^{c} (\al ) \cap I_{j} } = \es$.  Thus, by
Proposition~\ref{14} there exists an open cross-cut $\al_0$
in $I^{c}(\al)$ such that for every $D_j$ in group (ii), $[I^{c} (\al _{0}
) \cup \al_{0} ] \cap D_{j} = \es$.  It is clear that such $\al_0$ also
satisfies $[I^{c} (\al_{0}) \cup \al_{0} ] \subset I^{c} (\beta_{i})$ for
every $D_i$ in group (i) (see figure 8.)

In the event that all the disks belong to one or the other of the groups,
the modifications necessary in the above proof are minor and are left to
the reader.  $\Box$
\bigskip

\begin{figure}
\begin{center}~
\psfig{file=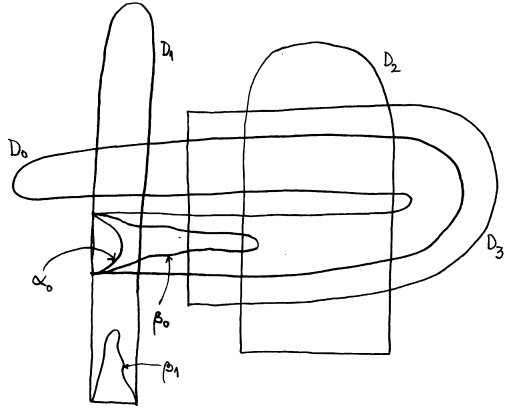,height=2.75in}
\end{center}
\caption{A cut collection where $D_{0} \succ D_{i}$ for every $i \in 
\mbox{\un{\it L}}$.}\label{f8}
\end{figure}

\begin{defn}
Let $\{(D_{i}, \beta_{i}) \}^{L}_{i=1}$ be a
cut collection, $S \subset \un{L}$ and $\var > 0$.  The collection $\{
\al_{i} \}_{i \in S}$ of disjoint open cross-cuts is said to be a
$(\var,c)$-{\em collection compatible with}\/ $\{ ( D_{i}, \beta_{i} ) \}
^{L}_{i=1}$ (see figure 9) if $\al_{i} \subset I^{c}(\beta_{i} ) \cap 
V_{\var} (C_{i})$ and
for every $i \in S$ and $j \in \un{L}$ such that $D_{i} \not\prec D_j$
either (i) or (ii) holds:

\begin{description}
\item[(i)] if ${\sr {\circ}{C_{i}} } \cap \ov{ I_{i} \cap I_{j} } \neq \es$
then $[I^{c} (\al_{i} ) \cup \al _{i} ] \subset I^{c} (\beta _{j} )$;

\item[(ii)]  otherwise $[I^{c} (\al_{i} ) \cup \al_{i}] \cap D_{j} = \es$.
\end{description}

The collection $\{\al_{i}\}_{i \in S}$ is called a $(\var,e)$-{\em
collection compatible with}\/ $\{(D_{i},$ $\beta _{i}) \}^{L}_{i=1}$ if
$\al_{i} \subset I^{e} (\beta_{i} ) \cap V_{\var} (E_{i})$ and for every $i
\in S$ and $j \in \un{L}$ such that $D_{i} \not\succ D_{j}$ either (iii)
or (iv) holds:

\begin{description}
\item[(iii)] if ${\sr{\circ}{E_{i}} } \cap \ov{ I_{i} \cap I_{j} } \neq \es$
then $[ I^{e} (\al_{i}) \cup \al_{i} ] \subset I^{e} (\beta _{j} )$;

\item[(iv)] otherwise $[I^{e} (\al_{i}) \cup \al_{i} ] \cap D_{j} = \es$.
\end{description}
\end{defn}

\begin{figure}
\begin{center}~
\psfig{file=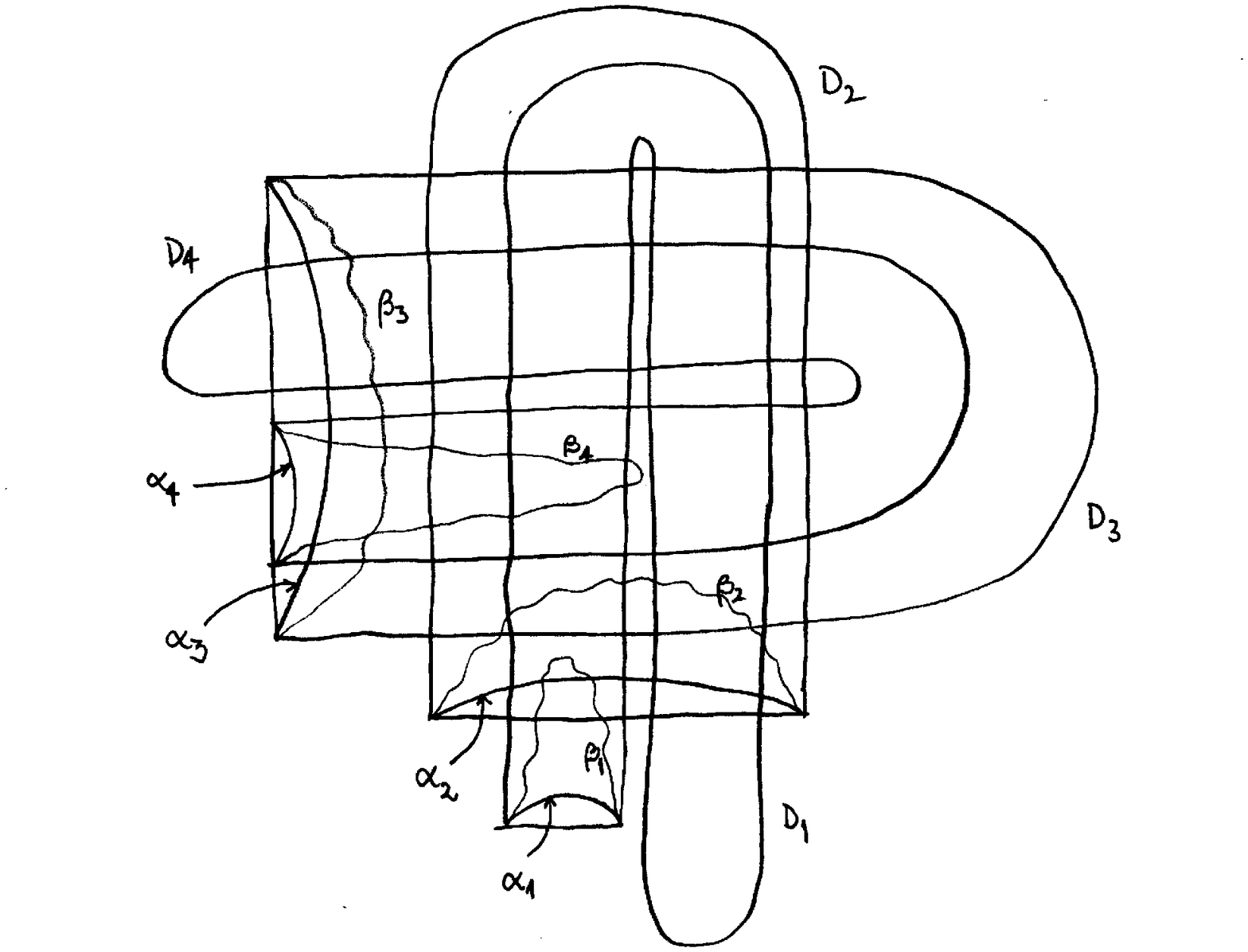,height=3in}
\end{center}
\caption{$\{\al_{1}, \al_{2} \}$ is a ($\var, e$)-collection and 
$\{\al_{3}, \al_{4} \}$ is a ($\var, c$)-collection, both compatible with 
$\{ (D_{i}, \beta_{i} ) \}^{4}_{i=1}$.} \label{f9}
\end{figure}

\noi {\sc Remarks}:  Notice that if $\{ \al_{i}\}_{i \in S}$ is a
$(\var,c)$-collection compatible with $\{ ( D_{i}, \beta_{i} ) \}$ and $\{
\gamma_{i} \}_{i \in S}$ is a collection of open cross-cuts joining the
vertices of $D_i$ and such that $\gamma _{i} \subset I^{c} (\al_{i} ),\ \{
\gamma_{i} \}_{i \in S}$ is also a $( \var,c)$-collection compatible with
$\{(
D_{i}, \beta_{i} ) \}$.  If moreover $\gamma_{i} \subset V_{\var '}(C_{i})$
then $\{ \gamma _{i} \}_{i \in S}$ is a $(\var ', c)$-collection.  The
analogous statement holds true for $(\var,e)$-collections.
\bigskip

\noi{\sc Warning}:  As the reader may have already noticed, statements
about $c$-``things'' and $e$-``things'' are ``dual'' to one another and
most proofs are totally analogous in both cases.  We will henceforward,
whenever there is nothing essentially different between the two, present
only the ``$c$-proof'' without further comments.
\bigskip

\begin{prop} \label{22}
Let $\{(D_{i}, \beta_{i}) \}^{L}_{i=1}$ be a
cut collection, $\{ \al_{i} \}_{i \in S}$ a $(\var,c)$-collec\-tion and $\{
\al_{i}' \}_{i \in S}$ a $(\var, e)$-collection both compatible with
$\{(D_{i}, \beta_{i} ) \}^{L}_{i=1}$.  If $i, j \in S$ are such that $ 
D_{i} \succ D_{j}$ and $ D_{i} \neq D_j$ then:

\begin{description}
\item[(i)] ${\sr{\circ}{C_{i}} } \cap \ov{ I_{i} \cap I_{j} } \neq \es$
implies $[I^{c} (\al_{i}) \cup \al_{i} ] \subset I^{c}(\al_{j} )$ and

\item[(ii)] ${\sr{\circ} {E_{j}} } \cap \ov{ I_{i} \cap I_{j} } \neq \es$
imples $[I^{e}(\al_{j}') \cup \al_{j}' ] \subset I^{e} (\al_{i}')$.
\end{description}
\end{prop}

\noi {\sc Proof}:  From the definition of $\succ$ and
Proposition~\ref{18}
it follows, under the hypotheses above, that $C_{i} \subset
C_{j}$ and $D_{i}, D_{j} |_{C_{i}}$ and from the definition of
$(\var,c)$-collection, that $[I^{c}(\al_{i} ) \cup \al_{i} ] \subset I^{c}
(\beta _{i} )$.  Therefore both $\al_i$ and $\al_j$ are open cross-cuts
in $I^{c}(\beta_{j})$.  Since they are assumed to be disjoint (by
definition), $\al_i$ joins the endpoints of $C_i, \ \al_j$ those of
$C_j$  and $C_{i} \subset C_j$, it must be the case that $[I^{c} (\al_{i})
\cup \al_{i} ] \subset I^{c} (\al _{j})$, as we wanted.  $\Box$
\bigskip

\begin{prop} \label{23}
Let $\{\al_{i} \}_{i \in S}$ be a ($\var ,c$)-collection
compatible with the cut collection $\{ ( D_{i}, \beta_{i} ) \}
^{L}_{i=1}$ and $\{\beta_{i}' \subset D_{i} \}^{L} _{i=1}$ a collection of
cross-cuts such that $\beta_{i}' \subset D^{e}(\beta_{i})$ for each $i \in
\un{L}$.  Then $\{\al_{i} \}_{i \in S}$ is also compatible with
$\{(D_{i}, \beta_{i}') \}^{L}_{i=1}$.  If above we change ($\var,c$)-
to ($\var,e$)- and $D^{e}(\beta_{i})$ to $D^{c}(\beta_{i})$ the resulting
statement is true.
\end{prop}

\noi {\sc Proof}:  Since the collection of $(c,e)$-disks remains
unchanged all there is to check is that if $i \in S$ and $j \in \un{L}$
are such that $D_{i} \succ D_{j}, \ {\sr{\circ}{C_{i}} } \cap \ov{ I_{i}
\cap I_{j} } \neq \es$ implies $[I^{c} (\al_{i}) \cup \al_{i} ] \subset 
I^{c} (\beta _{j} ' )$.  But by the ``closed'' version of
Proposition~\ref{17},
$\beta_{j}' \subset D^{e} (\beta_{j} )$ implies
that $I^{c} (\beta_{j}' ) \supset I^{c} (\beta _{j})$.  The result now
follows.  $\Box$
\bigskip

\begin{cor} \label{24}
Let $\{\al_{i} \}_{i \in S}$ and $\{
\al _{i}'\}_{i \in S'}$ be a $(\var,c)$- and a $(\var ',e)$-collection
respectively, both compatible with the cut collection $\{ (D_{i},
\beta_{i} ) \}^{L}_{i=1}.$  Then $\{\al_{i}\}_{i \in S}$ is a
$(\var,c)$-collection compatible with the cut collection $$\{(D_{i},
\beta_{i}), \ i \in \un{L} \sm S' \} \cup \{ ( D_{i}, \al_{i}'); \ i \in
S'\}$$ and $\{ \al_{i}'\}_{i \in S'}$ is a $(\var,e)$-collection
compatible with $$\{(D_{i}, \beta_{i}); \ i \in \un{L} \sm S \} \cup \{ (
D_{i}, \al_{i} ); \ i \in S \}. \ \ \Box$$
\end{cor}

\begin{prop} \label{25}
Under the hypotheses of Corollary~
\ref{24},
if $i \in S$ and $j \in S'$ are such that $D_{i} \not\prec D_{j}$, then
$\al_{i} \cap \al_{j} = \es$.
\end{prop}

\noi{\sc Proof}:   The proof is easy and is left to the reader.  $\Box$
\bigskip

We still have to show that $(\var , c)$- and $(\var ,e )$-collections
exist.  In the proof we will use the definition and the proposition below.
\bigskip

\begin{defn}
Let $\{D_{i}\}$ be a collection of $(c,e)$-disks
related by $\succ$.  We say $D_i$ and $D_j$ are $c$-equivalent, and write
$D_{i} \sim_{c} D_j$, if there exits $D_k$ in the collection such that
$C_{i}, C_{j} \subset C_k$ and $D_{i}, D_{k} |_{C_{i}}$ and $D_{j}, D_{k}
|_{C_{j}}$.  We define $e$-equivalence analogously by changing $c$-sides
to $e$-sides above, and denote it by $\sim_{e}$ (see figure 10.)
\end{defn}

\begin{figure}
\begin{center}~
\psfig{file=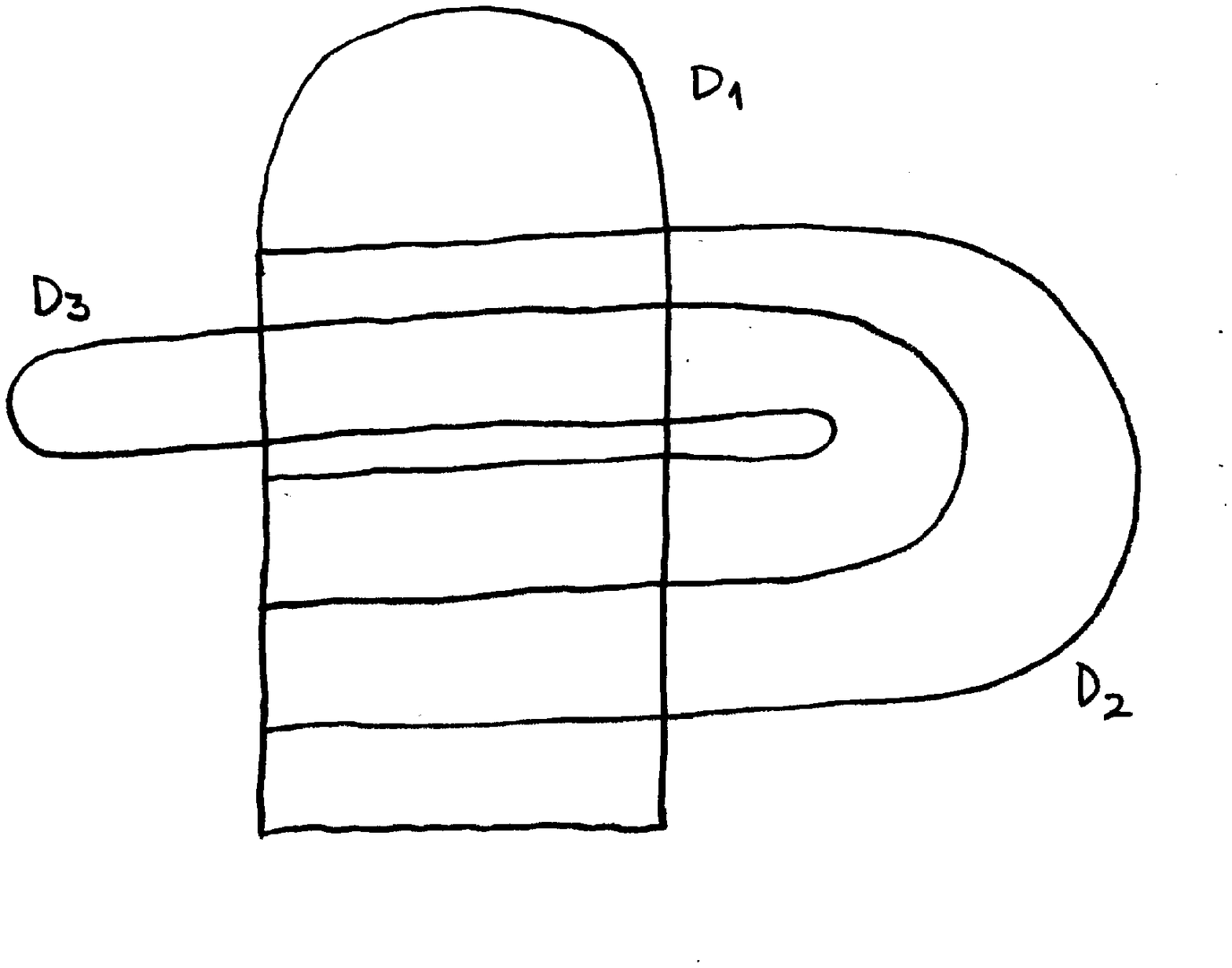,height=2in}
\end{center}
\caption{An equivalence class for $\sim_c$.  $D_1$ is the distinguished 
representative.} \label{f10}
\end{figure}

\noi{\sc Remark}:  Notice that by this definition, $D_{i} \sim_{c} D_j$
if $C_{i} \subset C_{j}$ and $D_{i}, D_{j} |_{C_{j}}$ or vice versa and
analogously for $\sim_e$.
\bigskip

\begin{prop} \label{26}
The
relations $\sim_c$ and $\sim_e$ defined above are equivalence relations.
If the collection $\{D_{i} \}^{L}_{i=1}$ is finite, each equivalence
class for $\sim_{c} \ (\sim_{e})$ has a distinguished representative
whose c-(e-)side contains the c-(e-)sides of all other disks in
its c-(e-)equivalence class.  Moreover, in each c-(e-)equivalence 
class the c-(e-)sides are unlinked in the c-side of
its distinguished representative.
\end{prop}

\noi {\sc Proof}:
 That $\sim_c$ is reflexive and symmetric
is clear.  In order to prove transitivity, assume $D_{i}
\sim_{c} D_{j}$ and $D_{j} \sim_{c} D_{k}$.  This means there exist
$D_{l}, \ D_{m}$ in the collection such that $C_{i}, C_{j} \subset C_{l}$
and $D_{i}, D_{l}|_{C_{i}}, \ D_{j}, D_{l}|_{C_{j}}$ and $C_{j}, C_{k}
\subset
C_{m}$ and $D_{j}, D_{m}|_{C_{j}}, \ D_{k}, D_{m}|_{C_{k}}$.  It follows
>from Proposition~\ref{13} that $D_{l},
D_{m}|_{C_{j}}$
and thus either $D_{l} \succ D_m$ or $D_{m} \succ D_{l}$.  We may assume
$D_{l} \succ D_m$, the other case being analogous.  Then, since $C_{j}
\subset C_l$ and $D_{l}, D_{m}|_{C_{j}}, \ {\sr {\circ}{C_{l}} } \cap \ov{
I_{l} \cap I_{m} } \neq \es$ and from the definition of $\succ$ and
Proposition~\ref{18} we can conclude that $C_{l}
\subset C_{m}$
and $D_{l}, D_{m}|_{C_{l}}$.  From this we see that $C_{i} \subset C_{m}$
and $D_{i}, D_{m}|_{C_{i}}$, which shows that $D_{i} \sim_{c} D_{k}$.
\bigskip


Consider now one $c$-equivalence class and let $D_i$ be an element in it
whose $c$-side is not strictly contained in the $c$-side of any other
disk in the same class.  If $D_{j}\sim_{c}D_{i}$ then it must be the case
that $C_{j} \subset C_i$ for otherwise there would exist $D_k$ in the
collection for which $C_{i}, C_{j} \subset C_{k}$ and $D_{i},
D_{k}|_{C_{i}}$ and $D_{j}, D_{k}|_{C_{j}}$.  But $D_{k} \sim_{c} D_{i}$
(see the remark just after the definition of $c$-equivalence) and if
$C_{j} \not\subset C_{i}, \ C_{k}$ contains $C_{i}$ strictly which is
contrary to our assumption.  This shows that for every $D_{j}$ such that
$D_{j} \sim_{c} D_{i}$ we have $C_{j} \subset C_{i}$ and $D_{i}, D_{j}
|_{C_{j}}$.  In order to see that the $c$-sides of disks in the
$c$-equivalence class of $D_{i}$ are unlinked in $C_{i}$ assume $D_{j}
\sim_{c} D_{k} \sim_{c} D_{i}$ and that $C_{j} \cap C_{k} \supset C$,
where $C$ is a closed arc.  Since $D_{i},
D_{j}|_{C_{j}}$ and $D_{j}, D_{k}|_{C_{k}}$ by Proposition~\ref{13}, it
follows that $D_{j}, D_{k}|_{C}$.  Then $I_{j} \cap I_{k} \neq \es$ and we
must have $D_{j} \succ D_{k}$ or $D_{k} \succ D_{j}$ and by (iii) in the
definition of $\succ$, $C_{j} \subset C_{k}$ or $C_{k} \subset C_{j}$. $\Box$

\bigskip

\noi{\sc Standing Convention}:  If the lower index in an indexed union or
collection is larger than the upper one we will take the union or
collection to be empty, so that $\dis{\bigcup ^{n-1}_{-n+1} }f^{k}(P) =
\es$ when $n=0$.  Also, recall that a bar under a positive integer
denotes the set of all positive integers smaller than or equal to it:
$\un{L} = \{1,2, \ldots, L \}$.  If $L=0$ we take $\un{L}$ to be the empty
set as well.

\bigskip

We now go on to prove the existence of $(\var, c)$- and $(\var,
e)$-collections (see figure 11.)

\begin{prop} \label{27}
Let $\{(D_{i}(k), \beta_{i}(k) ); \ k= -1,0,1 \ {\mbox {\rm and}} \ i \in
\un{L(k)} \}$ (where $L(k)$ is a nonnegative integer for each
$k=-1,0,1$) be a cut collection such that if $k < l$ then $D_{i}(k)
\not\succ D_{j}(l)$ for $i \in \un{L(k)}$ and $j \in \un{L(l)}$.  Then given 
$\var ,
\delta > 0$ there exist a $(\delta,e)$-collection $\{ \al_{i}(-1)
\subset D_{i}(-1) \}^{L(-1)}_{i=1}$ and a $(\var,c)$-collection $\{
\al_{j}(1) \subset D_{j}(1) \}^{L(1)}_{j=1}$ both compatible with $\{
(D_{i}(k), \beta_{i}(k) );\ k= -1,0,1 \ {\mbox{\rm and}} \ i \in
\un{L(k)} \}$.
\end{prop}

\noi  {\sc Proof}:  (See remark before the statement.)  We may assume,
without loss of generality, that the distinguished representatives in the
$c$-equivalence classes among $\{ D_{i}(1); \ i \in \un{L(1)} \}$ are the
first $n$ disks $D_{1}(1), \ldots, D_{n}(1)$.  For each $i \in \un{n}$
consider the cut collection $$\{ (D_{j}(k), \beta_{j}(k)); \ k=-1,0,1, \ j
\in \un{L(k)}, D_{j}(k) \not\succ D_{i}(1) \} \cup \{ ( D_{i}(1),
\beta_{i}(1) ) \}.$$  By Proposition~\ref{21} there exists an open
cross-cut $\al_{i}(1) \subset I^{c} (\beta_{i}(1) ) \cap V_{\var} (C_{i}(1))$
satisfying (i) and (ii) of that proposition (with $\al_{i}(1)$ in place
of $\al_{0}$.)  We do the same for every $i \in \un{n}$ obtaining $\{
\al_{i}(1) \}^{n} _{i=1}$.  These cross-cuts clearly satisfy (i) and
(ii) in the definition of $(\var,c)$-collections and $\al_{i}(1) \subset
I^{c} (\beta_{i}(1) ) \cap V_{\var}(C_{i}(1))$ by construction.  In order
to see they are disjoint, let $i,j \in \un{n}$.  If $I_{i}(1) \cap
I_{j}(1) = \es, \ \al_{i} (1) \cap \al_{j}(1)= \es$ since $\al_{i}(1)
\subset I_{i}(1)$ and $\al_{j}(1) \subset I_{j}(1)$.  If $I_{i}(1) \cap
I_{j}(1) \neq \es$, then either $D_{i}(1) \succ D_{j}(1)$ or $D_{j}(1)
\succ D_{i}(1)$, say, $D_{i}(1) \succ D_{j}(1)$.  It follows that
$\stackrel{\circ}{C_{i}} (1) \cap \ov{ I_{i}(1) \cap I_{j}(1) }= \es$ for
otherwise $C_{i}(1) \subset C_{j} (1)$ and $D_{i}(1), D_{j}(1)
|_{C_{i}(1)}$, which goes against our assumption that $C_{i}(1)$ was the
distinguished representative in its $c$-equivalence class.  From this we
can conclude that $[ I^{c} ( \al_{i} (1) ) \cup \al_{i} (1) ] \cap D_{j}
= \es$ and thus that $\al_{i}(1) \cap \al_{j} (1) = \es$.  Indeed we have
shown more, namely that $$[I^{c}(\al_{i}(1) ) \cup \al_{i}(1) ] \cap
[I^{c} ( \al_{j} (1) ) \cup \al_{j} (1)] = \es$$ for any $i,j \in \un{n}$.

\bigskip

We now look at the disks in one $c$-equivalence class.  By
Proposition~\ref{26} the $c$-sides of the elements in the class are
unlinked in the $c$-side of its distinguished representative, $D_{i}(1)$
say.  By Proposition~\ref{16} it is possible to find disjoint open
cross-cuts $\al_{j}(1) \subset I^{c}( \al _{i} (1) )$ joining the
endpoints of $C_{j}(1)$ such that $\al_{j}(1) \subset I^{c} ( \beta_{j}(1) ) 
\cap V_{\var} (
C_{j}(1) )$ for every $j$ such that $D_{j} \sim_{c} D_{i}$.  Doing this
for each $c$-equivalence class we find a collection of disjoint open
cross-cuts $\{ \al_{i}(1) \}^{L(1)}_{i=1}$ satisfying the conditions in
the definition of a $(\var,c)$-collection compatible with $\{ ( D_{i}(k),
\beta_{i}(k) ) \} \ \ \Box$.

\begin{figure}
\begin{center}~
\psfig{file=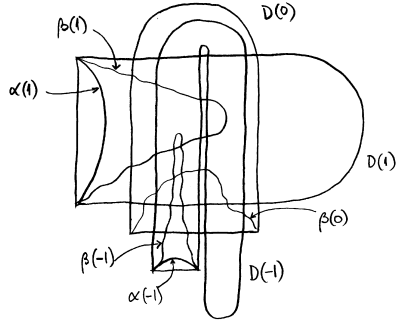,height=2.5in}
\end{center}
\caption{$\{ \al (1) \}$ is the ($\var,c$)-collection and $\{ \al (-1) 
\}$ is the ($\var,e$)-collection, both compatible with $\{ (D(k), 
\beta(k); k = -1,0,1 \}$} \label{f11}
\end{figure}

We will now introduce dynamics in our discussion and add to the definition
of $(c,e)$-disks a new requirement, as we promised earlier.  Let $f: \pi
\rightarrow \pi$ be a plane homeomorphism which we will have fixed for
the remainder of the time.
\bigskip

\noi (C,E) {\sc Dynamical Assumption}:   All $(c,e)$-disks henceforth will
be assumed to satisfy (i) and (ii):

\begin{description}

\item[(i)] $\dis{\lim_{n \rightarrow \infty} }$ diam $f^{n} (C) = 0$;

\item[(ii)] $\dis{\lim_{m \rightarrow - \infty} }$ diam $f^{m}(E) = 0$.

\end{description}

\noi The main purpose of the present work is to isotop away dynamics of
$f$ in a controlled manner.  We will now define sets within which it is
possible to do this, namely, to destroy all dynamics within them by an
isotopy which is identically equal to $f$ without them.  We call them
{\em pruning fronts}\/ after the work of Predrag Cvitanovi\'c~\cite{C}.

\begin{defn}
Let $\{D_{i} \}^{L} _{i=1}$ be a collection of $(c,e)$-disks (satisfying
the dynamical assumption above) such that (i), (ii) and (iii) hold:

\begin{description}

\item[(i)] $\succ$ can be extended by transitivity to a partial order 
on  $\{ D_{i} \}^{L} _{i=1}$
or, equivalently, there are no ``loops'' $D_{i_{1}} \succ D_{i_{2}} \succ
\ldots \succ D_{i_{n}} \succ D_{i_{1}}$;

\item[(ii)]  for every $n> 0$ and $i,j \in \un{L}, \ f^{n} (D_{i})
\not\prec D_{j}$;

\item[(iii)] for every $m <0$ and $i,j \in \un{L}, \ f^{m}(D_{i} )
\not\succ D_{j}$.

\end{description}

Such a collection will be called a {\em pruning collection}.  Its locus
$\ov{P} = \dis{\bigcup ^{L} _{i=1} } D_{i}$ (see \cite{C} and the
comments before the definition) will be called a {\em pruning front}.
\end{defn}

\noi {\sc Notation}:  We will use $\geq$ to denote the extension of
$\succ$ to a partial order and keep $\succ$ to denote the binary relation
as we defined previously.

\bigskip

Before we proceed, let us say a word about finite partially ordered
sets.  If $(X, \geq )$ is one such we define the set of {\em initial
elements}\/ of $X$ to be
$$ I(X) = \{x \in X; \ \forall y \in X, \ y \leq x \Longrightarrow y = x
\} $$
\noi It is easy to see that if $X$ is finite and nonempty, $I(X)$ is
nonempty and that no two distinct elements in $I(X)$ are related by
$\geq$.  Now let $X_{1} = I(X)$ and inductively set $X_{n} = I ( X \sm
\dis{\bigcup ^{n-1} _{i=1} } X_{i} )$.  From what we have said, $X_n$ is
nonempty if $X \sm \dis { \bigcup^{n-1}_{i=1} }X_{i}$ is nonempty.  Since
$X$ is finite, there exists $n \geq 1$ such that $X_{1}, X_{2}, \ldots,
X_{n}$ are all nonempty and for $m>n, \  X_{m} = \es$.  Clearly $X_{1},
\ldots, X_{n}$ is a partition of $X$ and if $X_i$ has $s_i$ elements we
can list the elements of $X = \{ x_{1}, x_{2}, x_{3}, \ldots, x_{L} \}$
so that the first $s_1$ elements are those in $X_1$, the next $s_2$
elements are those in $X_2$ and so on.  In this way the subscripts
reflect the partial order in the sense that if $i < j$ then $x_{i}
\not\geq x_{j}$.  Having said this we adopt the following
\bigskip

\noi{\sc Convention}:
Henceforth it will be assumed that the subscripts in a pruning collection
reflect the partial order $\geq$ in the sense that if $i < j$ then $D_{i}
\not\geq D_{j}$.  Notice that, in particular, if $i < j$ then $D_{i}
\not\succ D_{j}$.
\bigskip

We can now state a proposition containing one of the main ingredients in
the proof of the main theorem (see figure 12.)

\begin{prop}\label{28}
Let $\{D_{i} \}^{L}_{i=1}$ be a pruning collection and $\{\var_{n}
\}^{\infty}_{n=0}$ a sequence of positive numbers converging to zero.
Then there exists a collection $\{ \al_{i} (n) \subset f^{n}(D_{i}); \ i
\in \un{L}, \ n \in {\Bbb{Z}} \}$ of disjoint open cross-cuts joining the
vertices of $f^{n}(D_{i})$ such that (i) and (ii) below hold:

\begin{description}
\item[(i)]  For each $n \geq 1, \ \{\al _{i} (n); \ i \in \un{L} \}$ is a
$(\var_{n},c)$-collection compatible with
\begin{eqnarray*}
&&\{(f^{k}(D_{j}), \ \al_{j}(k));
\ j \in \un{L}, \ -n +1 \leq k \leq n-1 \} \\
&&\cup \{ ( f^{n} (D_{j}), \
f(\al_{j} (n-1) )); \ j \in \un{L} \}
\end{eqnarray*}

\item[(ii)]  For each $m \leq 0, \ \{ \al_{i} (m); \ i \in \un{L} \}$ is
a $(\var_{|m|},C)$-collection compatible with 
\begin{eqnarray*}
&&\{ (f^{k}(D_{j}),
\al_{j}(k)); \  j \in L,  m+1 \leq k \leq -m+1 \} \\
&& \cup \{ ( f^{m}(D_{j}
), f^{-1}  ( \al_{j} (m+1))); \ j \in \un{L} \}.
\end{eqnarray*}

\end{description}
\end{prop}

\noi{\sc Proof}:  We will let $m=-n+1$ and use induction on $n$.  In order
to prove the proposition for $n=1$, choose any collection $\{ \beta_{i}
\subset D_{i} \}^{L} _{i=1}$ of open cross-cuts joining vertices and
apply Proposition~\ref{27} with $L(0)=0$ (so that $\un{L(0)} = \es$ and
$\{ (D_{i}(0), \beta_{i} (0) ) \} = \es$) to the cut collection 
$${\mc{D}}
= \{ ( D_{i}, \beta_{i});\  i \in \un{L} \} \ \cup \{ ( f(D_{i}), \
f(\beta_{i})); \ i \in \un{L} \}$$ 
where $\{(D_{i}, \beta_{i} ) \}$ and 
$\{ (f (D_{i}), f(\beta_{i} )) \}$ play the roles of $\{(D_{i}(-1),
\beta_{i}(-1)) \}$ and $\{(D_{i}(1), \ \al_{i}(1) ) \}$ respectively in
the statement of that proposition, whereas $\var=\var_{1}$ and $\delta =
\var_{0}$.  By the definition of pruning collection, $f(D_{i}) \not\prec D_{j}$
for any $i, j \in \un{L}$ so that $\mc{D}$ satisfies the hypotheses and 
we can conclude there exist $\{\al_{i}(1) \}^{L}_{i=1}$ and $\{\al_{i}
\}^{L}_{i=1}$ a $(\var_{1},c)$- and a $(\var_{0},e)$-collection
respectively,
both compatible with $\mc{D}$.  Since $\al_{i} \subset I^{e}(\beta_{i} )$
and therefore $f(\al_{i}) \subset I^{e} (f (\beta_{i} ))$, by
Proposition~\ref{23}, and Corollary~\ref{24}, $\{\al_{i}(1) \}^{L} _{i=1}$
is a
$(\var_{1},c)$-collection compatible with $$\{ ( D_{i}, \al_{i} ); \ i \in
\un{L} \} \cup \{ ( f (D_{i} ), f( \al_{i} ) ); \ i \in \un{L} \}.$$  By
the same token $\{ \al_{i} \} ^{L}_{i=1}$ is a $(\var_{0},e)$-collection
compatible with $$\{(D_{i}, f^{-1} (\al_{i} (1) )); \ i \in \un{L} \} \cup \{
( f (D_{i}), a_{i}(1) ); i \in \un{L} \}.$$  That $\al_{i}(1)\cap \al_{j} =
\es$ for $i, j \in \un{L}$ is a consequence of Proposition~\ref{25}.
This proves the proposition for $n=1, \ m=0$.

\medskip

Assume we have constructed a collection $$\{\al_{i} (k); \ i \in \un{L},
\ -n +2 \leq k \leq n-1 \}$$ of disjoint open cross-cuts satisfying the
conclusions of the proposition.  Consider the cut collection
\begin{eqnarray*}
 {\mc{D}} &= &\{ ( f^{n} (D_{i}), \ f(\al_{i} (n-1) )); \ i \in \un{L} \} \\
&& \cup \ \{ (f^{k} (D_{i}), \al_{i}(k) ); \ i \in \un{L}, \
-n +2 \leq k \leq
n-1 \} \\
&& \cup \  \{ (f^{-n+1} (D_{i}), f^{-1} (\al_{i}(-n+2))); \ i
\in \un{L} \}
\end{eqnarray*}
\noi and apply Proposition~\ref{27} with $\{(D_{i} (1), \beta_{i}(1) )
\}, \
\{ ( D_{i} (0), \al_{i} (0) ) \}$ and \linebreak[3]
$\{ ( D_{i}$ $ (-1), \ \al_{i}(-1)) \}$
equal to the first, second and third collections respectively, in the
above union, letting $\var = \var_{n}$ and $\delta= \var_{|-n+1|}$.  From
the definition of pruning collection, $f^{n}(D_{i}) \not\prec f^{k}(D_{j})$ for
any $k < n$ and any $i, j \in \un{L}$ and $f^{-n+1}(D_{i}) \not\succ f^{k}
(D_{j})$ for any $k > -n+1$ and any $i,j \in \un{L}$, so that the hypotheses
of the proposition are satisfied.  We may then conlude there exist $\{
\al_{i}(n) \subset f^{n} (D_{i} ) \}^{L}_{i=1}$ and $\{\al_{i}(-n+1)
\subset f^{-n+1}(D_{i})\}^{L}_{i=1}$ a $(\var_{n}, c)$- and a
$(\var_{|-n+1|},e)$-collection respectively, both compatible with
$\mc{D}$.  From Corollary \ref{24}, $\{\al_{i}(n)
\}^{L}_{i=1}$ is compatible with 
\begin{eqnarray*}
&&\{( f^{k}(D_{i}), \al_{i} (k)); \ i \in
\un{L}, \ -n+1 \leq k \leq n-1 \} \\
&& \cup \{ (f^{n}(D_{i} ) , f(\al_{i}
(n-1))); \ i \in \un{L} \}
\end{eqnarray*} 
and $\{ \al_{i} (-n+1) \}^{L}_{i=1}$ is
compatible with 
\begin{eqnarray*}
&&\{ ( f^{k} (D_{i}), \al_{i}(k) ); \ i \in \un{L}, \ -n+2
\leq k \leq n \} \\
&& \cup \{ ( f^{-n+1} (D_{i}), f^{-1} ( \al_{i} (-n+2))); \
i \in \un{L} \}.
\end{eqnarray*}
  That $\al_{i}(n) \cap \al_{j}(k) = \es$ for $-n+1 \leq
k \leq n-1$ and $\al_{i} (-n+1) \cap \al_{j} (k) = \es$ for $-n+2 \leq k
\leq n$ is a consequence of Proposition~\ref{25}.  This finishes the
induction step and proves the proposition.  $\Box$

\bigskip

\begin{figure}
\begin{center}~
\psfig{file=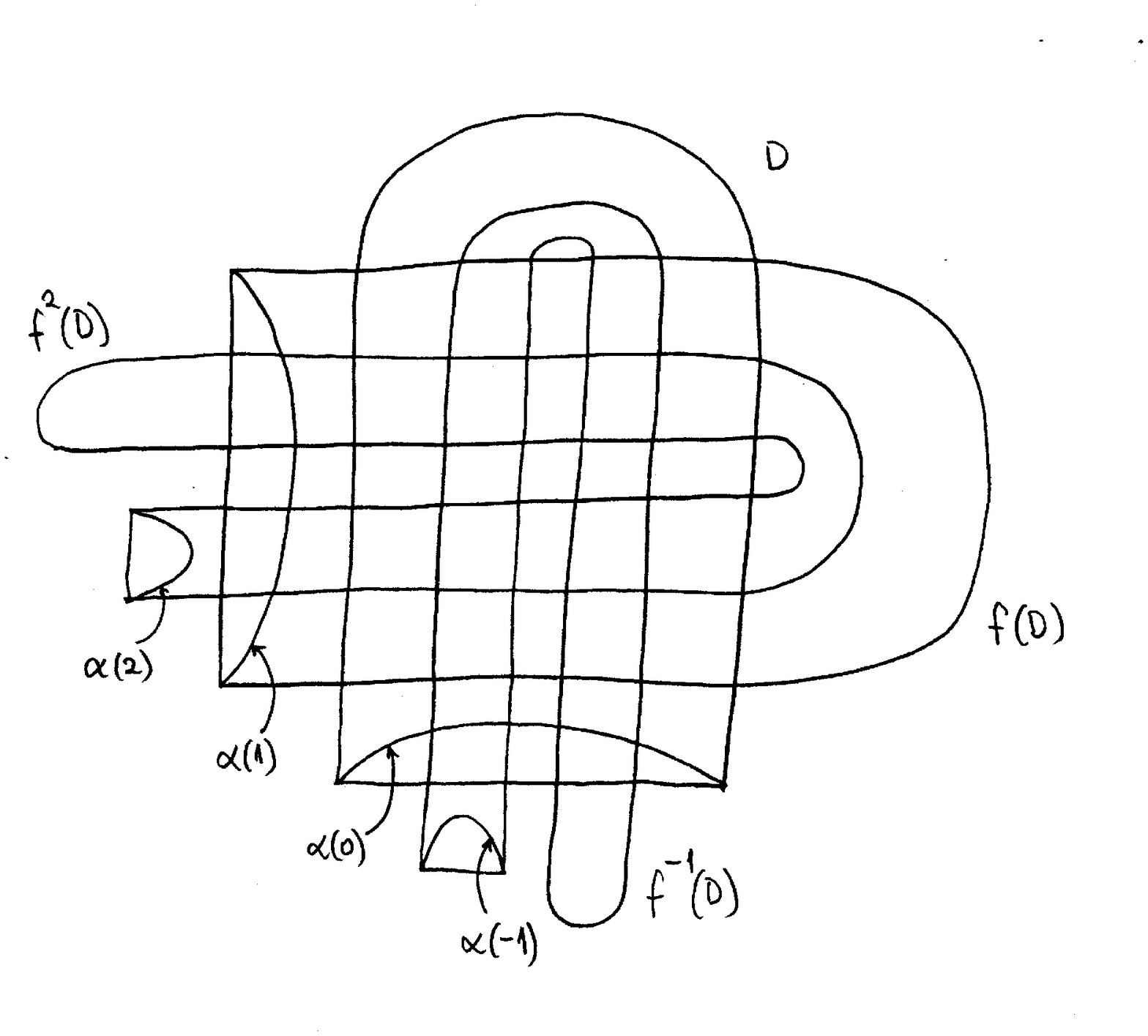,height=3in}
\end{center}
\caption{The first few $\al(n)$'s for a pruning collection containing 
only one ($c,e$)-disk $D$.} \label{f12} 
\end{figure}

\begin{cor}\label{29}
With the notation of Proposition~\ref{28}, for every $n \in \Bbb{Z}$,
 $\al_{i} (n) \subset I^{c} (f (\al_{i} (n-1)))$ and $\al_{i}(n) \subset
I^{e} (f^{-1} (\al_{i} (n+1)))$.
\end{cor}

\noi{\sc Proof}:  For $n \geq 1$, (i) of Proposition~\ref{28} implies
that $\al_{i}( \!n \!) \! \!\subset \!  I^{c} (f(\al_{i} ( \!n \!- \!1 \!) 
\!) \!)$ whereas (ii) 
implies that for $m \leq 0, \ \al_{i}(m) \subset$ $I^{e} ( f^{-1} (\al_{i}
(m+1)))$.  By Proposition \ref{17}, $f^{-1} (\al_{i} (m+1)) \subset I^{c}
(\al_{i} (m) )$ and applying $f$ to both sides we get $\al_{i} (m+1)
\subset f ( I^{c} ( \al_{i}(m))) = I ^{c}(f (\al_{i}(m)))$.  Letting $n=m+1$
we see that for $n \leq 1, \ \al_{i}(n) \subset I^{c} (f ( \al_{i}
(n-1)))$, which completes the proof of the first statement.  The second
is obtained from it using Proposition~\ref{17} (see figure 13.)   $\Box$
\bigskip

\begin{figure}
\begin{center}~
\psfig{file=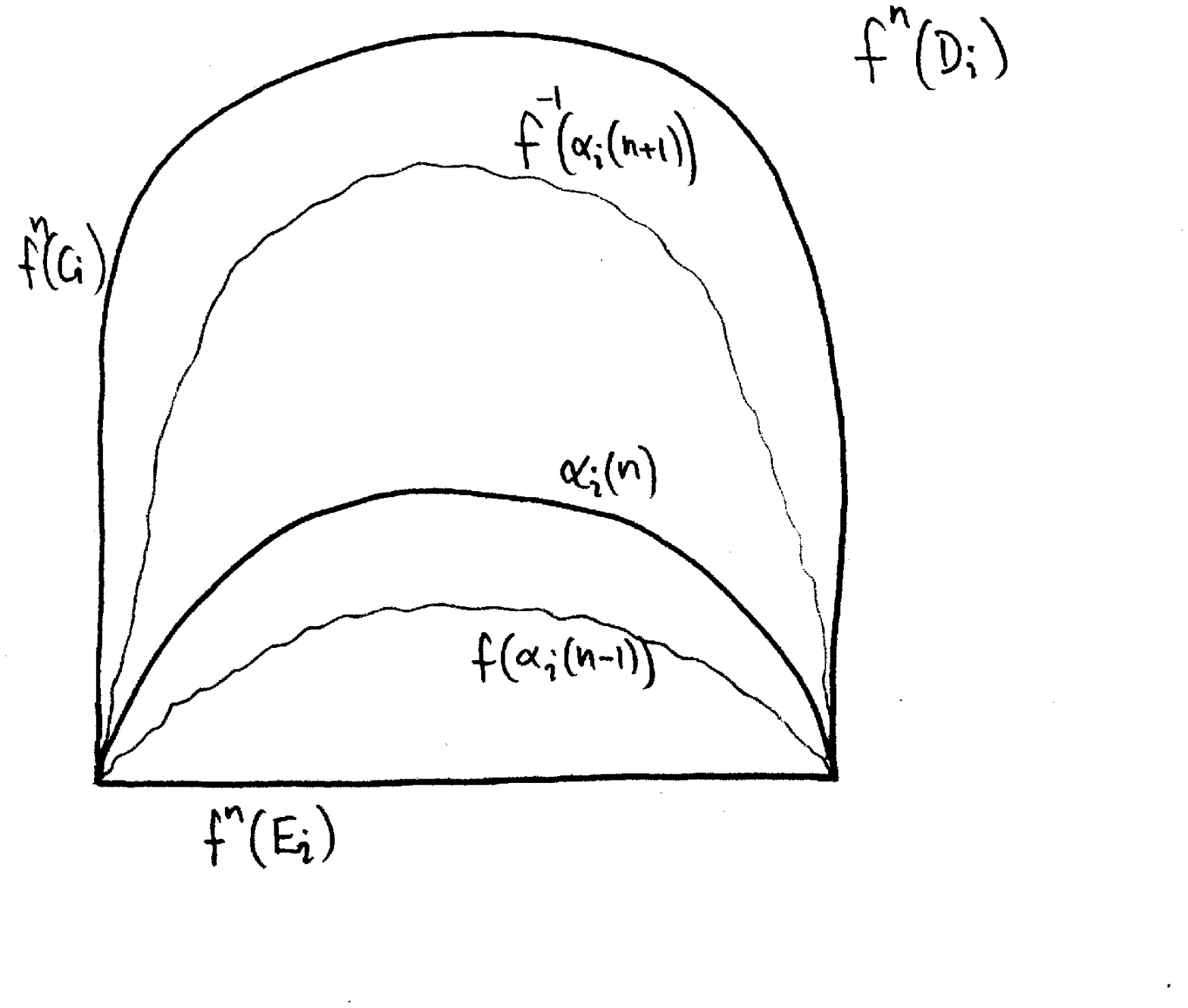,height=2.5in}
\end{center}
\caption{The $\al(n)$'s are chosen so that $\al_{i} (n) \subset I^{c} (f (\al_{i} (n-1)))$ and $\al_{i}(n) \subset
I^{e} (f^{-1} (\al_{i} (n+1)))$.} \label{f13}
\end{figure}

The next proposition is nothing but a ``fattened'' version of
Proposition~\ref{28} (see figure 14.)  We could have proven it together 
with Proposition~\ref{28} had we stated the ``fattened'' versions of the
propositions we proved before.  Although feasible, this would have been
rather cumbersome.  It is also possible to give a direct proof using 
the techniques we have used so far.  We leave it to the interested 
reader.
 \bigskip

\begin{prop}\label{30}
Let $\{\al_{i}(n); \ i \in \un{L}, \ n \in {\Bbb{Z}} \}$ be as in
Proposition~\ref{28}.  Then there exist collections of disjoint open
cross-cuts $\{\beta_{i}(n) \subset f^{n}(D_{i}); \ i \in \un{L}, \ n \in
{\Bbb{Z}} \}$ and $\{\gamma_{i}(n) \subset f^{n}(D_{i}); \ i \in \un{L}, \ n
\in {\Bbb{Z}} \}$ joining vertices such that:

\begin{description}

\item[(i)] $\beta_{i}(n) \subset I^{c} (\al_{i}(n))$ and $\gamma_{i}(n)
\subset I^{e} (\al_{i}(n))$;

\item[(ii)] for $n \geq 1, \ \{ \gamma_{i} (n); \ i \in \un{L} \}$ is a
$(\var_{n}, c)$-collection compatible with
\begin{eqnarray*}
&&\{(f^{k}(D_{i}),
 \ \beta_{i}(k)); \ i \in \un{L}, \ -n+1 \leq k \leq n-1 \} \\
&& \cup \{( f^{n}
(D_{i}), \ f(\beta_{i} (n-1))); \ i \in \un{L} \};
\end{eqnarray*}

\item[(iii)] for $m \leq 0, \ \{ \beta_{i} (m); \ i \in L \}$ is a
$(\var_{|m|},e)$-collection compatible with 
\begin{eqnarray*}
&&\{(f^{k}(D_{i}), \gamma_{i}
(k)); \ i \in \un{L}, \ m+1 \leq k \leq -m+1\} \\
&&\cup \{ ( f^{m} (D_{i} ),
f^{-1} (\gamma_{i} (m+1))); \ i \in \un{L} \}.\ \Box
\end{eqnarray*}

\end{description}
\end{prop}

\begin{figure}
\begin{center}~
\psfig{file=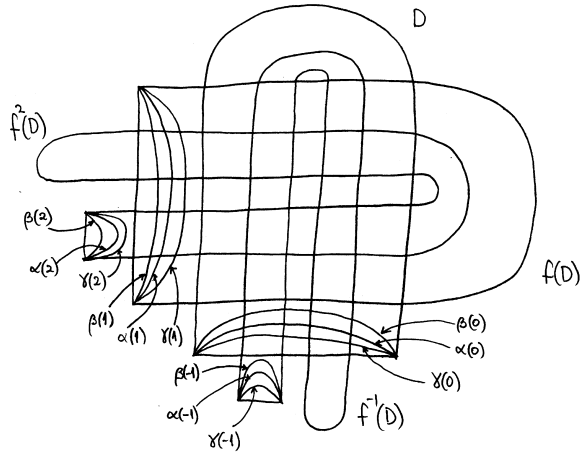,height=3in}
\end{center}
\caption{The first few $\gamma(n)$'s and $\beta(n)$'s for a 
pruning collection containing only one ($c,e$)-disk $D$.} \label{f14}
\end{figure} 

The corollary below is proved in the same way as Corollary~\ref{29}
(see figure~15.)
\bigskip

\begin{cor}\label{31}
With the notation of Proposition~\ref{30}, for every $n \in {\Bbb{Z}}$, 
$\gamma_{i}(n) \subset I^{c} (f(\beta_{i} (n-1)))$ and 
$\beta_{i}(n)
\subset I^{e} (f^{-1}( \gamma_{i} (n+1))). \ \ \Box$
\end{cor}

\begin{figure}
\begin{center}~
\psfig{file=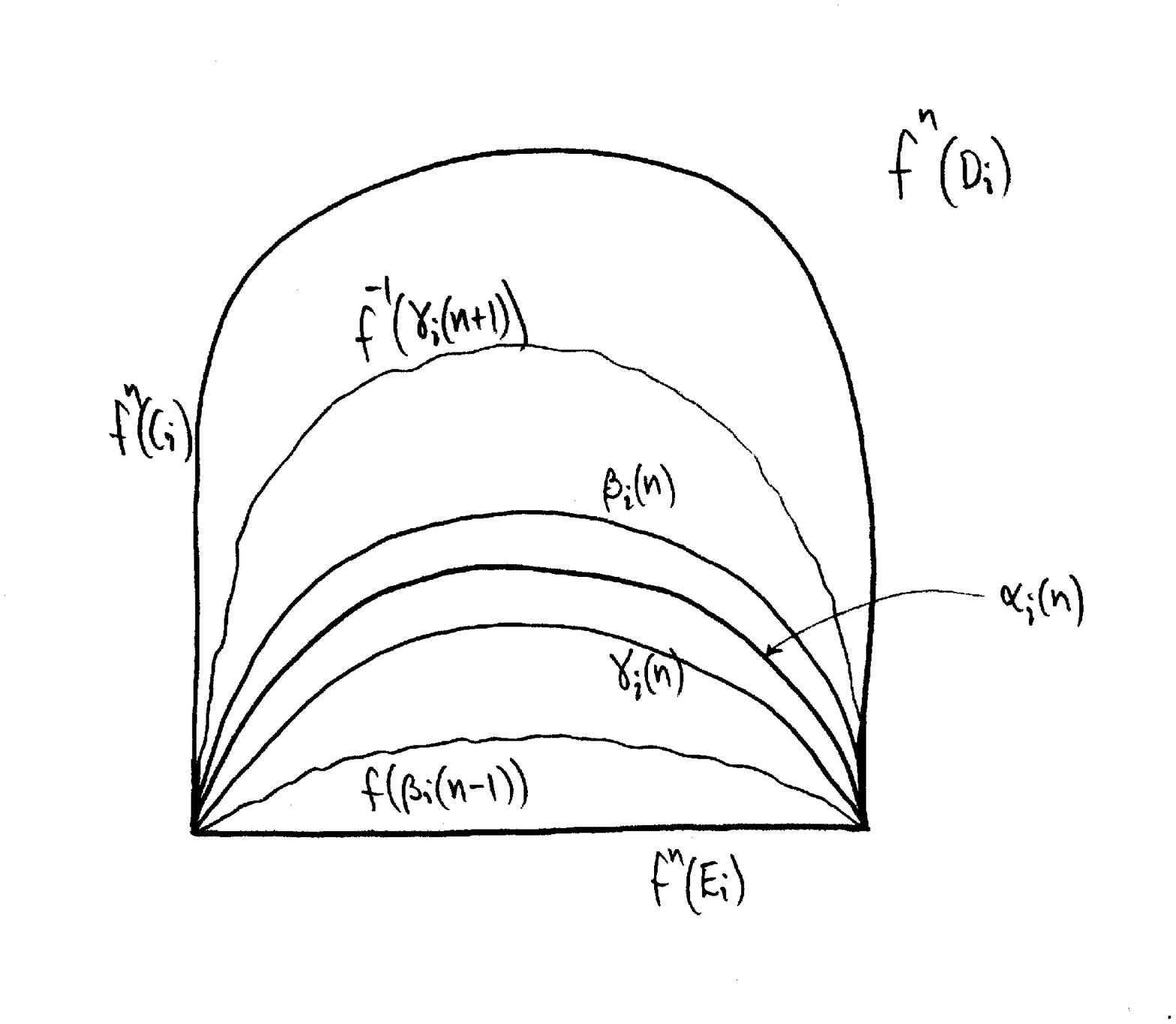,height=2.5in}
\end{center}
\caption{The $\beta_{i}(n)$'s and $\gamma_{i}(n)$'s  are chosen so that $\gamma_{i}(n) \subset I^{c} (f(\beta_{i} (n-1)))$ and $\beta_{i}(n)
\subset I^{e} (f^{-1}( \gamma_{i} (n+1)))$.}  \label{f15}
\end{figure}

The next proposition creates the sets in whose union will lie the
support of the isotopy we will construct to prove the main theorem.

\begin{prop}\label{32}
Let $\{\al_{i}(n) \}, \ \{\beta_{i}(n) \}$ and $\{ \gamma_{i} (n)\}$ be
as in Propositions \ref{28} and \ref{30}.
Then for every $n \in \Bbb{Z}$ and $i \in \un{L}, \ \ov{ f^{-1} (\beta _{i}
(n+1) ) \cup \gamma_{i} (n) }$ is a Jordan curve bounding a Jordan domain
${\mc{V}}_{i}(n)$ such that $${\mc{V}}_{i}(n) \supset f^{-1}(\al_{i} (n+1)) 
\cup \al_{i}
(n).$$  Moreover, $${\mc{V}}_{i}(n) = I^{c} (\gamma_{i}(n)) \cap I^{e} ( f^{-1}
(\beta_{i} (n+1))).$$
\end{prop}

\noi{\sc Proof}:  The proof is an easy exercise using (i) of
Proposition~\ref{30}, Corollary~\ref{31} and Proposition~\ref{17}
(see figure 16.)  $\Box$
\bigskip

\begin{figure}
\begin{center}~
\psfig{file=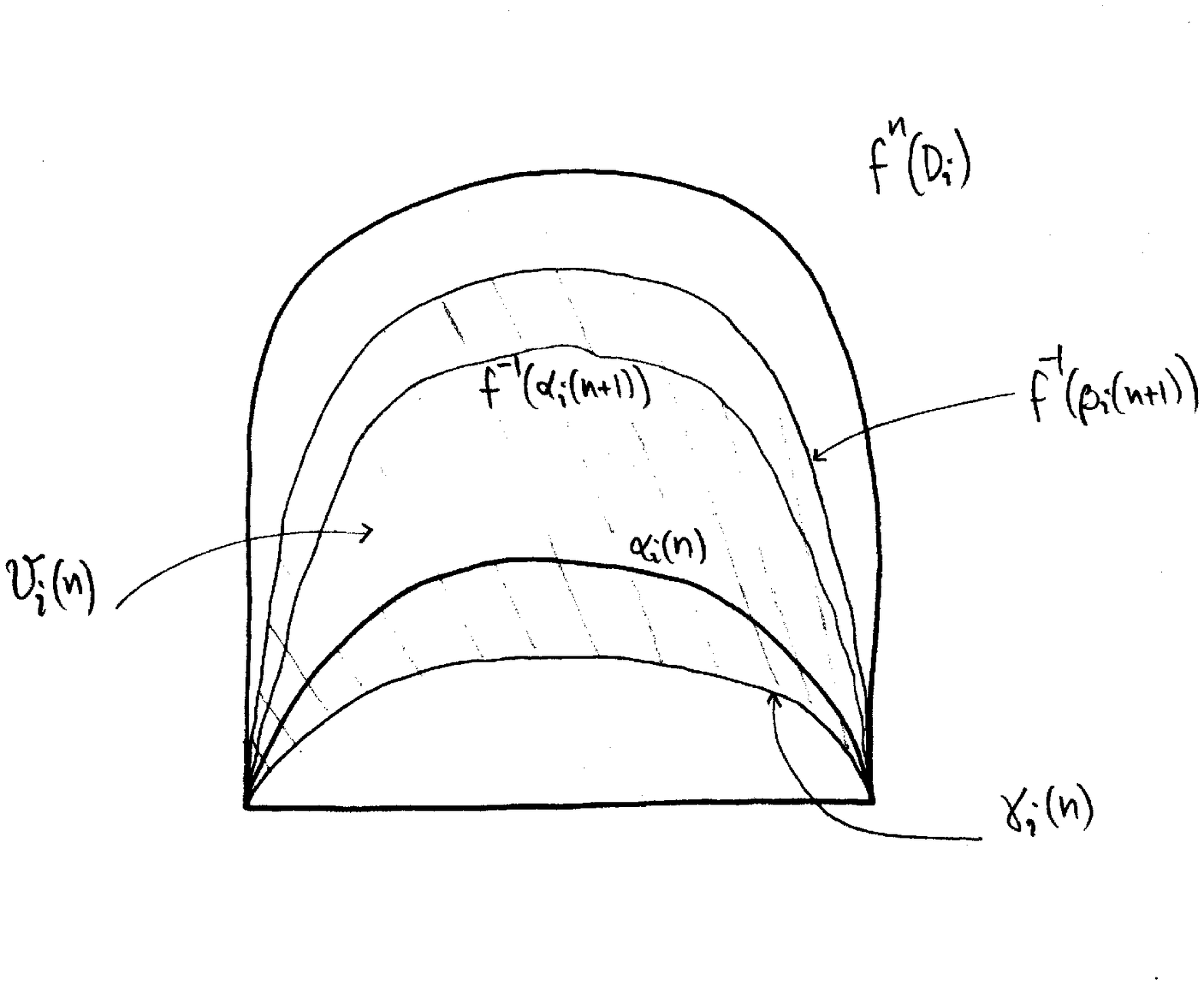,height=2.75in}
\end{center}
\caption{$\ov{f^{-1} (\beta_{i}(n+1))) \cup 
\gamma_{i} (n) }$ is a Jordan curve bounding the domain 
${\mc{V}}_{i}(n)$.} \label{f16}
\end{figure}

\begin{prop} \label{33}
Let $D_{1}, D_{2}$  be $(c,e)$-disks,  $D_{1} \not\prec D_{2}$ and
$\al_{1} \subset D_{1}$ and $\al_{2} \subset D_{2}$ be disjoint open
cross-cuts joining vertices.  Then $\al_{1} \cap I_{2} \subset I^{c}
(\al_{2})$ and $\al_{2} \cap I _{1} \subset I^{e} (\al_{1})$.
\end{prop}

\noi {\sc Proof}:  Since $D_{1} \not\prec D_{2}$ either $I_{1} \cap I_{2}
= \es$, in which case both statements are clearly true, or $D_{1} \succ
D_{2}$ and $D_{1} \neq D_{2}$.  If $D_{1} \succ D_{2}, \ C_{1} \cap I_{2}
= \es$ and since $\al_{1} \cap \al_{2} = \es, \ (\al_{1} \cup C_{1}) \cap
\al_{2}= \es$.  It follows, since $\al_{2}$ is connected, that either
$\al_{2} \subset I^{c}(\al_{1} ) $ or $\al_{2} \cap I^{c}(\al_{1})=
\es$.  We want to show that the latter is true, so we will assume
$\al_{2} \subset I^{c}(\al_{1})$ and reach a contradiction.  The
endpoints of $\al_{2}$ are the same as those of $E_2$ and since $E_{2}
\cap I_{1} = \es \ (D_{1} \succ D_{2} )$ and $\al_{1} \subset I_{1}$, if
$\al_{2} \subset I^{c} (\al_{1})$, it must be the case that the endpoints
of $\al_{2}$ lie on $C_1$.  But the endpoints of $\al_2$  coincide
with those of $C_2$ and, by (ii) in the definition of $\succ, \ C_{2}
\subset C_1$.  We claim that $C_{1} = C_{2}$, for if $C_2$ is strctly
contained in $C_1$, one of the endpoints of $\al_2$ lies in
$\stackrel{\circ}{C_{1}}$ and since $\al_{2} \subset I_{1} \cap I_2$,
(iii) in the definition of $\succ$ implies that $C_{1} \subset C_{2}$
which is a contradiction.  By Proposition~\ref{19} we see that $D_{1} =
D_{2}$ which is contrary to our hypothesis that $D_{1} \not\prec D_2$.

This contradiction shows that $\al_{2} \cap I^{c} (\al_{1}) = \es$ and
since $\al_{2} \cap \al_{1} = \es$ by hypothesis, we have shown that
$\al_{2} \cap I_{1} \subset I^{e}(\al_{1} )$.  The other statement is
proven analogously.  $\Box$
\bigskip

\begin{cor} \label{34}
Under the hypotheses of Proposition~\ref{33} $$I^{c} (\al_{1} ) \cap
I^{e} (\al_{2} )  = \es. \ \ \Box$$
\end{cor}

\begin{prop}\label{35}
Let $i,j \in \un{L}$ and $n, k \in \Bbb{Z}$:

\begin{description}

\item[(i)] if $f^{k} (D_{j}) \not\prec f^{n}(D_{i} )$ then $f^{-1}
(\al_{j} (k+1)) \cap \ov{ {\mc{V}}_{i} (n) } = \es$, and

\item[(ii)] if $f^{k}(D_{j}) \not\succ f^{n} (D_{i})$ then $\al_{j} (k)
\cap \ov{ {\mc{V}}_{i} (n) } = \es$.

\end{description}

\end{prop}

\noi{\sc Proof}:  From Proposition~\ref{30} (i) and
Proposition~\ref{17} it follows that $\al_{i}(n) \subset I^{e}
(\beta_{i}(n)) \cap I^{c} (\gamma_{i}(n))$.  From Proposition~\ref{30}
it also follows that $\beta_{j}(k) \cap \gamma _{i} (n)= \es$ for any
$i,j \in \un{L}, \ k,n \in \Bbb{Z}$.  Assume $f^{k}(D_{j}) \not\succ
f^{n}(D_{i})$.  By Corollary~\ref{34}, we see that $I^{e} (\beta_{j} (k))
\cap I^{c} (\gamma_{i} (n)) = \es$.
Since ${\mc{V}}_{i}(n) \subset I^{c}(\gamma_{i} (n)), \ \al_{j}(k) \subset 
I^{e}
(\beta_{j}(k) )$ and $I^{e}(\beta_{j}(k) )$ is open we can conclude that
$\ov{ {\mc{V}}_{i} (n)} \cap
\al_{j}(k)= \es$.  This proves (ii).  In order to prove (i) assume $f^{k}
(D_{j}) \not\prec f^{n}(D_{i})$.  Then $f^{k+1}(D_{j}) \not\prec
f^{n+1}(D_{i})$ and, as above, we can conclude that $I^{e} (\beta_{i}
(n+1)) \cap I^{c}(\gamma_{j} (k+1))=\es$.  It follows that 
$$I^{e}(f^{-1}
(\beta_{i} (n+1))) \cap I^{c}(f^{-1} ( \gamma_{j} (k+1))) = \es$$ and
since 
$${\mc{V}}_{i}(n) \subset I^{e} (f^{-1} (\beta_{i} (n+1))),$$
then 
$$ f^{-1}
(\al_{j} (k+1)) \subset I^{c} ( f^{-1} ( \gamma_{j} (k+1))). $$ 
This latter being an open set, we see that 
$$f^{-1} (\al_{j}(k+1)) \cap \ov{ 
{\mc{V}}_{i} (n) } = \es.$$  This completes the proof.  $\Box$

\bigskip

\begin{prop} \label{70}
With the notation above:
\begin{description}
\item[(i)] for $n \geq 1$ and $ -n +1 \leq k \leq n, \ f^{k} (C_{i} ) 
\cap f^{n} (I_{j}) \subset I^{c} ( \gamma_{j} (n) )$;

\item[(ii)] for $m \leq 0$ and $m \leq k \leq -m+1, \ f^{k} (E_{i} ) \cap 
f^{m}(I_{j} ) \subset I^{e} ( \beta_{j} (m))$.

\end{description}
\end{prop}

\noi{\sc Proof}:  From Proposition \ref{30} we know that $\{ \gamma_{j} (n) 
\}^{L}_{j=1} $ is compatible with $$\{ ( f^{k} (D_{i} ), \beta_{i} (k) ) ; 
i \in \un{L},-n +1 \leq k \leq n-1 \} \cup \{ (f^{n} (D_{i} ), 
f(\beta_{i} (n-1))); i \in \un{L} \}.$$  If $f^{k} (C_{i}) \cap f^{n} 
(I_{j} ) = \es$ there is nothing to prove.  Otherwise, $f^{n}(D_{j}) 
\succ f^{k} (D_{i} ) $ and therefore either $[ I^{c} ( \gamma_{j} (n)) 
\cup \gamma_{j} (n) ] \subset I^{c} (\beta_{i} (k) )$ or $[ I^{c} ( 
\gamma_{j} (n)) \cup \gamma_{j} (n) ] \cap f^{k} (D_{i} ) = \es$.  Since 
$f^{k} (C_{i} ) \subset f^{k} (D_{i} ) \cap {\mc{C}} I^{c} ( \beta _{i} 
(k))$, the conclusion of (i) follows.  $\Box$
\bigskip

\begin{cor} \label{67}
For $k \geq 1$, $f^{k}(C_{i})$ and $f^{-k}(E_{i})$ are disjoint from 
${\mc{V}}_{j} (n)$ for every $i, j \in \un{L}$ and every $n \in \Bbb{Z}$.
\end{cor}

\noi{\sc Proof}:  If $k > n$ by the definition of pruning collection 
$f^{k}(D_{i}) \not\prec f^{n}(D_{j} )$ which implies that $f^{k} (C_{i}) 
\cap f^{n} (I_{j}) = \es$.  Since ${\mc{V}}_{j}(n) \subset f^{n}(I_{j})$ 
this proves the result for $k > n$.  If $1 \leq k \leq n$ by Proposition
\ref{70}, $f^{k}(C_{i}) \cap {\mc{V}}_{j}(n) \subset I^{e} ( \gamma_{j}(n) 
)$ whereas by Proposition~\ref{32}, ${\mc{V}}_{j}(n) \subset I^{c} 
(\gamma_{j} (n) )$, which completes the proof of $f^{k} (C_{i}) \cap 
{\mc{V}}_{j}(n) = \es$ if $k \geq 1, \ j \in \un{L}$ and $n \in 
\Bbb{Z}$.  

If $n >k$ we again have by the definition of pruning 
collection that $f^{n} (D_{j}) \not\prec f^{k}( D_{i})$, which implies 
that $f^{k} (E_{i} ) \cap f^{n} (I_{j} ) = \es$ and thus that $f^{k} 
(E_{i}) \cap {\mc{V}}_{j}(n) = \es$.  If $m \leq k \leq 0 $, by 
Proposition \ref{70}, $f^{k} (E_{i} ) \cap f^{n} (I_{j} ) \subset I^{c} 
(\beta_{j}(n) )$ 
which implies that, if $n \leq k \leq -1$, $$f^{k} (E_{i} ) \cap f^{n} 
(I_{j} ) \subset f^{-1} (I^{c} (\beta_{j} (n+1))) = I^{c} (f^{-1} 
(\beta_{j} (n+1))).$$  By Proposition \ref{32}, ${\mc{V}}_{j}(n) \subset 
I^{e} (f^{-1} (\beta_{j} (n+1)))$ and thus $f^{k}(E_{i} ) \cap 
{\mc{V}}_{j} (n) = \es$ if $ k \leq 1, \ j \in \un{L}, \ n \in 
\Bbb{Z}$.  This completes the proof.  $\Box$
\bigskip

%% file: isotopies.tex
\section{Isotopies}

\begin{defn}
Let $X,Y$ be topological spaces.  By an {\em isotopy} we mean a
continuous map $H:X \times [0,1] \rightarrow Y$ such that the ``slice''
map $H_{t}: X \rightarrow Y, \ H_{t}(x)= H(x,t)$ is a homeomorphism for
each $t \in [0,1]$.  If $f,g: X \rightarrow Y$ are homeomorphisms, we
say $f$ and $g$ are {\em isotopic} if there exists an isotopy $H: X\times 
[0,1]
\rightarrow Y$ such that $H(x,0)=f(x)$ and $H(x,1)=g(x)$ for every $x \in
X$.

The {\em support} of an isotopy $H$ is by definition (see the remark 
below) the set
$$ \mbox{\rm supp } H = {\mc{C}} \{x \in X; \ H(x,t)=H(x,0) \ \forall t
\in [0,1] \} $$
\noi where, as usual, $\mc{C}$ stands for complement.

If $f:X \rightarrow X$ is a homeomorphism we define the support of $f$ as
$$ \mbox{\rm supp }f = {\mc{C}} \{ x \in X; \ f(x)= x \} $$
\end{defn}

\noi{\sc Remark}:  Notice that our definition of support is not the usual 
one in that we are not taking closures.  Supports of isotopies and 
homeomorphisms are therfore {\em open} sets.
\bigskip

The following proposition is a straightforward exercise in point set
topology and we omit the proof.

\begin{prop} \label{36}
Let $H:X\times [0,1] \rightarrow X$ be an isotopy of the identity, i.e.,
$H(x,0)=x$ for every $x \in X$.  If $x \in$ supp $H$, then $H(x,t)$ and
$x$ belong to the same path component of supp $H. \ \ \Box$
\end{prop}

\noi{\sc Remark:}  If $X$ is locally path-connected, the path components
of supp $H$ coincide with its connected components, since supp $H$ is open.
\bigskip

\begin{defn}
Let $G$ be a collection of subsets of a metric space.  We call $G$ a {\em
null collection} if for every $\var > 0$ only finitely many elements of
$G$ have diameter greater than $\var$.
\end{defn}

The lemma below is true in greater generality than we state and is part
of the folklore of hyperbolic geometry, geodesic laminations, etc.  The
proof we give is somewhat sketchy but is rather elementary.
\bigskip

\begin{lem}\label{37}
Let $I\!\!D$ denote the unit disk $\{ x \in I\!\!R^{2}; \ ||x|| \leq 1 \}$,
and $\{ \al_{n} \}^{\infty}_{n=1}$ a null collection of closed
cross-cuts, disjoint except possibly at endpoints, no two $\al_{n}$'s
sharing both endpoints.  For each $n \geq 1$, let $\gamma_n$ be the closed
arc of circle perpendiular to $S^{1}= \{x \in I\!\!R^{2}; \ ||x||=1 \}$
with the same endpoints as $\al_n$.  Then there exists a homeomorphism
$\zeta: I\!\!D \rightarrow I \!\! D$ such that $\zeta|_{S^{1} }$ is the
identity and $\zeta(\al_{n} ) = \gamma_{n}$.
\end{lem}

\noi{\sc Proof:}  From the hypotheses that the $\al_{n}$'s are interior
disjoint and no two share both endpoints it follows that the cross-cuts
$\gamma_n$ are interior disjoint and the correspondence $\al_{n} \rightarrow
\gamma_n$ is one-to-one in the sense that if $\al_{n} \neq \al_{m}$ then
$\gamma_{n} \neq \gamma_{m}$.  Moreover, $\{ \gamma_{n} \}^{\infty}_{n=1}
$ is a null collection, since given $\var > 0$ only finitely many pairs
of endpoints of the $\al_n$'s can be more than $\var$ apart, which
implies that only finitely many $\gamma_{n}$'s have diameter greater than
$\var$.

Let $\psi_{n}: \gamma_{n} \rightarrow \al_{n}$ be a homeomorphism
extending the identity homeomorphism between the endpoints of $\gamma_n$
and $\al_n$, for each $n \geq 1$, and define the map $\psi$ as
$$ \psi = {\mbox{\rm id}} \cup \bigcup^{\infty}_{n=1} \psi_{n}: S^{1} \cup
\bigcup ^{\infty}_{n=1} \gamma_{n} \longrightarrow S^{1} \cup \bigcup
^{\infty}_{n=1} \al_{n}$$
where id: $S^{1} \longrightarrow S^{1}$ is the 
identity homeomorphism. $\psi$ is well defined since the interiors
$\stackrel{\circ}{\gamma_{n} }$ are disjoint and $\psi_n$ is the identity
at the endpoints of $\gamma_n$.  We claim $\psi$ is a homeomorphism.
>From what we have said above, $\psi$ is clearly one-to-one and onto.  All
there remains to show is that $\psi$ is continuous.  Let $\{x_{k} \}$ be
a sequence in $S^{1} \cup \dis{\bigcup ^{\infty} _{n=1} } \gamma_{n}$ and
assume $x_{k} \rightarrow x$.  We want to show that $\psi (x_{k})
\rightarrow \psi (x)$.  If there exists $n$ such that all but finitely
many points $x_{k}$ lie in $\gamma_n$, then for $k_0$ sufficiently large
$x_{k} \in \gamma_{n}$ for every $k \geq k_0$ and since $\gamma_n$ is
closed $x \in \gamma_n$.  It follows that for $k \geq k_0, \ \psi(x_{k})
= \psi_{n} (x_{k} ) \rightarrow \psi_{n}(x) = \psi(x)$ since $\psi_n$ is
continuous.  If there is no $\gamma_n$ containing all but finitely many
$x_{k}$'s, we can choose a subsequence $x_{k_{j}} \in \gamma_j$ so that
different points lie in different $\gamma_{j}$'s.  Since $\{ \gamma_{n}
\}$ is a null sequence, diam $\gamma_{j} \rightarrow 0$ as $j
\rightarrow \infty$ and, since $x_{k_{j}} \in \gamma_j$ and $x_{k_{j}}
\rightarrow x$, for any sequence $y_{j} \in \gamma_{j}, \ y_{j}
\rightarrow x$.  In particular, if $p_{j}, q_{j}$ are the endpoints of
$\gamma_{j}, \ p_{j}, q_{j} \rightarrow x$.  This shows that $x \in
S^1$.  Also, the cross-cuts $\al_j$, whose endpoints are $p_{j}, q_{j}$,
are all distinct, since the $\gamma_{j}$'s are, and since $\{ \al_{n} \}$
is a null family and $p_{j}, q_{j} \rightarrow x$, for any sequence
$z_{j} \in \al_{j}, \ z_{j} \rightarrow x$.  We then have $\psi
(x_{k_{j}}) = \psi_{j} (x_{k_{j}} ) = z_{j} \in \al_{j}$ and $z_{j}
\rightarrow x = \psi (x)$ since $x \in S^1$.  This shows that $\psi$ is a
homeomorphism.  Assume for a moment we have shown that every component of
the complement of $S^{1} \cup \dis{\bigcup ^{\infty}_{n=1} } \gamma_{n}$
in $I \!\!D$ is a Jordan domain.  Let $U$ be one such and $\p U =J$.  $J$
is a Jordan curve in $S^{1} \cup \dis{\bigcup^{\infty}_{n=1} }
\gamma_{n}$ and thus $\psi(J)$ is a Jordan curve in $S^{1} \cup
\dis{\bigcup^{\infty}_{n=1} }\al_n$.  We claim that the Jordan domain $V$
bounded by $\psi(J)$ is a component of the complement of $S^{1} \cup
\dis{\bigcup^{\infty}_{n=1} } \al_n$ in $I \!\! D$.  It is clear that $V
\subset \{x;\ ||x|| < 1 \}$ so that $V \cap S^{1} = \es$.  If $V \cap
\al_{j} \neq \es$ for some $\al_j$, then $\stackrel{\circ}{\al_{j}}$,
which is connected and disjoint from $S^{1} \cup \dis{\bigcup_{n \neq j}}
\al_{n} \supset \p V$, is contained in $V$ and its endpoints in $\p V$.
But this implies that the endpoints of $\gamma_j$ lie on $J$ which in turn
implies $\gamma_{j} \subset \ov{U}$.  Since we assumed $U$ to be in the
complement of $S^{1} \cup \dis{\bigcup^{\infty}_{n=1} } \gamma_{n}, \
\gamma_{j} \subset J = \p U$.  This would then contradict the hypothesis
that no two $\al_n$'s shared both endpoints.  This shows that if $U$ is a
component of the complement of $S^{1} \cup \dis{\bigcup ^{\infty}_{n=1} }
\gamma_n$ in $I \!\! D$ whose boundary is a Jordan curve $J, \ \psi(J)$
is a Jordan curve in $S^{1} \cup \dis{\bigcup^{\infty}_{n=1} } \al_n$
bounding a component $V$ of the complement of $S^{1} \cup \dis{
\bigcup^{\infty}_{n=1} } \al_{n} $ in $I \!\! D$.  So if every
component $U$ of
the complement of $S^{1} \cup \dis{ \bigcup^{\infty}_{n=1}} \gamma_{n}$
in $I \!\! D$
is a Jordan domain we can use Theorem~\ref{4} to extend $\psi$ to a
homeomorphism $\tilde{\psi} : I \!\!D \rightarrow I \!\! D$ and $\zeta =
\tilde{\psi}^{-1}$ will satisfy the conclusions of the lemma.

In order to see that the components $U$ of $I \!\! D \sm \left[S^{1} \cup
\dis{ \bigcup ^{\infty}_{n=1} } \gamma_{n}\right]$ are Jordan domains let
$\gamma_n$ be a cross-cut such that $\gamma_{n} \subset \p U$.  Such a
$\gamma_n$ must exist unless $\{ \gamma_{n} \} = \es$ in which case the
statement is trivial.  By a conformal mapping, map $I \!\! D$ onto the
upper half plane $I \!\! H$ so that $\gamma_n$ maps onto $S^{1} \cap I
\!\! H$ and $U$ maps onto $U' \subset \{x; ||x|| < 1 \} \cap I \!\! H$.
It is now not hard to see that $\p U' \sm S^{1}$ is the graph of a
continuous function $g: (-1,1) \rightarrow [0,1)$ such that $|g(x)| <
\sqrt{1-x^{2} }$ for every $x \in (-1,1)$.  This proves that $\p U'$
is a Jordan curve and therefore so is $\p U$ and completes the proof of
the lemma.   $\Box$
\bigskip

\begin{cor}\label{38}
Let $J$ be a Jordan curve and $\{\al_{n} \}^{\infty}_{n=1} $ and $\{
\beta _{n} \}^{\infty}_{n=1}$ two null collections of interior disjoint
cross-cuts in $D$, the closed disk bounded by $J$.  Assume that no two
elements of each collection share both endpoints and that the endpoints
of $\al_i$ and $\beta_i$ coincide.  Then there exists a homeomorphism
$\zeta: D \rightarrow D$ such that $\zeta|_{J}$ is the identity,
$\zeta(\al_{n}) = \beta_{n}$ and $\zeta$ is isotopic to the identity
through an isotopy with support in $I$, the interior of $D$.

\end{cor}

\noi{\sc Proof:}  Let $f: D \rightarrow I \!\! D$ a homeomorphism and
$\zeta_{\al}: I \!\!D \rightarrow I \!\! D$ and $\zeta_{\beta}: I \!\! D
\rightarrow I \!\! D$ homeomorphisms ``straightening'' $\{f( \al_{n} )
\}$ and $\{ f(\beta_{n} ) \}$, which exist by Lemma~\ref{37}.  Set
$\zeta = \zeta^{-1}_{\beta} \circ \zeta_{\al}$.  It is not hard to check
that $\zeta|_{J} =$ id and $\zeta (\al_{n}) = \beta_{n}$.  That $\zeta$
is isotopic to the identity is a consequence of Theorem \ref{5}.  $\Box$
\bigskip

\begin{cor} \label{39}
Let $J$ be a Jordan curve and $\al, \beta$ cross-cuts in $D$ having the
same endpoints.  Then there exists an isotopy of the identity taking $\al$
to $\beta$ with support in $I$.
\end{cor}

\noi{\sc Proof:}  The collection $\{\al \}$ with a single element is a null
collection so it is possible to apply Lemma~\ref{37} and
Corollary~\ref{38}.  $\Box$
\bigskip

\noi{\sc Notation:}
Let $D_{1},D_{2}$ be closed disks, $D_{1} \subset D_{2}$ and $D_{1},
D_{2}|_{L}$, where $L \subset \p D_{1} \cap \p D_{2}$ is an arc.  If
$D_{1} \sm L \subset I_{2}$, the interior of $D_{2}$, we will write
$D_{1} \subset D_{2} |_{L}$.

\bigskip

\begin{lem} \label{40}
Let $\psi: D \rightarrow D$ be a homeomorphism onto its image so that
$\psi(D) \subset D|_{\psi(L)}$, where $L$ is a closed arc, $\psi(L)
\subset L$ and $p \in L$ is a fixed point such that $\psi^{n}(x)
\rightarrow p$ for every $x \in L$.  Then there exists an isotopy $h: D
\times [0,1] \rightarrow D$ of the identity such that $h|_{\p D}$= id and if 
$\zeta( \cdot ) = h
( \cdot, 1 ) $, then $(\psi \circ \zeta ) ^{n} (x) \rightarrow p$ for
every $x \in D$.
\end{lem}

\noi{\sc Proof:}  We will construct a null collection $\{ \al_{n} \}
^{\infty}_{n=1}$ of disjoint open cross-cuts in $I$ with the following
properties:

\begin{description}
\item[(i)] $\al_n$ has the same endpoints as $\psi^{n}(L)$;

\item[(ii)]  if $I_{n-1}$ is the Jordan domain bounded by $\al_{n-1} \cup
\psi^{n-1}(L), \ \al_n$ is a cross-cut in $I_{n-1} \cap \psi (I_{n-1} )$,
for $n \geq 2$;

\item[(iii)]  $\al_{n} \subset V_{\frac{1}{n} }(\psi^{n} (L) )$.

\end{description}

Set $\al_{1} = \psi (\p D \sm L )$ and $D_{1} = \psi (D)$.  Notice that
$D_{1} \subset D|_{\psi (L)}$ implies $\psi (D_{1}) \subset D_{1}
|_{\psi^{2}(L) }$.  By Proposition~\ref{16}, it is possible to find
$\al_{2} \subset \psi (I_{1}) \cap I_{1} = \psi (I_{1})$, an open
cross-cut joining the endpoints of $\psi^{2} (L)$ such that $\al_{2}
\subset V_{\frac{1}{2}}(\psi^{2}(L))$.

Assume we have constructed $\al_{1}, \al_{2}, \ldots, \al_{n}$ satisfying
(i), (ii) and (iii) above.  Since $\al_{n} \subset \psi(I_{n-1}) \cap
I_{n-1}$, and $\al_n$ has the same endpoints as $\psi^{n}(L), \ D_{n}
\subset \psi (D_{n-1} ) |_{\psi^{n}(L) }$ and $D_{n} \subset D_{n-1}|
_{\psi^{n}(L)}$.  This latter implies that $\psi(D_{n}) \subset \psi
(D_{n-1}) |_{\psi^{n+1}(L) }$
and since $\psi^{n+1} (L)
\subset \psi^{n}(L)$, by Proposition~\ref{13}, it follows that
$D_{n}, \psi (D_{n} ) | _ {\psi^{n+1}(L)} $.  By Proposition~\ref{16},
there exists $\al_{n+1} \subset I_{n} \cap \psi (I_n)$ an open cross-cut with
the same endpoints as $\psi^{n+1}(L)$ such that $\al_{n+1} \subset
V_{\frac{1}{n+1}} (\psi^{n+1} (L) )$.  By induction, we construct the
collection $\{ \al_{n} \}^{\infty} _{n=1} $.  That $\{ \al_{n}
\}^{\infty}_{n=1}$ is a null collection follows from the fact that
$\al_{n} \subset V_{\frac {1}{n}} (\psi^{n} (L) )$ and diam $\psi^{n} (L)
\rightarrow 0$.  That the $\al_n$'s are disjoint is clear since $\al_{n}
\subset I_{n-1}$ for every $n \geq 1$.  Notice also that no two $\al_n$'s
share both endpoints.  This is so because the endpoints of $\al_n$ are
the same as those of $\psi^{n}(L)$ and if $\psi^{n} (L)$ and $\psi^{m}
(L)$ shared both endpoints, $L$ would contain more than one fixed point.

Let $\beta_{n} = \psi^{-1}(\al_{n+1} )$.  The collection $\{ \beta _{n}
\} ^{\infty} _{n=1}$ is clearly a null collection of disjoint open
cross-cuts no two of which share both endpoints.  Also, for each $n \geq
1, \ \al_n$ and $\beta_n$ have the same endpoints.  By
Corollary~\ref{38} there exists an isotopy of the identity $h : D
\times [0,1] \rightarrow D$ such that if $\zeta (\cdot) = h ( \cdot, 1), \
\zeta (\al_{n} ) = \beta_n$.  Then $\psi \circ \zeta (\al_{n}) = \psi
(\beta_{n}) = \al_{n+1}$ and since $\psi ( \psi^{n} (L)) = \psi^{n+1}
(L)$ we see that $\psi \circ \zeta (D_{n} ) = D_{n+1} $.  But diam $D_{n}
\rightarrow 0$ as $n \rightarrow \infty$ and therefore it follows that
$(\psi \circ \zeta )^{n} (x) \rightarrow p, \ \forall x \in D$ as $n
\rightarrow \infty$ as we wanted.  $\Box$
\bigskip

\begin{cor} \label{41}
For $i=1, \ldots, n$ let $D_i$ be closed disks with disjoint interiors and 
$L_{i} \subset \p
D_i$ a closed arc.  Let $\psi: \pi \rightarrow \pi$ be a homeomorphism of
the plane such that $\psi(L_{i} ) \subset L_{i+1}$ and $\psi(D_{i})
\subset D_{i+1} |_{\psi(L_{i}) }$, where we let the indices ``wrap
around'', i.e., we set $n+1$ to be $1$.  Assume $\psi^{n}|_{D_{1}}:
(D_{1}, L_{1} ) \rightarrow (D_{1}, L_{1} )$ satisfies the hypotheses of
Lemma~\ref{40}.  Then there exists an isotopy $h: \pi \times [0,1]
\rightarrow \pi$ of the identity such that supp $h \subset D_1$ and if
$\zeta ( \cdot) = h (\cdot, 1), \ ( \psi \circ \zeta)^{kn} (x)
\rightarrow p$ as $k \rightarrow \infty$ for every $x \in D_1$ where $p
\in L_1$ is the fixed point of $\psi^{n}|_{L_{1}}$.
\end{cor}

\noi{\sc Proof:}  The proof is straightforward using Lemma~\ref{40} and
we omit the details.  $\Box$
\bigskip

%% file: theorem.tex
\section{The Proof of the Main Theorem}

In what follows $f: \pi \rightarrow \pi$ will be a uniformly continuous
homeomorphism of the plane and $\{ D_{i} \}^{L} _{i=1}$ a pruning
collection for $f$.  As we pointed out before, we may and will assume
that the subscripts reflect the partial order $\geq$ in $\{ D_{i}
\}^{L}_{i=1}$, in the sense that, if $i > j$ then $D_{i} \not\leq D_j$.
In particular, if $i > j $ then $D_{i} \not\prec D_j$.

\bigskip

\begin{defn}
For each $i \in L$ we define four numbers $n(i),\  N(i), \ m(i),$ $ M(i)
\in {\Bbb{Z}} \cup \{ \pm \infty \}$ as follows:

\begin{description}
\item[(i)] $n(i)$ is the smallest integer $\geq 1$ such that $f^{n(i)}
(D_{i} ), f(D_{j} ) | _{f^{n(i) }(C_{i})}$ for some $j \in \un{L}$ or
$n(i) = \infty$ if $f^{k} (D_{i} ), f(D_{j} ) \not|_{f^{k}(C_{i} ) }$ for
every $k \geq 1$ and $j \in \un{L}$;

\item[(ii)]  $N(i) = \left\lceil \frac{n(i)}{2} \right\rceil$, i.e., the
smallest integer greater than or equal to $\frac{n(i)}{2}$, if $n(i) <
\infty$ or $N(i)= \infty$ if $n(i)= \infty$;

\item[(iii)]  $m(i)$ is the largest integer $\leq 0$ such that $f^{m(i)}
(D_{i} ), D_{j} | _{f^{m(i)} (E_{i} )}$ for some $j \in \un{L}$ or $m(i)=
- \infty$ if $f^{k} (D_{i} ), D_{j} \not| _{f^{k}(E_{i})}$ for
every $k
\leq 0$ and $j \in \un{L}$;

\item[(iv)]  $M(i) = \left\lceil \frac{m(i)}{2} \right\rceil$ if $m(i) >
- \infty$ or $M(i) = - \infty$ if $m(i) = - \infty$.

\end{description}

\end{defn}

The following proposition is a straightforward consequence of the
definitions and we omit the proof.
\bigskip

\begin{prop}\label{66}
If $n(i), \ N(i), \ m(i)$ and $M(i)$ are finite the following holds for
each $i \in \un{L}$:

\begin{description}

\item[(i)]  $n(i) = 2 N(i) - \delta$ and $m(i) = 2 M(i) - \delta '$ where
$\delta, \delta'= 0$ or $1$;

\item[(ii)]  $f^{N(i)} (D_{i} ),\ f^{-N(i) + \delta +1 } (D_{j}) |
_{f^{N(i) } (C_{i} )}$ 
for some $j \in \un {L}$ but for $-N(i) + \delta
+2 \leq k \leq N(i) -1$, 
\[\ f^{N(i)} (D_{i}), f^{k} (D_{j}) \not|_{f^{N(i) } (C_{i})}\] 
for any $j \in \un{L}$;

\item[(iii)]  for $1 \leq n < N(i), \ -n+1 \leq k \leq n-1, \
f^{n}(D_{i}), f^{k} (D_{j}) \not| _{f^{n}(C_{i}) }$ for any $j \in \un{L}$;

\item[(iv)]  $f^{M(i)}(D_{i}), f^{-M(i) + \delta'} (D_{j}) |_{f^{m(i)}
(E_{i})}$ for some $j \in \un{L}$ but for $M(i) +1 \leq k \leq -M(i)+\delta'
+1$,
\[ f^{M(i)} (D_{i}), f^{k} (D_{j}) \not| _{f^{M(i)}(E_{i})}\]
for any $j \in \un{L}$;

\item[(v)]  for $M(i) < m \leq 0$ and $m +1 \leq k \leq -m+1, \ f^{m}
(D_{i} ), f^{k}(D_{j} ) \not| _{f^{m}(E_{i} ) }$ for any $j \in \un{L}$.
$\Box$

\end{description}
\end{prop}

Recall that we defined $c$- and $e$-equivalence relations in a collection
$\{D_{i} \}^{L}_{i=1}$ of $(c,e)$-disks and in Proposition~\ref{26}
proved that the equivalence classes have distinguished representatives.
The following proposition is again an easy consequence of the definitions.
\bigskip

\begin{prop}\label{42}
For each $i \in \un{L}, \ n(i) > 1$ if and only if $D_i$ is the
distinguished representative in its $c$-equivalence class in $\{D_{i}
\}^{L} _{i=1}$.  Likewise, $m(i) < 0$ if and only if $D_i$ is the
distinguished representative in its $e$-equivalence class in $\{D_{i}
\}^{L}_{i=1}$.  $\Box$
\end{prop}

We now start the construction of the isotopy for the proof of the main
theorem. If the pruning collection contains only one $(c,e)$-disk $D_1$
and $N(1)=\infty$ and $M(1)=-\infty$, most of what is presented
>from here to the end of this section is very much simplified.
We suggest that the reader concentrate on this case upon a first reading.

\bigskip

Recall that ${\mc{V}}_{i}(n) \subset f^{n} (I_{i} )$ is a Jordan domain
containing $\al_{i}(n)$ and $f^{-1} (\al_{i} (n+1) )$ as cross-cuts with
the same endpoints.  Using Corollary~\ref{31} construct, for each $i
\in \un{L}$ and $M(i) \leq n < N(i)$ an isotopy $k_{i,n} : \pi \times [0,1]
\rightarrow \pi$ of the identity such that supp $k_{i,n} \subset
{\mc{V}}_{i}(n)$ and $k_{i,n} (\al_{i} (n),1 ) = f^{-1} (\al_{i} (n+1) 
)$.  If
$n < M(i)$ or $n \geq N(i)$ we let $k_{i,n} \equiv$ identity.  Set
$\zeta_{i,n} (\cdot) = k_{i,n} (\cdot , 1)$.  For $n \in {\Bbb{Z}}$ define
$$k_{n}(x,t) =\left\{
\begin{array}{ll}
k_{1,n} (x,Lt),   & t \in \left[ 0, \dis\frac{1}{L} \right] \\
       \\
\zeta_{1,n} (k_{2,n} (x, Lt -1)),  & t \in \left[ \dis\frac{1}{L},
\dis\frac{2}{L} \right] \\
         \\
\zeta_{1,n} \circ \zeta_{2,n} (k_{3,n} (x, Lt -2)),  & t \in \left[
\dis\frac{2}{L}, \dis\frac{3}{L} \right] \\
           \\
\vdots   & \vdots \\
             \\
\zeta_{1,n} \circ \zeta_{2,n} \circ \ldots \circ \zeta _{L-1,n}
(k_{L,n}(x, Lt -L+1)), & t \in \left[ \dis\frac{L-1}{L}, 1 \right]
\end{array}
\right.  $$
\noi and let $\zeta_{n} (\cdot) = k_{n} (\cdot, 1)$.  Now let $r_{0} =
k_0$ and for $n \geq 1$
$$
r_{n} (x,t) = \left\{
\begin{array}{lr}
k_{-n} (x, 2t),   & t \in \left[ 0, \dis\frac{1}{2} \right] \\
    \\
\zeta_{-n} (k_{n} (x, 2t -1) ),   & t \in \left[ \dis\frac{1}{2}, 1
\right] \\
\end{array}
\right.
$$
\noi and set $\rho_{n} (\cdot) = r _{n} (\cdot, 1 )$ for $n \geq 0$.

Recall that the locus $\ov{P} = \dis\bigcup ^{L} _{i=1} D_i$ of a pruning
collection $\{ D_{i} \} ^{L} _{i=1}$ was called a pruning front.  We will
denote the union of the interiors $\dis\bigcup^{L}_{i=1} I_i$  by 
$P$. 
\bigskip

\begin{prop}\label{43}
The isotopies $r_n$ just defined have the following properties:

\begin{description}

\item[(i)] supp $r_{n} \subset [ f^{n} (P) \cup f^{-n} (P) ] \sm
\dis\bigcup  ^{n-1}_{-n+1} f^{k} (\ov{P} )$ for every $n \geq 0$ so that
if $n \neq m, \ \mbox{\rm supp } r_{n} \cap \mbox{\rm supp }r_{m} = \es$;

\item[(ii)]  since $f$ is uniformly continuous, the diameters of the
connected components of supp $r_{n}$ converge to $0$ as $n \rightarrow
\infty$;

\item[(iii)]  for each $i \in \un{L}$, if $n < N(i), \ \rho_{n} (\al_{i}
(n) ) = f^{-1} (\al_{i} (n+1))$ and if $-n \geq M(i)$,  $\rho_{n} (\al_{i}
(-n) ) = f^{-1} (\al_{i} (-n +1))$.

\end{description}

\end{prop}

\noi{\sc Proof:}  From the definition of $k_n$ it is clear that for $n \in
{\Bbb{Z}}$, $$\mbox{\rm supp }k_{n} \subset \bigcup \{ {\mc{V}}_{i} (n); \ M(i) 
\leq n < N(i) \}$$
so that for $n \geq 0$
$$ \mbox{\rm supp }r_{n} \subset \bigcup \{ {\mc{V}}_{i} (n); \ 1 \leq n < 
N(i) \} \cup \bigcup \{ {\mc{V}}_{i} (-n); \ M(i) \leq -n \leq 0 \} $$
\noi and since ${\mc{V}}_{i}(n) \subset f^{n} (I_{i})$, it is clear 
that 

$$\mbox{\rm supp } r_{n} \subset f^{n} (P) \cup f^{-n}(P)= \dis \bigcup^{L}_{i=1} f^{n}
(I_{i} ) \cup \dis \bigcup^{L}_{i=1} f^{-n}(I_{i} ).$$   There is nothing
more to prove for $n=0$ (recall that $\dis \bigcup^{-1}_{1} f^{k}(P) =
\es$, by our convention) and we may assume that $n \geq 1$ (see figure~17.)

If $1 \leq n < N(i)$, it follows from Proposition \ref{66} that $$f^{n}
(\!D_{i}\!), f^{k} (\!D_{j}\!) \not|_{f^{n}(\!C_{i}\!) }$$ for any $j \in 
\un{L}$ and
$-n +1 \leq k \leq n-1$.  Since $\{ \gamma_{i}(n) \}^{L}_{i=1}$ is a
$(\var_{n},c)$-collection compatible with $$\{ ( f^{k} (D_{j}), \beta_{j}
(k) ) : \ j \in \un{L}, \ -n +1 \leq k \leq n-1 \},$$ by
Proposition~\ref{30}, we must have $\left[ I ^{c} (\gamma_{i} (n) ) \cup
\gamma_{j} (n) \right] \cap f^{k} (D_{j} ) = \es$ for every $j \in
\un{L}$ and $-n+1 \leq k \leq n-1$.  But ${\mc{V}}_{i} (n) \subset I^{c}
(\gamma_{i} (n) )$ by Proposition~\ref{32} and taking the union
over $j \in \un{L}$ and $-n+1 \leq k \leq n-1$ we see that
$$ {\mc{V}}_{i}(n) \cap \dis\bigcup^{n-1}_{-n+1} f^{k} (\ov{P} ) = \es $$
\noi from which it follows that
$$ \bigcup \{ {\mc{V}}_{i} (n);\ 1 \leq n < N(i) \} \cap 
\dis\bigcup^{n-1}_{-n+1} f^{k} (\ov{P}) = \es \ . $$

If $M(i) < m \leq 0$, it follows from Proposition~\ref{66} that
$$f^{m}(\!D_{i}\!), f^{k}(\!D_{j} \!) \not|_{f^{m}(\!E_{i}\!)}$$ for any $j 
\in \un{L}$
and $m+1 \leq k \leq -m+1$.  Again by Proposition~\ref{30} $\{\beta_{i}
(m) \}^{L} _{i=1}$ is a $(\var_{|m|},e)$-collection compatible with $$\{ (
f^{k} (D_{j} ), \beta _{j} (k) ); \ j \in \un{L}, \ m+1 \leq k \leq -m+1
\},$$ which implies that $[ I^{e} (\beta _{i} (m) ) \cup \beta _{i} (m) ]
\cap f^{k} (D_{j} ) = \es$ for every $j \in \un{L}$ and $m+1 \leq k \leq
-m+1$, and thus that $[I^{e} (f^{-1} (\beta _{i} (m) )) \cup f^{-1} (\beta
_{i} (m) ) ] \cap f^{k} (D_{j}) = \es$ for every $j \in \un{L}$ and $m
\leq k \leq -m$.  Letting $m = -n +1$ and noticing that ${\mc{V}}_{i}(-n)
\subset I^{e} (f^{-1} (\beta _{i} (-n+1)))$, by Proposition~\ref{32},
what we have just seen implies that for $M(i) \leq -n < 0$
$$ {\mc{V}}_{i}(-n) \cap \dis\bigcup ^{n-1}_{-n+1} f^{k} (\ov{P}) = \es $$
\noi from which it follows that
$$ \bigcup \{ {\mc{V}}_{i} (-n); \ M(i) \leq -n < 0 \} \cap \dis
\bigcup^{n-1}_{-m+1} f^{k} ( \ov{P} ) = \es \ . $$

This finishes the proof of (i).

\bigskip

In order to prove (ii) notice that supp $r_{n} \subset \dis{ \bigcup
^{L}_{i=1} } \{ {\mc{V}}_{i} (n) \cup {\mc{V}}_{i} (-n) \}$ and from
Propositions~\ref{30} and \ref{32}, for $n \geq 1, {\mc{V}}_{i} (n) \subset
I^{c} ( \gamma_{i} (n) ) \subset V_{\ep_{n}}(f^{n}(C_{i}) ) $ and 
$${\mc{V}}_{i}
(-n) \subset I^{e} (f^{-1} (\beta _{i} (-n +1 ))) = f^{-1}(I^{e}(\beta_{i}
(-n+1))) \subset f^{-1} (V_{\ep_{n-1}} (f^{-n+1}(E_{i} ))).$$

>From the $(c, e)$ dynamic assumption, diam $f^{n}(C_{i}) \rightarrow 0$ as
$n \ra \infty$ and diam $f^{m} (E_{i}) \ra 0$ as $m \ra - \infty$.  Since
$\ep_{n} \ra 0$, it is clear that diam ${\mc{V}}_{i} (n) \ra 0$, as $n \ra
\infty$ and from the uniform continuity of $f$ we can also conclude that
diam ${\mc{V}}_{i} (-n) \ra 0$ as $n \ra \infty$.  It is now easy to see 
that the
connected components of $\dis{ \bigcup^{L} _{i=1} } \{ {\mc{V}}_{i} (n) \cup
{\mc{V}}_{i} (-n) \}$ have diameters converging to zero as $n \ra \infty$.  
This proves (ii).
\bigskip

Let us now look at (iii).  From the way we indexed the pruning
collection, if $i > j, \ D_{i} \not\prec D_j$ which implies $f^{n}(D_{i})
\not\prec f^{n} (D_{j})$ for any $n \in \Bbb{Z}$.  From
Proposition~\ref{35} it follows that $f^{-1} (\al_{i} (n+1)) \cap 
{\mc{V}}_{j}
(n) = \es$.  Similarly, if $l > i$ the same proposition implies that
$\al_{i}  (n) \cap {\mc{V}}_{l}(n)=\es$.
Since each $k_{i,n}$ is an isotopy of the identity with support contained
in ${\mc{V}}_{i}(n)$, and $\zeta_{i,n}(\cdot) =k_{i,n} (\cdot, 1)$, we have 
supp
$\zeta_{i,n} \subset {\mc{V}}_{i} (n)$ and, from what we said above, we see 
that if $j < i$, $\ \zeta_{j,n} (f^{-1} (\al_{i} (n+1))) = f^{-1} (\al_{i}
(n+1))$ and that if $l > i$, $ \zeta_{l,n} (\al_{i}(n) )= \al_{i}(n)$.
Thus, for any $M(i) \leq n < N(i)$
\begin{eqnarray*}
\zeta_{n} (\al_{i} (n) ) & = &\zeta_{1,n} \circ \ldots \circ \zeta_{i,n}
\circ \ldots \circ \zeta_{L,n} (\al_{i}(n) ) \\
  &= &\zeta_{1,n} \circ \ldots \circ \zeta_{i,n} (\al_{i}(n)) \\
  & = &\zeta_{1,n} \circ \ldots \circ \zeta_{i-1,n} (f^{-1} (\al_{i}
(n+1))) \\
  & = & f^{-1} (\al_{i} (n+1)).
\end{eqnarray*}

>From the definition of pruning collections, $f^{-n} (D_{i}) \not\succ f^{n}
(D_{j})$ for any $n \geq 1$ and any $i,j \in \un{L}$ and, by
Proposition~\ref{35}, it follows that $f^{-1} (\al_{i}(n+1) )$  $\cap
{\mc{V}}_{j} (-n) = \es$ and ${\mc{V}}_{i} (n) \cap \al_{j} (-n) = \es$.  
Thus we can
conclude that for any $i \in \un{L}, f^{-1} (\al_{i} (n+1)) \cap {\mbox
{\rm supp } }\zeta_{-n} = \es$ and that $\al_{i}(-n) \cap$ supp
$\zeta_{n} = \es$, for $n \geq 1$.  Therefore if $1 \leq n < N(i)$,
\begin{eqnarray*}
\rho_{n} (\al_{i} (n)) & = &\zeta_{-n} \circ \zeta_{n} (\al_{i} (n)) \\
  &= &\zeta_{-n} (f^{-1} (\al_{i} (n+1)) \\
  &= &f^{-1}( \al_{i} (n+1)) \\
\end{eqnarray*}
\noi and if $M(i) \leq -n \leq -1$,
\begin{eqnarray*}
\rho_{n} (\al_{i} (-n) ) &= &\zeta_{-n} \circ \zeta_{n} (\al_{i} (-n) ) \\
& = &\zeta_{-n}( \al_{i} (-n)) \\
& = & f^{-1} (\al_{i} (-n+1)). 
\end{eqnarray*}

This completes the proof since, for $n=0, \ \rho_{0} = \zeta_{0}$ and this
case had already been taken care of.  $\Box$

\bigskip

\begin{figure}
\begin{center}~
\psfig{file=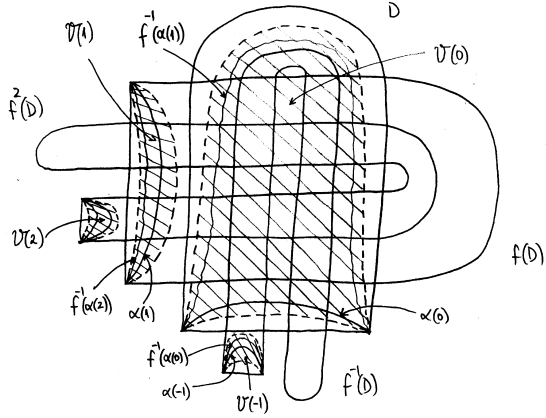,height=3in}
\end{center}
\caption{The first few ${\mc{V}}(n)$'s for a pruning collection 
with only one ($c,e$)-disk $D$.}
\label{f17}
\end{figure}

\begin{cor} \label{44}

The sequence $R_{n} = \dis\bigcup^{n}_{i=0}  r_{n}$ is a Cauchy
sequence in the uniform topology and converges to an isotopy $R: \pi
\times [0,1] \ra \pi$.  If we set $\rho(\cdot) = R (\cdot, 1)$, for each $i
\in \un{L}$ and $M(i) \leq n < N(i)$, $\rho(\al_{i} (n) ) = f^{-1} (\al_{i}
(n+1))$.  Moreover supp $R \subset \bigcup \{ {\mc{V}}_{i} (n); \ i \in 
\un{L}, \ M(i) \leq n < N(i) \}$.
\end{cor}

\noi{\sc Proof}:  
Given $\var > 0$, by Proposition~\ref{43}, there exits $K$ large 
enough so that all the connected components of supp $r_{m}$ have diameter 
smaller than $\var$ if $m \geq K$.  Let $n > m \geq K$.  We 
then have
\begin{eqnarray*}
d (R_{m}, R_{n} ) & =  & \sup_{(x,t)} d (R_{m}(x,t), R_{n} (x,t)) \\
  & = & \sup_{(x,t)} d (R_{m} (x,t), [ R_{m} \cup \dis\bigcup ^{n} 
_{m+1} r_{i} ] (x,t) ) \\
& = & \sup_{(x,t)} d (x, \dis \bigcup^{n}_{m+1} r_{i} (x,t)) \\
& < & \var 
\end{eqnarray*}
\noi where the last inequality is a consequence of Proposition~\ref{36}.  
This shows that $R_n$ is a Cauchy sequence.  The remaining statements 
are readily proven and we leave them to the reader.  $\Box$
\bigskip

\begin{prop} \label{45}
Let $R$ and $\rho$ be as in Corollary~\ref{44}.  Then for each $i \in 
\un{L}$ we have:

\begin{description}
\item[(i)] $\rho (D^{c} (\al_{i}(n))) = D^{c} (f^{-1}(\al_{i}(n+1))) \quad
{\mbox{\rm for }} 1 \leq n < N(i)$ and

\item[(ii)] $\rho (D^{e} ( \al_{i}(m))) = D^{e} (f^{-1} (\al_{i} (m+1))) 
\quad \mbox{\rm for } M(i) \leq m < 0$.
\end{description}
\end{prop}

\noi {\sc Proof}:
Notice that supp $R \subset \bigcup \{ {\mc{V}}_{i}(n); \ i \in \un{L}, \ 
M(i) \leq 
n < N(i) \}$.  By Corollary~\ref{67}, for $n \geq 1, f^{n} (C_{i}) \cap \
\mbox{\rm supp } R = \es$ and by Corollary~\ref{44}, if $1 \leq n < 
N(i), \rho(\al_{i}(n)) = f^{-1} (\al_{i}(n+1))$.  Therefore
\begin{eqnarray*}
\rho(f^{n} (C_{i} ) \cup \al_{i} (n) ) & = & \rho (f^{n} (C_{i}) ) \cup \rho 
(\al_{i} (n)) \\
& = & f^{n} (C_{i}) \cup f^{-1} (\al_{i} (n+1)) 
\end{eqnarray*}
But 
$$f^{n} (C_{i}) \cup \al_{i} (n) = \p D^{c} (\al_{i} (n))$$ 
\noi and 
$$f^{n} 
(C_{i} ) \cup f^{-1} ( \al_{i} (n+1)) = \p D^{c} ( f^{-1} ( \al_{i} 
(n+1)))$$  

This completes the proof of (i).  (ii) is proven analogously.  
$\Box$

\bigskip

\begin{defn}
For each $n \geq 0$, let $\psi_{n} = f \circ \rho_{n}$, $\Psi_{n} = \dis 
\bigcup^{n} _{i=0} \psi_{i}$ and $\Psi = f \circ \rho$, i.e., $\psi_{n} ( 
\cdot) = f \circ r_{n} (\cdot, 1)$, $\Psi_{n} (\cdot ) = f \circ R_{n} ( 
\cdot, 1)$ and $\Psi (\cdot ) = f \circ R (\cdot, 1).$
\end{defn}

Recall that if $\xi: X \ra X$ is a homeomorophism we defined
$$ \mbox{\rm supp } \xi = {\mc{C}} \{ x \in X; \quad \xi(x) = x \}. $$

\begin{lem} \label{68}
Let $\xi, \eta: X \ra X$ be homeomorphisms so that supp $\xi \subset A$ 
and supp $\eta \subset B$.  Then $A \cup B = A \cup \xi \circ \eta (B)$.
\end{lem}

\noi{\sc Proof}:
First notice that if supp $\xi \subset A$ then $\xi(A) = A$ since $\xi 
({\mc{C}}A) = {\mc{C}}A$ and $\xi$ is a homeomorphism.  Therefore, since 
supp $\xi \circ \eta \subset A \cup B$ we have
$$ A \cup B = \xi \circ \eta ( A \cup B) = \xi (A) \cup (B) = A \cup \xi 
(B) = A \cup \xi \circ \eta (B) \quad \Box $$

\begin{prop} \label{46}
For $n \geq 0$,
\begin{description}
\item[(i)]  $f^{n} (P) \cup f^{-n} (P) = \rho _{n} (f^{n} (P) ) \cup 
f^{-n} (P)$;
\item[(ii)] $f^{n} (P) \cup f^{-n} (P) = f^{n} (P) \cup \rho^{-1}_{n} 
(f^{-n}(P) )$.
\end{description}
\end{prop}

\noi{\sc Proof}:
For $n=0$, supp $\rho_{0} \subset P$ and the result follows.  For $n \geq 
1, \rho_{n} = \zeta_{-n} \circ \zeta_{n}$ and supp $\zeta_{-n} \subset 
f^{-n} (P)$ and supp $\zeta_{n} \subset f^{n}(P)$.  The results now follow 
as easy applications of Lemma~\ref{68} (see figure 18.)  $\Box$
\bigskip

\begin{figure}
\begin{center}~
\psfig{file=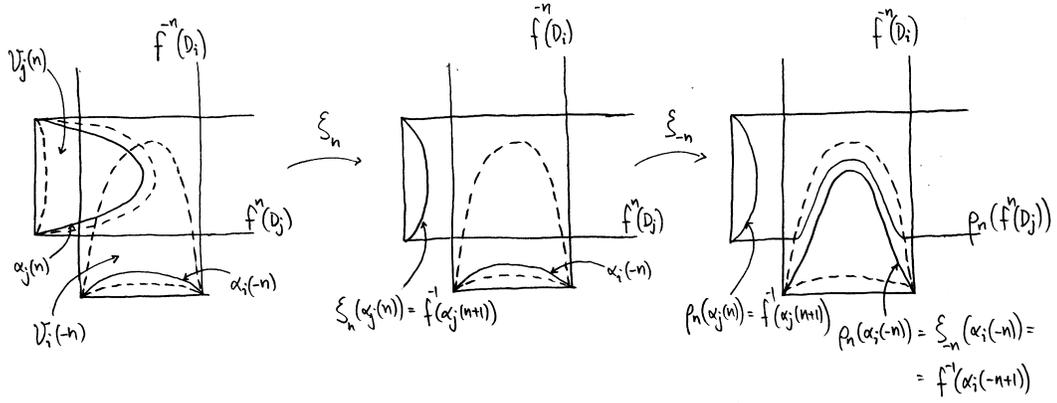,height=2.5in}
\end{center}
\caption{The homeomorphism $\rho_n$.}
\label{f18}
\end{figure}

The following corollary is immediate from the definition of $\psi_{n}$.

\begin{cor} \label{47}
For $n \geq 1$,

\begin{description}
\item[(i)]  $f^{n+1}(P) \cup f^{-n+1} (P) = \psi_{n} (f^{n} (P) ) \cup 
f^{-n+1} (P) $;
\item[(ii)]  $f^{n} (P) \cup f^{-n} (P) = f^{n} (P) \cup \psi^{-1}_{n} 
(f^{-n+1} (P) ). \quad \Box$

\end{description}
\end{cor}

We now state and prove an important technical proposition to be used 
later.  We will use the following
\bigskip

\begin{defn}  
Let $P(0) =P$ and define inductively $P(\!n\!) \!=\! \Psi_{n-1} (\!P 
(\!n-1 \!) \!)$ and $P(-n) = \Psi^{-1}_{n} (P (-n+1))$, for every $n \geq 1$.
\end{defn}

\begin{prop} \label{48}
With the notation above $P(1) = f(P)$ and 

\begin{description}
\item[(i)] for $n \geq 1, \quad \dis\bigcup^{n-1}_{-n+1} f^{k}(P) = \dis 
\bigcup^{n-1}_{-n+1} P(k)$;
\item[(ii)]  for $n \geq 2, \quad \dis \bigcup^{n}_{-n+2} f^{k}(P) = \dis 
\bigcup^{n}_{-n+2} P (k)$.

\end{description}
\end{prop}

\noi{\sc Proof}:
We will use induction on $n$.  For $n=1$, (i) states that $P=P(0)$, which 
is just the definition, whereas $P(1) = \Psi_{0}(P) = \psi_{0}(P) = f 
\rho_{0}(P)$ and, since supp $\rho_{0} \subset P, \ \rho_{0} (P) =P$, which 
shows that $P(1) = f(P)$ (see figure 19.)

We  now show that $\dis \bigcup^{2}_{0} f^{k}(P) = \dis \bigcup^{2}_{0} 
P(k)$, but before we start, let us point out that, from the definitions of 
$\psi_n$ and $\Psi_n$, the following is clear, for each $n \geq 0$:

\begin{description}
\item[a)] $\Psi_{n} = f$ in the complement of supp $R_{n} \subset \dis 
\bigcup^{n}_{-n} f^{k} (P)$;
\item[b)] $\psi_{n} =f$ in the complement of supp $\rho_{n} \subset [ f^{n} 
(P) \cup f^{-n} (P) ] \sm \dis \bigcup^{n-1}_{-n+1} f^{k} ( \ov{P} )$;
\item[c)]  $\Psi_{n} = \left\{ \begin{array}{ll}
              \Psi_{n-1} &\mbox{\rm within } \dis \bigcup^{n-1}_{-n+1} 
f^{k} (P) \\
              \psi_{n} &\mbox{\rm  without } \dis\bigcup^{n-1}_{-n+1} 
f^{k} (P)      \end{array}   \right.; $

\item[d)] $\Psi^{-1}_{n} = \left\{ 
    \begin{array}{ll}
     \Psi^{-1}_{n-1} &\mbox{\rm within } \dis\bigcup^{n}_{-n+2} f^{k}(P) 
= \Psi_{n-1} \left( \dis\bigcup^{n-1}_{-n+1} f^{k} (P) \right) \\
     \psi^{-1}_{n} & \mbox{\rm without } \dis\bigcup^{n}_{-n+2} f^{k} (P)
= \Psi_{n-1} \left( \dis\bigcup^{n-1}_{-n+1} f^{k} (P) \right)    
 \end{array}   \right. .  $
\end{description}

Having said this, let us go back to the proof of $\dis\bigcup^{2}_{0}
f^{k}(P) = \dis\bigcup^{2}_{0}P(k)$.  Notice that from c) above we have
$$ P(2) = \Psi_{1} (P(1)) = \left\{ \begin{array}{ll}
     \Psi_{0}(P(1) ) \quad \mbox{\rm within } & P \\
     \psi_{1} (P(1)) \quad \mbox{\rm without } & P
  \end{array}
   \right.  $$
\noi and, since we have seen that $P(0)=P$ and $P(1) = f(P)$,
\begin{eqnarray*}
\Psi_{1}(P(1))  &= &\Psi_{1} (P(1) \cap P(0) ) \cup \Psi_{1} (P(1) \sm
P(0) ) \\
   & = &[\Psi_{1} (P(1)) \cap \Psi_{0}(P(0)) ] \cup [ \Psi_{1} (f(P) \sm
P ) ] \\
   & = &[ P(2) \cap P(1) ] \cup [ \psi_{1} (f(P) \sm P) ] \\
 & = &[ P (2) \cap f(P) ] \cup [ \psi_{1} (f(P)) \sm f (P) ]
\end{eqnarray*}
\noi where the last equality is a consequence of b) above.  Thus
\begin{eqnarray*}
\dis\bigcup^{2}_{0} P(k) &= &\Psi_{1} (P(1)) \cup \dis\bigcup^{1}_{0}
P(k) \\
  &= &\Psi_{1} (P(1)) \cup \dis\bigcup^{1}_{0} f^{k}(P) \\
  & = &\psi_{1} (f(P)) \cup \dis\bigcup^{1}_{0} f^{k}(P) \\
  & = &\dis\bigcup^{2}_{0} f^{k} (P) \\
\end{eqnarray*}
\noi where the last equality is a consequence of Corollary~\ref{47}
(i), with $n=1$.
\bigskip

We now show that $\dis\bigcup^{1}_{-1}f^{k}(P) = \dis\bigcup^{1}_{-1}
P(k)$.  From d) above we have
$$ P(-1) = \Psi^{-1}_{1}(P(0))= \left\{ \begin{array}{l}
\Psi^{-1}_{0} (P(0)) \quad \mbox{\rm within } f(P)= P(1) \\
        \psi^{-1}_{1}(P(0)) \quad \mbox{\rm without } f (P) = P(1)
        \end{array}
        \right.  $$
\noi so that
\begin{eqnarray*}
\Psi^{-1}_{1} (P(0)) &= &\Psi^{-1}_{1} (P(0) \cap P(1)) \cup
\Psi^{-1}_{1} (P(0) \sm P(1)) \\
  & = & [ \Psi^{-1}_{1} (P(0)) \cap \Psi^{-1}_{0} (P(1)) ] \cup
\Psi^{-1}_{1} ( P \sm f(P))  \\
   & = & [ P(-1) \cap P(0) ] \cup [ \psi^{-1}_{1} (P \sm f(P)) ] \\
   & = & [P(-1) \cap P ] \cup [ \psi^{-1}_{1}(P) \sm P ]
\end{eqnarray*}
\noi where the last equality is a consequence of b).  From this we see that
\begin{eqnarray*}
\dis\bigcup^{1}_{-1} P(k) &= &\dis\bigcup^{1}_{0} P(k) \cup
\Psi^{-1}_{1} (P(0)) \\
      \\
&= &\dis\bigcup^{1}_{0} f^{k}(P) \cup \Psi^{-1}_{1} (P(0))\\
        \\
& = &\dis\bigcup^{1}_{0} f^{k} (P) \cup \psi^{-1}_{1} (P)\\
          \\
  &= &\dis \bigcup^{1}_{-1} f^{k}(P) \\
\end{eqnarray*}
\noi where the last equality is again a consequence of
Corollary~\ref{47} (ii), with $n=1$.  This completes the proof of (i)
and (ii) for $n=2$.  Suppose we have proven that (i) and (ii) hold for $2
\leq n \leq N$.  From this assumption the assertions below follow:

\begin{description}
\item[1)]  $\dis \bigcup^{n}_{-n+1} f^{k}(P) = \dis\bigcup^{n}_{-n+1}
P(k)$, for $2 \leq n \leq N$, by just taking the union of (i) and (ii).

\item[2)] $f^{n}(P) = P(n)$ and $f^{-n+1}(P) = P (-n+1)$ in the
complement of $\dis\bigcup^{n-1}_{-n+2} f^{k}(P) =
\dis\bigcup^{n-1}_{-n+2} P(k)$, for $0 \leq n \leq N$.  This can be seen
as follows: by (i), $\dis\bigcup^{n}_{-n+2} f^{k}(P) = \dis \bigcup^{n}
_{-n+2} P(k)$ and by 1), $\dis \bigcup^{n-1} _{-n+2} f^{k}(P) =
\dis\bigcup^{n-1}_{-n+2} P(k)$.  Then
$$ f^{n}(P) \cup \dis \bigcup^{n-1}_{-n+2} f^{k} (P) = P (n) \cup \dis
\bigcup^{n-1}_{-n+2} P(k)\ . $$
\noi It then follows that $f^{n}(P)=P(n)$ in the complement of $\dis
\bigcup^{n-1}_{-n+2} f^{k}(P) = \dis \bigcup ^{n-1}_{-n+2} P(k)$.  The
other part is proven similarly.

\item[3)]  $\Psi_{n} (P(j))= P(j+1)$ for any $-n \leq j \leq n, \ 0 \leq
n \leq N$.  For notice that $\Psi _{n} = \Psi_{|j|}$ in $\dis
\bigcup^{|j|}_{-|j|}f^{k}(P) = \dis\bigcup^{|j|}_{-|j|} P(k) \supset
P(j)$.  Thus $\Psi_{n}(P(j)) = \Psi_{|j|} (P(j)) = P(j+1)$ from the
definition of $P(j)$.  This reasoning is valid for $-n \leq j \leq n, \ 0
\leq n < N$.  For $n = N$ what remains to be shown is that $\Psi_{N}
(P(N)) = P(N+1)$ and $\Psi_{N} (P(-N)) = P (N+1)$ or equivalently $P(-N)
= \Psi_{N}^{-1} (P(-N+1))$.  But these are just the definitions again.

\end{description}

We now proceed to prove (i) and (ii) for $N+1$.  We start with (ii) $\dis
\bigcup ^{N+1}_{-N+1}P(k) = \dis \bigcup^{N+1}_{-N+1} f^{k}(P)$.

>From c) in the beginning of the proof

$$ P(N+1) = \Psi_{N} (P(N)) = \left\{ \begin{array}{l}
         \Psi_{N-1} (P(N)) \quad \mbox{\rm within } \dis
\bigcup^{N-1}_{-N+1} f^{k}(P)  \\
  \qquad = \dis \bigcup^{N-1}_{-N+1}P(k) \\
     \psi_{N} (P(N)) \quad \mbox{\rm without } \dis \bigcup^{N-1}_{-N+1}
f^{k}(P) \\
   \qquad = \dis \bigcup^{N-1}_{-N+1} P(k) \\
   \end{array}
\right.   $$

Thus,
\begin{eqnarray*}
\Psi_{N} (P(N)) &= &\Psi _{N} \left( P(N) \cap \dis\bigcup^{N-1}_{-N+1} P 
(k) \right) \cup \Psi_{N} \left(P (N)
\sm \dis\bigcup^{N-1}_{-N+1} P(k)\right) \\
& = &\left[ \Psi_{N} (P(N))
\cap \Psi_{N-1} \left( \dis\bigcup^{N-1}_{-N+1}
P(k) \right) \right] \cup \\ 
& & \Psi_{N}  \left(f^{N}(P) \sm 
\dis\bigcup^{N-1}_{-N+1} f^{k}(P) \right) \\
&= &\left[ P (N+1) \cap \dis\bigcup^{N}_{-N+2} P(k) \right] \cup \psi
_{N} \left(f^{N}(P) \sm \dis\bigcup^{N-1}_{-N+1} f^{k}(P) \right) \\
&= &\left[ P (N+1) \cap \dis\bigcup^{N}_{-N+2} f^{k}(P) \right] \cup
\left[ \psi_{N} (f^{N} (P) ) \sm \dis\bigcup^{N}_{-N+2} f^{k} (P) \right]
\end{eqnarray*}
\noi where we used 2) in the second equality, 3) in the third and b) from the
beginning in the forth, not to mention the induction hypothesis here and
there.  From this it follows that
\begin{eqnarray*}
\dis\bigcup^{N+1}_{-N+1} P(k) &= &\Psi_{N} (P(N)) \cup \dis\bigcup^{N}
_{-N+1} P(k) \\
&= & \Psi_{N} (P(N)) \cup \dis\bigcup^{N}_{-N+1} f^{k}(P) \\
&= & \psi_{N} (f^{N} (P)) \cup \dis\bigcup^{N}_{-N+1} f^{k}(P) \\
&= & \dis\bigcup^{N+1}_{-N+1} f^{k} (P)
\end{eqnarray*}
\noi where the last equality comes from Corollary~\ref{47} (i) with $n=N$.

We now prove (i) $\dis\bigcup^{N}_{-N} P(k) = \dis\bigcup^{N}_{-N} 
f^{k}(P)$. From d) we have
$$ P(-N) = \Psi^{-1}_{N} (P(-N+1)) = \left\{ \begin{array}{l}
   \Psi^{-1}_{N-1} (P(-N+1)) \  \mbox{\rm within } \\ 
   \quad \dis\bigcup^{N}_{-N+2}
   f^{k}(P) =  \dis\bigcup^{N}_{-N+2} P(k)   \\
  \psi^{-1}_{N} (P(-N+1)) \ \mbox{\rm without } \\ 
 \quad \dis\bigcup^{N}_{-N+2}
 f^{k} (P) =   \dis\bigcup^{N}_{-N+2} P(k)
\end{array} \right. $$

Thus
\begin{eqnarray*}
\Psi^{-1}_{N} (P(-N+1)) &= &\Psi^{-1}_{N} \left( P(-N+1) \cap \dis\bigcup
^{N}_{-N+2} P(k) \right) \cup \\ 
& & \Psi^{-1}_{N} 
 \left(P (-N+1) \sm \dis\bigcup
^{N}_{_N+2} P(k) \right) \\
&= &\left[ \Psi_{N}^{-1}(P(-N+1)) \cap \Psi^{-1}_{N-1} \left(\dis\bigcup
^{N}_{-N+2} P(k) \right) \right] \cup \\ 
& & \Psi^{-1}_{N} 
 \left( f^{-N+1} (P)
\sm \dis\bigcup ^{N}_{-N+2} f^{k}(P) \right) \\
&= &\left[ P (-N) \cap \dis\bigcup^{N-1}_{-N+1} P(k) \right] \cup \\ 
& & \psi^{-1}_{N} 
 \left( f^{-N+1}(P) \sm \dis\bigcup^{N}_{-N+2} f^{k} (P) 
\right)\\
&= &\left[ P(-N) \cap \dis\bigcup^{N-1}_{-N+1} f^{k}(P) \right] \cup \\
& &  \left[ \psi^{-1}_{N} ( f^{-N+1} (P) ) \sm \dis\bigcup 
^{N-1}_{-N+1} f^{k} (P) \right]
\end{eqnarray*}
\noi where we have used 2) in the second equality, 3) in the third, b) in the
fourth and the induction hypothesis.

Therefore
\begin{eqnarray*}
\dis\bigcup^{N}_{-N} P(k) &= &\dis\bigcup^{N}_{-N+1} P(k) \cup \Psi^{-1}
_{N} (P(-N+1)) \\
&= &\dis\bigcup^{N} _{-N+1} f^{k}(P) \cup \Psi^{-1}_{N} (P(-N+1)) \\
&= &\dis\bigcup^{N}_{-N+1} f^{k} (P) \cup \psi^{-1}_{N} (P(-N+1)) \\
&= &\dis\bigcup^{N}_{-N} f^{k}(P)
\end{eqnarray*}

\noi where the last equality comes from Corollary~\ref{47} (ii) with $n=N$.
This completes the proof.  $\Box$
\bigskip

\begin{figure}
\begin{center}~
\psfig{file=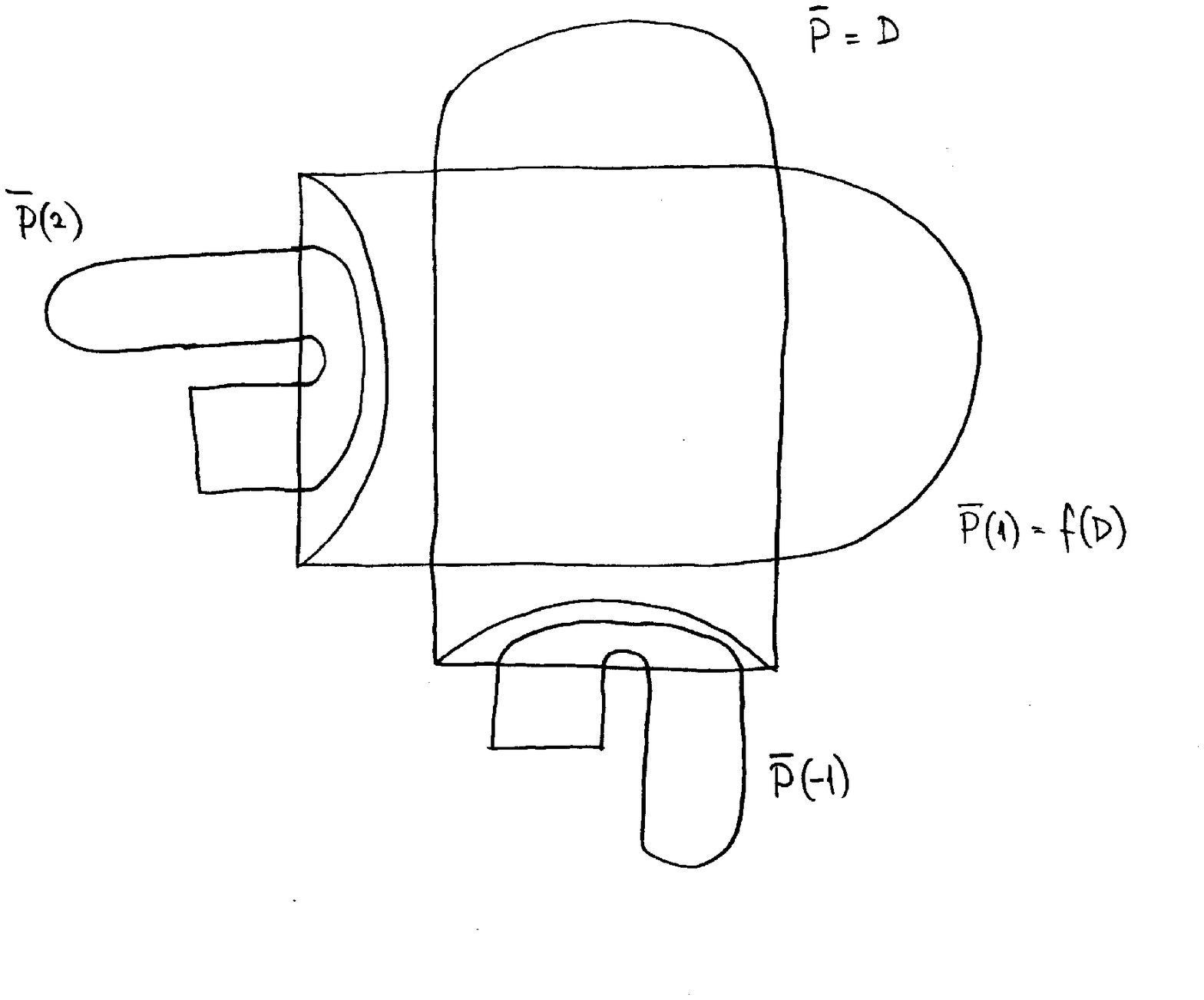,height=3.25in}
\end{center}
\caption{$\ov{P}(k), k=- 1, 0,1,2$, for a pruning collection 
containing only one ($c,e$)-disk $D$.}
\label{f19}
\end{figure}

\begin{cor}\label{49}
For $n \geq 1$

\begin{description}
\item[(i)]  $\dis\bigcup^{n}_{n+1} f^{k}(P) = \dis\bigcup^{n}_{-n+1}P(k)$;
\item[(ii)] $f^{n}(P) = P(n)$ and $f^{-n+1}(P) = P(-n+1)$ in the complement
of $$\dis\bigcup ^{n-1}_{-n+2} f^{k}(P) = \dis\bigcup^{n-1}_{-n+2} P(k).$$
\end{description}
\end{cor}

\noi{\sc Proof:}  The proof is the same as that given for 1) and 2) in the
proof of Proposition~\ref{48}.  $\Box$
\bigskip

\begin{cor}\label{50}
If $\Psi$ is as we defined above, $P(k) = \Psi(P(k-1))$ for every $k \in
{\Bbb{Z}}$, that is $\{P(k); \ k \in {\Bbb{Z}} \}$ is an orbit under $\Psi$.
\end{cor}

\noi{\sc Proof:}  Just notice that $\Psi = \Psi_{n}$ in
$\dis\bigcup^{n}_{-n} f^{k}(P)$ and argue like in the proof of 3) in
Proposition~\ref{48}.  $\Box$
\bigskip

We are now going to define new closed disks $A_{i}, \ i \in \un{L}$ whose
union is still the closed pruning front $\ov{P}$.  We will see that the
cross-cut $\al_{i}(0) \subset D_{i}$ is also a cross-cut in $A_i$ and
divides it into two disks $A^{c}_{i}$ and $A^{e}_{i}$ (see figure 
20.)  These will have
some disjoint/nested properties we will make precise later and will be
useful in the proof of the theorem.
\bigskip

\begin{defn}
Let $A_{L} = A_{L}(0)= D_{L}$ and, for $1 \leq i \leq L$, set $A_{i} =
A_{i}(0) = \zeta^{-1} _{L,0} \circ \ldots \circ
\zeta^{-1}_{i+1,0}(D_{i})$.  Then define inductively for $n \geq 1$, $ \
A_{i}  (n) = \Psi (A_{i} (n-1))$ and $A_{i}(-n) = \Psi^{-1} (A_{i}(-n+1))$.
\end{defn}

\begin{figure}
\begin{center}~
\psfig{file=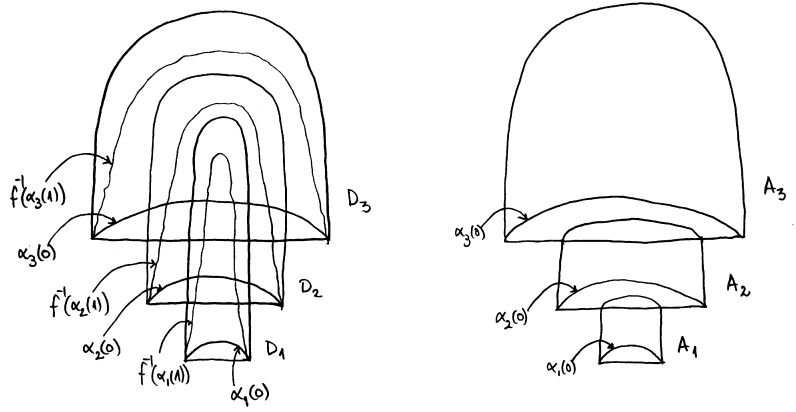,height=2.75in}
\end{center}
\caption{The $D_{i}$'s and the $A_i$'s}
\label{f20}
\end{figure}

\begin{prop}\label{51}
For $l \in \un{L}, \ \dis \bigcup^{L}_{i=l} A_{i} = \dis\bigcup^{L}
_{i=l} D_{i}$.  In particular $\dis\bigcup^{L} _{i=1} A_{i} = \ov{P}$.
\end{prop}

\noi{\sc Proof:}  By definition $A_{L} = D_{L}$.  Assume we have shown
that $\dis \bigcup^{L} _{i=l+1} A_{i} = \dis \bigcup ^{L} _{i=l+1}
D_{i}$.  Then
\begin{eqnarray*}
\dis\bigcup^{L} _{i=l} A_{i} & = & \dis \bigcup^{L} _{i=l+1} A_{i} \cup
A_{l} \\
& = & \dis\bigcup^{L}_{i=l+1} D_{i} \cup \zeta^{-1} _{L,0} \circ \ldots
\circ \zeta^{-1} _{l+1,0} (D_{l}) \\
& = & \dis\bigcup^{L}_{i=l} D_{i}
\end{eqnarray*}

\noi where the last equality holds because supp $\zeta^{-1}_{L,0} \circ
\ldots \circ \zeta^{-1}_{l+1,0} \subset \dis\bigcup^{L}_{i=l+1}D_{i}. \ \
\Box$

\bigskip

\begin{cor} \label{53}
For every $n \in {\Bbb{Z}}, \ \ov{P}(n) = \dis \bigcup^{L} _{i=1} A_{i}(n). \
\  \Box$
\end{cor}

\begin{prop} \label{52}
For $n \geq 1 $ and $i \in \un{L}$,

\begin{description}
\item[(i)]  $\zeta_{n} \left( \dis\bigcup _{j \leq i} f^{n} (D_{j})
\right) = \dis\bigcup_{j \leq i} f^{n} (D_{j})$;

\item[(ii)]  $ \zeta_{-n} \left( \dis \bigcup_{j \geq i} f^{-n}(D_{j})
\right) = \dis \bigcup _{j \geq i} f^{-n} (D_{j})$.
\end{description}

\end{prop}

\noi{\sc Proof:}  If $k > i \geq j$ then $D_{k} \not\prec D_j$ and we
have seen that for $n \geq 1$,  $I^{c} ( \gamma_{k} (n) )$ is either
contained in $f^{n}(I_{j})$ or it is disjoint from $f^{n}(D_{j})$.  If
$I^{c} (\gamma_{k} (n))$ is contained in $f^{n}(I_{j})$ it is because
$f^{n}(D_{k}), \ f^{n}(D_{j} ) | _{f^{n}(C_{k}) }$ and therefore
$f(D_{k}), \ f(D_{j}) |_{f(C_{k})}$.  By Proposition~\ref{42}, $N(i)
=1$ and it follows that $k_{i,n} \equiv$ identity.  If $I^{c} (\gamma_{k}
(n))$ is disjoint from $f^{n}(D_{j})$ so is ${\mc{V}}_{k}(n)$, since 
${\mc{V}}_{k}(n) \subset I^{c}( \gamma_{k} (n) )$.  Either way we see that 
(supp $\zeta_{k,n} ) \cap f^{n}(D_{j} ) = \es$.  Thus
\begin{eqnarray*}
\zeta_{n} \left( \dis \bigcup _{j \leq i} f^{n} (D_{j} ) \right) & = &
\zeta _{1,n} \circ \ldots \circ \zeta_{L,n} \left( \dis \bigcup _{j \leq
i} f^{n} (D_{j} ) \right) \\
& = & \zeta_{1,n} \circ \ldots \circ \zeta_{i,n} \left( \dis\bigcup _{j
\leq i} f^{n} (D_{j}) \right) \\
& = & \dis \bigcup _{j \leq i} f^{n} (D_{j})
\end{eqnarray*}
\noi where the last equality holds because supp $\zeta_{1,n} \circ \ldots
\circ \zeta_{i,n} \subset \dis \bigcup _{j \leq i} f^{n} (D_{i})$.  This
proves (i).  (ii) is proven analogously.  $\Box$
\bigskip

The next proposition and corollary are analogous to
Proposition~\ref{46} and Corollary~\ref{47}.  The proofs use
Proposition~\ref{52} but are otherwise completely similar.  We omit them.

\begin{prop} \label{54}
For $n \geq 1$ and $i \in \un{L}$,

\begin{description}
\item[(i)]  $\rho_{n} \left( \dis \bigcup _{j \leq i} f^{n} (D_{j}) \right)
\cup f^{-n} (P) = \dis \bigcup _{j \leq i} f^{n} (D_{j}) \cup f^{-n} (P) $;
\item[(ii)]  $ f^{n}(P) \cup \rho^{-1}_{n} \left( \dis \bigcup _{j\geq i}
f^{-n} (D_{j} ) \right) = f^{n} (P) \cup \dis \bigcup _{j \geq i} f^{-n}
(D_{j})$.  $\Box$
\end{description}

\end{prop}

\begin{cor} \label{55}
For $n \geq 1$ and $i \in \un{L}$,

\begin{description}
\item[(i)]  $\psi _{n} \left( \dis \bigcup_{j \leq i} f^{n} (D_{j} )
\right) \cup f^{-n+1} (P) = \dis \bigcup _{j \leq i} f^{n+1} (D_{j} )
\cup f^{-n+1} (P)$;

\item[(ii)]  $f^{n}(P) \cup \psi^{-1}_{n} \left( \dis \bigcup _{j \geq i}
f^{-n+1} (D_{j} ) \right) = f^{n} (P) \cup \dis \bigcup _{j \geq i}
f^{-n} (D_{j} )$.  $\Box$
\end{description}

\end{cor}

We can now state and prove a proposition which sharpens
Proposition~\ref{48} somewhat.  Although the proof goes along the same
lines as that of Proposition~\ref{48} we present it for completeness.
\bigskip

\begin{prop} \label{56}
For $n \geq 1$ and $i \in \un{L}$ we have
\begin{description}

\item[(i)]  $\dis \bigcup_{j \leq i} f^{n} (D_{j}) \cup \dis
\bigcup^{n-1} _{-n+2} f^{k} (P) = \dis \bigcup _{j \leq i} A_{j} (n) \cup
\dis \bigcup ^{n-1} _{ -n+2} P(k)$;

\item[(ii)]  $ \dis \bigcup_{j \geq i} f^{-n+1} (D_{j}) \cup \dis \bigcup
^{n-1}_{-n+2} f^{k} (P) = \dis \bigcup _{j \geq i} A_{j} (-n+1) \cup \dis
\bigcup^{n-1}_{-n+2} P (k).$

\end{description}

\end{prop}

\noi {\sc Proof:}
The proof is by induction on $n$.  Notice that for $n=1$, (ii) above is
just Corollary~\ref{53}.  In order to prove (i) with $n=1$, observe
that, since $A_{i} = A_{i}(0) \subset \ov{P}, \ A_{i}(1) =
\Psi(A_{i}(0)) = \psi_{0} (A_{i}(0)) = f \circ \rho_{0} (A _{i}(0))$.  Thus
\begin{eqnarray*}
A_{i}(1) & = & f \circ \rho_{0} ( \zeta^{-1} _{L,0} \circ \ldots \circ
\zeta^{-1}_{i+1,0} (D_{i}) ) \\
& = & f \circ ( \zeta_{1,0} \circ \ldots \circ \zeta_{L,0} ) \circ (
\zeta^{-1}_{L,0} \circ \ldots \circ \zeta ^{-1} _{i+1,0} (D_{i} )) \\
& = & f \circ ( \zeta_{1,0} \circ \ldots \circ \zeta_{i,0}) (D_{i})
\end{eqnarray*}

Reasoning as in the proof of Corollary~\ref{53}, it is
easy to prove that $\dis \bigcup_{j \leq i} A_{j} (1) = \dis \bigcup_{j
\leq i} f(D_{j} )$ which is (i) for $n=1$.

Assume we have shown that 
$$\dis \bigcup_{j \leq i} A_{j} (n) \cup \dis
\bigcup ^{n-1}_{-n+2} P(k) = \dis \bigcup _{j \leq i} f^{n} (D_{j}) \cup
\dis \bigcup^{n-1}_{-n+2} f^{k} (P) . $$  
\noi Then, since we know that $\dis
\bigcup^{n-1}_{-n+1} P(k) = \dis \bigcup^{n-1}_{-n+1} f^{k} (P) $ and
$\dis \bigcup^{n}_{-n+1} P(k) = \dis \bigcup^{n}_{-n+1}$ $f^{k}(P)$ by
Proposition~\ref{48}, just like in the proof of that proposition we
can see that $$\dis \bigcup_{j \leq i} A_{j}(n) \sm \dis
\bigcup^{n-1}_{-n+1} P(k) = \dis \bigcup_{j \leq i} f^{n}(D_{j}) \sm \dis
\bigcup ^{n-1}_{-n+1}f^{k}(P).$$  Using all this information we have
\begin{eqnarray*}
\dis\bigcup_{j \leq i} A_{j} (n+1) \cup \dis \bigcup ^{n} _{-n+1} P(k) &
= & \left[ \dis \bigcup_{j \leq i} A_{j} (n+1) \sm \dis \bigcup^{n}
_{-n+2} P(k) \right] \! \cup \!\dis \bigcup^{n} _{-n+1}\! P(k) \\
& = & \left[ \Psi \left( \dis \bigcup _{j \leq i} A_{j} (n) \sm \dis
\bigcup^{n-1}_{-n+1} P(k) \right) \right] \cup \dis \bigcup ^{n} _{-n+1}
P(k) \\
& = & \left[ \Psi \left( \dis \bigcup_{j \leq i} f^{n} (D_{j} ) \sm \dis
\bigcup^{n-1} _{-n+1} f^{k} (P) \right) \right] \cup \dis\bigcup
^{n}_{-n+1} f^{k} (P) \\
& = & \left[ \psi_{n} \left( \dis \bigcup _{j \leq i} f^{n} (D_{j})
 \sm \dis \bigcup^{n-1} _{-n+1} f^{k} (P) \right) \right] \cup \dis\bigcup
^{n} _{-n+1} f^{k} (P) \\
& = & \left[ \psi_{n} \left( \dis \bigcup _{j \leq i} f^{n} (D_{j})
\right) \sm \dis \bigcup^{n} _{-n+2} f^{k} (P) \right] \cup \dis\bigcup
^{n} _{-n+1} f^{k} (P) \\
& = & \psi _{n} \left( \dis \bigcup_{j\leq i} f^{n} (D_{j}) \right)\! \cup \!
\dis \bigcup ^{n} _{-n+1} f^{k} (P) \\
& = & \dis \bigcup_{j \leq i} f^{n+1} (D_{j} ) \cup \dis \bigcup ^{n}
_{-n+1} f^{k}(P)
\end{eqnarray*}
\noi where the last equality holds by Corollary~\ref{55}, (i).
Statement (ii) is proven analogously and we leave it to the interested
reader.  $\Box$
\bigskip

\begin{cor} \label{57}
For $n \geq 1$ and $i \in \un{L}$ we have:

\begin{description}
\item[(i)] $A_{i} (n) = f^{n} (D_{i})$ in the complement of 
$$\dis
\bigcup _{j < i} f^{n}(D_{j}) \cup  \dis \bigcup^{n-1}_{-n+2} f^{k}(P) =
\dis \bigcup _{j < i} A_{j} (n) \cup \dis \bigcup^{n-1}_{-n+2} P(k) ; $$
\item[(ii)]  $A_{i} (-n+1) = f^{-n +1} (D_{i})$ in the complement of 
$$\dis
\bigcup _{j > i} f^{-n+1} (D_{j}) \cup \dis \bigcup ^{n-1} _{-n+2} f^{k}
(P) = \dis \bigcup_{j >i} A_{j} (-n+1) \cup \dis \bigcup^{n-1}_{-n+2}
 P (k) . \ \Box$$
\end{description}

\end{cor}

\begin{prop} \label{58}
For each $i \in \un{L}, \ \al_{i}(0)$ is a cross-cut in $A_i$ and divides
$A_i$ into two closed disks $A^{c}_{i}$ and $A^{e}_{i}$ bounded by
$\rho^{-1}_{0} ( C_{i}) \cup \al_{i} (0)$ and $\al_{i} (0) \cup E_{i}$,
respectively.
\end{prop}

\noi {\sc Proof:}  In the proof of Proposition~\ref{43} (iii),
we have shown that for each $i \in \un{L}$,  $\rho_{0} (\al_{i}(0)) =
\zeta_{0} (\al_{i}(0) ) = \zeta_{1,0} \circ \ldots \circ \zeta_{i, 0}
(\al_{i}(0))$.  Thus we see that
\begin{eqnarray*}
\al_{i}(0) & = &\rho^{-1}_{0} ( \zeta_{1,0} \circ \dots \circ \zeta_{i,0}
(\al_{i}(0))) \\
&= & \zeta^{-1}_{L,0} \circ \ldots \circ \zeta^{-1}_{1,0} \circ \zeta_{1,0}
\circ \dots \circ \zeta_{i,0}(\al_{i}(0)) \\
&= & \zeta^{-1}_{L,0} \circ \dots \circ \zeta^{-1}_{i+1} (\al_{i}(0))
\end{eqnarray*}
\noi so that $\al_{i}(0)$ is left fixed by $\zeta^{-1}_{L,0} \circ \dots
\circ \zeta^{-1}_{i+1,0}$.  Since $\al_{i}(0)$ is a cross-cut in $D_{i}$ and
$A_{i} = \zeta^{-1}_{L,0} \circ \ldots \circ \zeta^{-1}_{i+1,0} (D_{i}), \
\al_{i}(0)$ is also a cross-cut in $A_i$.

We now show that $A_i$ is bounded by $\rho^{-1}_{0} (C_{i}) \cup E_{i}$ which
will complete the proof of the proposition.  Notice that if $j \leq i$, 
$C_{i} \cap
I_{j} = \es$ so that $\rho^{-1}_{0} (C_{i}) = \zeta^{-1}_{L,0} \circ \ldots
\circ \zeta^{-1}_{i+1,0} (C_{i} )$.  On the other hand, for $j \geq i$, 
$I_{j}
\cap E_{i} = \es$ so that $\zeta ^{-1}_{L,0} \circ \ldots \circ \zeta^{-1}
_{i+1,0}(E_{i}) = E_{i}$.  This shows that $\zeta^{-1}_{L,0} \circ \ldots 
\circ
\zeta^{-1}_{i+1,0} (C_{i} \cup E_{i} ) = \rho^{-1}_{0} (C_{i}) \cup E_i$, as
we wanted.  $\Box$
\bigskip

\begin{prop} \label{59}
For each $i \in \un{L}$, (i) and (ii) hold:

\begin{description}
\item[(i)]  for $n \geq 1, \ A_{i}(n)$ is bounded by the Jordan curve 
$$f^{n}
(C_{i}) \cup \Psi^{n} (E_{i} );$$
\item[(ii)]  for $m \leq 0, \ A_{i} (m)$ is bounded by the Jordan curve  
$$\Psi
^{m} (\rho^{-1}_{0} (C_{i})) \cup f^{m}(E_{i}) = \Psi^{m-1} (f(C_{i})) 
\cup f^ {m}(E_{i}).$$
\end{description}

\end{prop}

\noi{\sc Proof:}
In the proof of Proposition~\ref{58} we saw that $A_{i}(0)$ is bounded by
$\rho_{0}^{-1}(C_{i}) \cup E_{i} = \psi^{-1}_{0} (f(C_{i} )) \cup E_{i} =
\Psi ^{-1} (f (C_{i}) ) \cup E_{i}$.
Therefore $A_{i}(1) = \Psi (A_{i}(0))$ is bounded by $\Psi ( \Psi^{-1} (f(C_{i}
)) \cup E_{i} ) = f(C_{i}) \cup \Psi (E_{i})$, which proves (i) and (ii) for
$n=1$ and $m=0$ respectively.  The general result is now proved by induction
using Corollary~\ref{67} to guarantee that $\Psi^{n}(f(C_{i}) ) = f^{n+1}
(C_{i})$ for $n \geq 0$ and that $\Psi^{m}(E_{i}) = f^{m} (E_{i})$ for $m 
\leq 0. \ \ \Box$
\bigskip

\begin{defn}
Let $A^{c}_{i}(0) = A ^{c}_{i}$ and $A^{e}_{i}(0) = A^{e}_{i} $ as in
Proposition~\ref{58} and define inductively for $n \geq 1$, $A^{c(e)}_{i} (n)
= \Psi (A^{c(e)}_{i} (n-1))$ and $A^{c(e)}_{i} (-n) = \Psi^{-1} (A^{c(e)}_{i}
(-n+1))$ (see figure 21.)
\end{defn}

\begin{figure}
\begin{center}~
\psfig{file=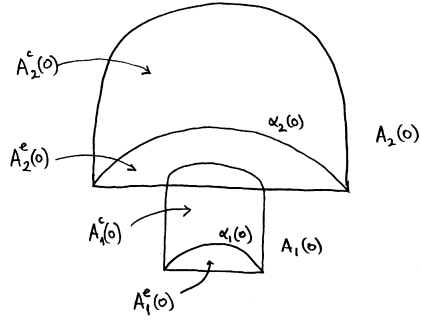,height=2.25in}
\end{center}
\caption{$A^{c}_{i}(0)$ and $A^{e}_{i}(0)$ for $i= 1,2$.}
\label{f21}
\end{figure}

\begin{prop} \label{60}
With the notation just introduced we have:

\begin{description}
\item[(i)]  $A^{c}_{i}(n) = D^{c} (\al_{i}(n))$ for $1 \leq n \leq N(i)$;
\item[(ii)]  $ A^{e}_{i}(m) = D^{e}(\al_{i}(m))$ for $M(i) \leq m \leq0$.
\end{description}
\end{prop}

\noi{\sc Proof:}
That $A^{e}_{i} = D^{e}( \al_{i}( 0))$ is a direct consequence of Proposition
~\ref{58}, since $A^{e}_{i}$ is bounded by $\al_{i}(0) \cup E_{i}$ which is the
same curve that bounds $D^{e} (\al_{i}(0))$.  On the other hand, $A^{c}_{i}(0)$
is bounded by $\rho^{-1}_{0}(C_{i}) \cup \al_{i} (0)$ and it follows that
$A^{c}_{i} (1) = \Psi (A^{c}_{i} (0))$ is bounded by $$\Psi (\rho^{-1}_{0}
(C_{i}) \cup \al_{i}(0)) = \psi_{0} ( \rho^{-1}_{0} (C_{i}) \cup \al_{i}(0))
= f(C_{i}) \cup \al_{i} (1).$$  This shows that $A^{c}_{i}(1) = D^{c}(\al_{i}
(1))$.

Assume we have shown that $A^{c}_{i}(n) = D^{c}(\al_{i}(n))$ for $n < N(i)$.
 Then, using Propositon~\ref{45} (i), we see that
\begin{eqnarray*}
A^{c}_{i} (n+1) & = & \Psi (A^{c}_{i}(n)) \\
&= & f\rho (D^{c} (\al_{i} (n))) \\
&= &f (D^{c} (f^{-1} (\al_{i} (n+1)))) \\
&= &D^{c} (\al_{i} (n+1))
\end{eqnarray*}
\noi This proves (i).  (ii) is proven similarly.  $\Box$
\bigskip

Let $n(i), N(i), m(i)$ and $M(i)$ be as we defined them in the beginning of
this section.  By Proposition~\ref{66}
if $n(i), m(i)$ are finite then $n(i) = 2N(i) - \delta, m(i) = 2M(i) - \delta '
$, where $\delta, \delta ' = 0 $ or $1$.  Moreover, $$f^{N(i)} (D_{i}), 
\ f^{-N(i) + \delta +1} (D_{j}) | _{f^{N(i)} (C_{i}) }$$  and 
$$f^{M(i)}(D_{i}), \ f^{-M(i) +
\delta '} (D_{l} ) |_{f^{M(i)} (E_{i}) }$$
 for some $j, l \in \un{L}$.  
Recall
also that if $D_{1}, D_{2}|_{L}$ and $D_{1} \backslash L \subset I_{2}$ we write
$D_{1} \subset D_{2} |_{L}$.  We can now state

\begin{prop} \label{61}
With the above notation for each $i \in \un{L}$, (i) and (ii) hold:

\begin{description}
\item[(i)]  If $n(i) < \infty$ and $j \in \un{L}$ is largest such that 
$$f^{N(i)}
(D_{i}), f^{-N(i)+ \delta+1 } (D_{j} ) |_{f^{N(i)} (C_{i}) }$$ 
then 
$$A^{c}_{i}
(N(i)) \subset A_{j} (-N(i) + \delta +1) | _{f^{N(i)} (C_{i})};$$

\item[(ii)]  if $m(i) > - \infty$ and $j \in \un{L}$ is smallest such that 
$$f^{M(i)
} (D_{i}), f^{-M(i)+ \delta '} (D_{j}) |_{f^{M(i)} (E_{i}) }$$  
then 
$$A^{e}_{i}
( M(i)) \subset A_{j} (-M(i) + \delta ' ) | _{f ^{M(i)} (E_{i}) }. $$
\end{description}
\end{prop}

\noi {\sc Proof:}
By Proposition~\ref{60} we know that $A^{c}_{i}(N(i)) = D^{c}( \al_{i}
(N(i)))$.  We have to show that $f^{N(i)} (C_{i}) \subset \p A _{j} ( -N(i) +
\delta +1)$ and that 
\begin{eqnarray*}
 A^{c}_{i} (N(i) ) \sm f^{N(i)} (C_{i}) &= &D^{c} 
(\al_{i} (N(i))) \sm  f^{N(i)} (C_{i}) \\ 
 & = &I^{c} (\al_{i} (N(i) ))
 \cup \al_{i} (N(i))  \\ 
& \subset &I (A_{j} (-N(i) + \delta +1))
\end{eqnarray*} 
where the last set is the
interior of $A_{j} (-N(i) + \delta + 1 )$ and the second equality is 
just the definition and only the last inclusion needs proof (see figure 
22.)

Let us first show that 
$$f^{N(i)} (C_{i}) \subset \p A_{j} ( -N(i) + 
\delta +1 ). $$
 By assumption 
$$f^{N(i)} (D_{i}), f^{-N(i) + \delta +1 } (D_{j}) | 
_{f^{N(i)} (C_{i})}$$ 
which implies that 
$$f^{N(i)} (C_{i}) \subset f^{-N(i) + \delta 
+ 1} (C_{j}).$$  
If $n(i) =1$, then $N(i)=1$ and $\delta=1$, so that by 
Proposition~
\ref{59} and the above we have $f(C_{i}) \subset f (C_{j}) \subset \p A_{j}
(1)$.  If $n(i) > 1$, applying $f^{N(i) - \delta -1}$ to the inclusion 
$f^{N(i)} (C_{j}) \subset f^{-N(i) + \delta +1}(C_{j})$ we get
$f^{n(i)-1} (C_{i}) \subset C_{j}$.  From
Corollary~\ref{67} we know that $f^{n} (C_{i}) \cap$ supp $R = \es$ for
$n \geq 1$.  In particular, $f^{n(i) -1}(C_{i}) \cap$ supp $\rho_{0} = 
\es$, and
it follows that $f^{n(i)-1} (C_{i} ) \subset \rho_{0}^{-1} (C_{j}) \subset \p A
_{j} (0)$, this last inculsion coming from Proposition~\ref{59}.  Also, if
$k < n(i)$ then $\Psi^{-k} (f^{n(i)} (C_{i})) = f^{-k} (f^{n(i)}(C_{i}) ) =
f^{n(i)-k} (C_{i})$, so that
\begin{eqnarray*}
 \Psi^{-N(i) + \delta + 1} (f^{n(i) -1} (C_{i})) & = &f^{N(i)} (C_{i}) \\
& \subset
&\Psi^{-N(i) +\delta + 1} ( \rho_{0}^{-1} ( C_{i}) )\\ 
&\subset & \p A _{j} ( -N(i) + \delta +1 )
\end{eqnarray*}
\noi as we wanted.

In order to see that $I^{c} (\al_{i} ( N(i))) \cup \al_{i}(N(i)) \subset I 
(A_{j}
( -N(i) + \delta +1))$ first notice that, since $ \{ \al_{l} (N(i) ) 
\}^{L}_{l=1}$ is\
a ($\var_{N(i)},c$)-collection compatible with $$\{ (f^{k} (D_{j} ), 
\al_{j} (k)
); \ j \in \un{L}, \ -N(i) + 1 \leq k \leq N(i) -1 \}$$ 
and that
$$ f^{N(i)} (D_{i} ),
f^{-N(i) + \delta + 1} (D_{j}) |_{f^{N(i)} (C_{i})}$$
then $$ I^{c} (\al_{i} 
(N(i))) \cup \al_{i} (N(i))  \subset f^{-N(i) + \delta +1} (I_{j}).$$

We will now show that
$$[I^{c} (\al_{i}(N(i))) \cup \al_{i} (N(i)) ] \cap \left[ \bigcup_{l > j} 
f^{-N
(i) + \delta +1} (D_{l} ) \cup \bigcup^{N(i) - \delta -1} _{-N(i) + \delta +2}
f^{k} (\ov{P} ) \right] = \es .$$
\noi This is so because by assumption $f^{N(i)} (D_{i}), f^{-N(i) + \delta + 1}
(D_{l}) \not|_{f^{N(i)} (C_{i})}$ for $l >j$ and from Proposition~\ref{66}
(ii), $f^{N(i)} (D_{i}), f^{k} (D_{j} ) \not|_{f^{N(i)} (C_{i}) }$ for any
$j \in\un{L}$ and $-N(i) + \delta +2 \leq k \leq N(i) -1$.
This together with the aforementioned compatiblity of $\{ \al_{i} (N(i)) \} ^{L}
_{l=1}$ are exactly what we need in order to verify the equation above.  By
Corollary~\ref{57} (ii), 
$$A_{j} (-N(i) +\delta +1) = f^{-N(i) + \delta +1} (D_{j})$$ 
in the complement of
$$ \bigcup_{l > j} f^{-N(i) + \delta +1} (D_{l}) \cup \bigcup^{N(i) - 
\delta -1} _{-N(i) + \delta +2 } f^{k} (P) $$
\noi which shows that $$I^{c} (\al_{i} (N(i))) \cup \al_{i}(N(i))  \subset
A_{j} (-N (i) + \delta + 1).$$  
We leave it for the reader to show that 
it is possible
to put $I(A _{j} ( -N(i) + \delta + 1))$ in place of $A_{j}( -N(i) + \delta +1 )
$ in the inclusion above.  $\Box$
\bigskip

\begin{figure}
\begin{center}~
\psfig{file=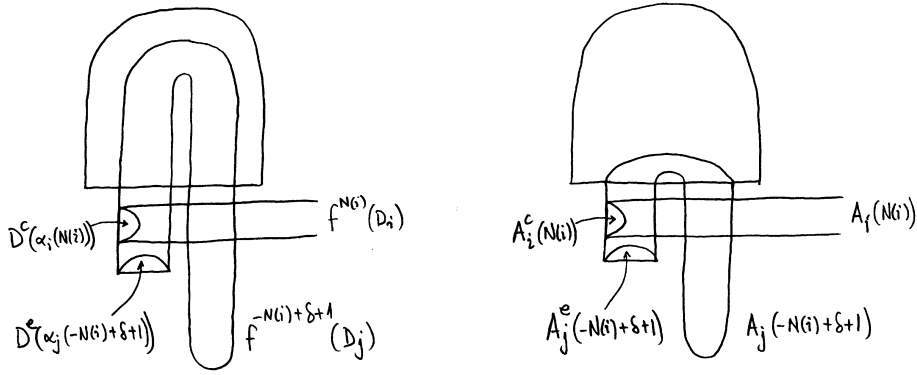,height=2.5in}
\end{center}
\caption{A possible configuration for $f^{N(i)} (D_{i}), \  
f^{-N(i) + \delta + 1} (D_{j})$ and $A_{i} (N(i))$, \ $A_{j}(-N(i) + 
\delta + 1 )$ and $A^{c}_{i} (N(i))$.}
\label{f22}
\end{figure}

\begin{prop} \label{62}
Under the hypotheses of Proposition~\ref{61} (i) and (ii) respectively, (i$'$)
and (ii$'$) below hold:

\begin{description}
\item[(i$'$) ] $A^{c}_{i}(N(i)) \subset A^{c}_{j} ( -N(i) + \delta + 1) | _{f^
{N(i)} (C_{i} )}$;

\item[(ii$'$) ]  $A^{e}_{i} (M(i)) \subset A_{j} (-M(i) + \delta ') | 
_{f^{M(i)} (E_{i}) }$.

\end{description}
\end{prop}

\noi{\sc Proof:}
By Proposition~\ref{61}, $A^{c}_{i} (N(i) ) \subset A_{j} ( -N(i) + 
\delta +1)
|_{f^{N(i)} (C_{i} ) }$.  Therefore all we need to prove is that 
$$[I^{c} ( \al
_{i} (N(i) )) \cup \al_{i} (N(i) ) ] \cap A^{e} _{j} ( -N(i)  + \delta +1 ) 
= \es.$$

There are two cases to be considered:  $M(j) \leq -N(i) +\delta +1$ and $M(j)
> -N(i) + \delta +1$.  If $M(j) \leq -N(i) + \delta + 1$, by Proposition~\ref
{60}, 
$$A^{e}_{j} (-N(i) + \delta +1 ) = D^{e} ( \al_{j} ( -N(i) + \delta 
+1 ))$$
and, since $\{ \al_{l} (N(i)) \}^{L}_{l=1}$ is a ($\var_{N(i)},c$)-collection
compatible with 
$$\{ (f^{k} (D_{j}), \al_{j} (k) ); j \in \un{L}, \ -N(i)+1
\leq k \leq N(i) -1 \},$$
then 
$$[I^{c} (\al_{i}(N(i))) \cup \al_{i} (N(i)) ] \subset
I^{c} ( \al_{j} (-N(i) + \delta + 1 ))$$ 
so that 
$$[ I^{c} (\al_{i} (N(i) )) \cup
\al_{i} (N(i)) ] \cap D^{e} (-N(i) + \delta + 1 ) = \es,$$ 
as we wanted.

If $M_{j} > -N(i) + \delta +1$, there exists $l \in \un{L}$ such that
$$f^{M(j)} (D_{j}), f^{-M(j) + \delta '}(D_{l}) | _ {f^{M(j)} (E_{j}) }$$
where $m(j) = 2 M(j)
+ \delta '$, and, assuming $l$ is the smallest such, by Proposition~\ref
{62} (ii), we can conclude that $A^{e}_{j}(M(j)) \subset A_{l} (-M(j) + \delta
')$.  Therefore
\begin{eqnarray*}
A^{e}_{j} (-N(i) + \delta +1) & = & \Psi ^{-M(j) -N(i) + \delta +1 }
(A^{e}_{j} (M(j))) \\
& \subset & \Psi^{-M(j) -N(i) + \delta +1}
(A_{l} ( -M(j) + \delta ' ) ) \\
& = & A_{l} (-m(j) -N(i) + \delta +1 ).
\end{eqnarray*}

>From $M(j) \geq -N(i) + \delta +2$ we have
\begin{eqnarray*}
-m(j) -N(i) + \delta + 1 & = & -2 M(j) + 2 \delta ' -N(i) + \delta + 1 \\
& \leq & 2N(i) - 2 \delta -4 +2 \delta ' -N(i) + \delta + 1 \\
& = & N(i) - \delta -3 -2 \delta ' \\
& \leq & N(i) - \delta -1
\end{eqnarray*}

\noi If $m(j) > 0$, then $-m(j) -N(i) + \delta +1 \geq -N(i) + \delta +2$, 
so that 
$$A_{l} (-m(j) -N(i) + \delta +1 ) \subset \dis\bigcup^{N(i) - 
\delta -1}
_{-N(i) + \delta +2} \ov{P} (k).$$  
If $m(j) = 0$, then $D_{j}, 
D_{l}|_{E_{j}}$,
which implies that $D_{l} \succ D_{j}$ and therefore that $l > j$.  With this
we have shown that
\begin{eqnarray*}
A_{l} (-m(j) -N(i) + \delta + 1) & \subset & \bigcup _{l > j} A_{l}
(-N(i) + \delta + 1 ) \cup \bigcup^{N(i) - \delta -1} _{-N(i) + \delta +2}
\ov{P} (k) \\
& = & \bigcup_{l> j} f^{-N(i) + \delta +1} (D_{l} ) \cup \bigcup ^{N(i)-
\delta -1} _{-N(i) + \delta + 2} f^{k} (\ov{P})
\end{eqnarray*}

\noi where this last equaltiy is a conseqence of Proposition~\ref{56} (ii).
But, from the proof of Proposition~\ref{61}, we have seen that $I^{c} (\al_{i}
(N(i))) \cup \al_{i} (N(i))$ does not intersect the set after the equal
sign just above.  This finishes the proof.  $\Box$
\bigskip

\begin{cor} \label{63}
With the same notation as above, for every $i \in \un{L}$ the following 
holds:

\begin{description}
\item[(i) ]  if $n(i) < \infty$, then for every $j \in \un{L}$ such that 
$$f^
{N(i)} (D_{i}), f^{-N(i) + \delta +1} (D_{j})| _{f^{N(i)} (C{i}) }$$
we have
$$A^{c}_{i} (N(i)) \subset A^{c}_{j}(-N(i) + \delta + 1 ) | _ {f^{N(i) } 
(C_{i}) };$$

\item [(ii) ] if $m (i) >- \infty$, then for every $j \in \un{L}$ such that 
$$f^
{M(i)} (D_{i} ) \ f^{-M(i) + \delta '} ( D_{j}) | _ 
{f^{M(i)} (E_{i}) }$$
we have
$$ A^{e}_{i} (M (i) ) \subset A^{e}_{j} ( -M (i) + \delta ') | _{f^{M(i)} 
(E_{i})}.$$
 
\end{description}
\end{cor}

\noi {\sc Proof:}
We have shown that if $j \in \un{L}$ is largest such that 
$$f^{n(i)}(D_{i} ) ,
f(D_{j}) | _ {f^{N(i)} (C_{i} ) }$$ 
(which is equivalent to the condition in (i) 
above) then the desired inclusion holds.  Let $l \in \un{L}$ be such that
$f^{n(i)} (D_{i} ), f(D_{l} ) | _{f^{n(i) } (C_{i}) }$, $l \neq j$ (and thus
$l < j$.)  By Proposition~\ref{13}, $f(D_{j}), f(D_{l}) | _ {f^{n(i) } 
(C_{i} ) }
$ and since $l < j$, we must have $f(D_{l}) \prec f(D_{j})$.  Since $\{ 
\al_{i}
(1) \}_{i=1}^{L}$ is a ($\var,c$)-collection, it follows that $[ I^{c} 
( \al
_{j}(1) ) \cup \al_{j} (1) ] \subset I^{c} ( \al _{l} (1))$ and by 
Proposition
~\ref{60} this is equivalent to $A^{c}_{j}(1) \subset A^{c}_{l}(1) | _ 
{f(C_{j}
)} $ (see figure 23.)  Taking the $\Psi^{-N(i) + \delta}$-image of this 
latter inclusion, we get
$$ A ^{c} _{j} (-N(i) + \delta + 1) \subset A^{c}_{l} (-N(i) + \delta +1 ) 
| _{ \Psi ^{-N(i) + \delta } (f (C_{j}) ) }.$$

Notice that $\Psi^{-N(i) + \delta} (f (C_{j}) ) = \Psi^{-N(i) + \delta +1} (
\rho^{-1}_{0} ( C_{j} ))$ and that in the proof of Proposition~\ref{62} we
showed that $f^{N(i)}(C_{i}) \subset \Psi ^{-N(i) + \delta + 1 } ( \rho^{-1}
_{0}(C_{j}))$.  Therefore
$$A^{c}_{i} (N(i)) \subset A^{c}_{j} (-N(i) + \delta +1) | _ 
{f^{N(i)}(C_{i})} $$
\noi and
$$A^{c}_{j} (-N(i) + \delta + 1) \subset A^{c} _{l}( -N(i) + \delta +1 ) |
_{\Psi ^{-N(i) + \delta + 1} ( \rho^{-1}_{0}( C_{i}) ) } $$
\noi imply that
$$ A^{c}_{i} (N(i)) \subset A^{c}_{l} ( -N(i) + \delta + 1) | _ {f^{N(i)}
(C_{i}) } $$
\noi as we wanted.  $\Box$
\bigskip

\begin{figure}
\begin{center}~
\psfig{file=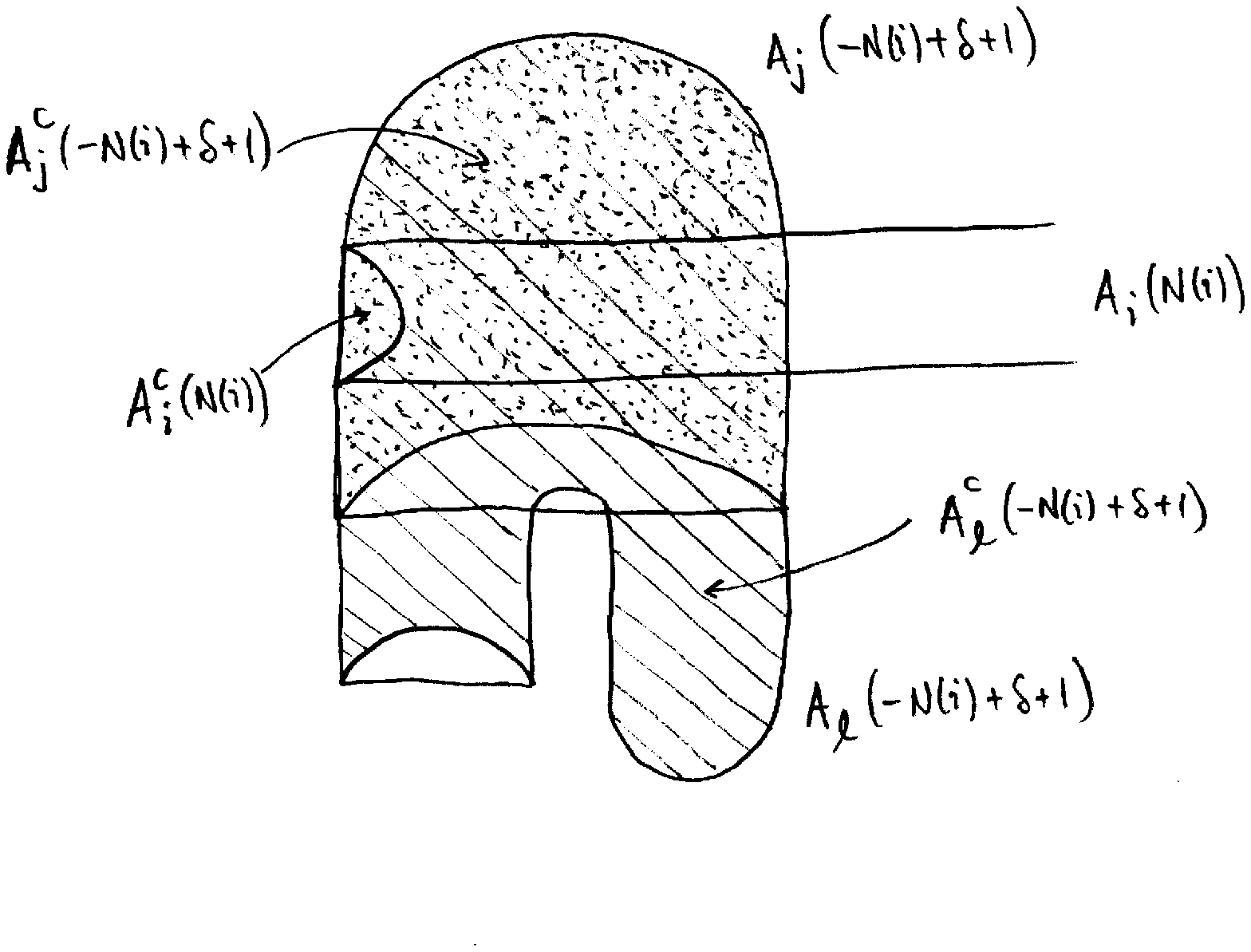,height=2.25in}
\end{center}
\caption{An example of $A_{i} (N(i)), \ A_{j} (-N(i) + \delta + 
1 ) $ and $A_{l} (-N (i) + \delta + 1)$.}
\label{f23}
\end{figure}

\begin{prop} \label{64}
$\Psi$ has the following properties:
\begin{description}

\item[(i)]  if $n(i) = \infty$, then $\Psi^{n} (A_{i}^{c})$ has interior
disjoint from $\ov{P}$ for every $n > 0$;

\item[(ii)]  if $n(i) < \infty$, then for every $j \in \un{L}$ such that $f^
{n(i)} (D_{i}), f(D_{j}) | _{f^{n(i)} (C_{i}) }$
$ \Psi ^{n(i)} (A^{c}_{i} ) \subset A^{c}_{j}(1) | _ {f^{n(i)} (C_{i})}$ and
$f^{n(i)} (C_{i}) \subset f(C_{j})$

\item[(iii)]  if $m(i) = - \infty$, then $\Psi^{m} (A^{e}_{i} ) $ has interior
disjoint from $\ov{P}$ for every $m < 0$;

\item[(iv)]  if $m(i) > - \infty$, then for every $j \in \un{L}$ such that
$f^{m(i)} (D_{i}), D_{j} | _{f^{m(i)} (E_{i}) }$ $\Psi^{m(i)} (A^{e}_{i} )
\subset A^{e}_{j} | _{f^{m(i)} (E_{i}) }$ and $f^{m(i)} (E_{i}) \subset E_j$.

\end{description}
\end{prop}

\noi {\sc Proof:}
If $n(i) = \infty$, then $N(i) = \infty$ and by Proposition~\ref{60} we see
that $\Psi^{n} (A^{c}_{i} ) = A^{c}_{i} (n) = D^{c} ( \al_{i}(n) )$ for every
$n \geq 1$.  Since $f^{n} (D_{i}), f^{k} (D_{j}) \not| _{f^{n}(C_{i}) } $ for
$-n+1 \leq k \leq n-1$, 
$$[I^{c}( \al_{i}(n) ) \cup \al_{i} (n) ] \cap \dis
\bigcup^{n-1}_{-n+1} f^{k} ( \ov{P}) = \es.$$  
This being true for every 
$n \geq
1$, we se that $[ I^{c} ( \al_{i}( n)) \cup \al_{i}(n) ] \cap \ov{P} = \es$ for
every $n \geq 1$, which proves (i).  (ii) is immediate from 
Corollary~\ref{63}.
(iii) and (iv) are analogous and we omit the proofs.  $\Box$
\bigskip

We can now state and prove the main theorem.

\begin{theorem} [Main Theorem] \label{65}
Let $f: \pi \ra \pi$ be a homeomorphism of the plane, $\{D_{i} 
\}^{L}_{i=1}$ a
pruning collection and $P = \dis \bigcup^{L}_{i=1} I_{i}$, where $I_i$ is the
interior of the disk $D_i$.  Then there exists an isotopy $H: \pi \times [0,1]
\ra \pi$ of the identity such that supp $H \subset \dis \bigcup_{k \in{\Bbb{Z}}
}f^{k} (P)$, and if we set $f_{P}(\cdot) = f \circ H ( \cdot, 1)$, every 
point of $P$ is wandering under $f_P$.
\end{theorem}

\noi{\sc Proof:}
Construct a directed graph $G_c$ as follows:  its vertices are the integers
$\{ i \in \un{L}; \ n(i) > 1 \}$ and there is a directed vertex from $i$ to
$j$ if $n(i) < \infty$ and $f^{n(i)} (D_{i}), f(D_{j}) | _{f^{n(i)}(C_{i})}$.
Since we have taken only $i \in \un{L}$ for which $n(i) > 1$, it is easy 
to see
that from each vertex there is at most one outgoing edge (or none, if $n(i) =
\infty$).  A {\em loop} in the directed graph consists of an ordered set of
distinct vertices $\{i_{1} < i_{2} < \ldots < i_{l} \}$ such that there is a
directed edge from $i_{r}$ to $i_{r+1}$, for $1 \leq r \leq l$ where we let the
indices ``wrap around'', i.e., $l + 1$ ``=''$1$.  Since there is at most one
edge emanating from each vertex, and the vertices in a loop are ordered and
distinct, it follows that two loops are either equal or disjoint.  Let ${\mc{L}
} = \{ i_{1}, \ldots , i_{l} \}$ be a loop in $G_c$, which for now we will
represent by just its subscripts $\{1, \ldots, l \}$ so that the notation is
not too awkward.  By definition, we have
$$ f^{n(r)} (D_{r}), \ f(D_{r+1}) | _{f^{n(r)}(C_{r})} \ \mbox{\rm for } 1 
\leq r \leq l$$
\noi and by Proposition~\ref{64}
$$ \Psi^{n(r)} (A^{c}_{i}) \subset A_{r+1} (1) |_ {f^{n(r)}(C_{r}) } 
\ \mbox{\rm for } 1 \leq r \leq l$$

\noi from which it follows that
$$ \Psi^{ \sum^{l}_{r=1}  n(r) - (l-2) } (A^{c}_{1}(1 )) \subset 
A^{c}_{1}(1) |_{f^{\sum^{l}_{r=1}  n(r) - (l-1) } (C_{i} ) } $$

For a loop ${\mc{L}} = \{i, \ldots , i_{l} \}$ let $n ({\mc{L}} ) = \dis\sum
^{l}_{r=1} n(i_{r}) - (l -2)$.  By Lemma~\ref{40} and Corollary~\ref{41},
there exits an isotopy $h_{{\mc{L}} }$ of the identity with supp $h_{ {\mc{L}
} }$  $\subset I (A^{c} _{i_{1}} (1) )$ such that, if $\zeta_{ {\mc{L} } 
}(\cdot)
= h _{ {\mc{L}}} (\cdot, 1 )$, then $(\Psi \circ \zeta_{ {\mc{L}} })^{k n( {\mc
{L}} ) }(x) \rightarrow p$ for every $ x \in A^{c}_{i} (1)$, where
$p$ is the fixed point of $\Psi^{n( {\mc{L}}) }$ in $f^{n( {\mc{L} } ) +1}
(C_{i})$.  We then construct isotopies $h_{{\mc{L}} }$ for each loop ${\mc{L}}$
in $G_c$.  Since the vertices of $G_c$ where integers $i \in \un{L}$ for which
$n(i) > 1$, the supports of isotopies associated to different loops are
disjoint.  Let $h_c$ be the union of all these isotopies.  By
construction supp $h_{c} \subset \dis \bigcup^{l}_{i=1} I (A^{c}_{i} (1)) =
\dis \bigcup^{L}_{i=1} I^{c} ( \al_{i} (1) )$.

In an analogous manner, we construct a directed graph $G_e$ whose 
vertices are
$\{ i \in \un{L}; m(i) < 0 \}$ and for each loop $\mc{L}$ in $G_e$, we
construct an isotopy $h_{\mc{L}}$ of the identity, with support in
$A^{e}_{i_{1}}(1)$, playing the analogous role for $\Psi^{-1}$ as the above
ones played for $\Psi$.  Let $k_e$ denote the union of these isotopies, for it
is again easy to check that they have disjoint supports, and define $h_{e} =
\Psi^{-1} \circ k_{e}^{-1} \circ \Psi$, i.e., for each fixed $t$, $h_{e}(x,t)
= \Psi^{-1} ( k_{e}^{-1} ( \Psi (x),t) )$.  $h_e$ is also an isotopy 
of the
identity and since supp $k_{e} \subset \dis \bigcup^{L}_{i=1} I ( A^{e}_{i}
(1))$, 
$$\mbox{\rm supp } h_{e} \subset \dis \bigcup^{L}_{i=1} I (A^{e}_{i} (0)) = 
\dis
\bigcup ^{L}_{i=1} I^{e} ( \al_{i} (0)).$$  
>From this it follows that supp
$h_{e} \cap$ supp $h_{c} = \es$ and we let $h = h_{c} \cup h_e$ and $\zeta
(\cdot ) = h ( \cdot, 1 )$.  Finally set
$$ H(x,t) = \left\{ \begin{array}{ll}
 h(x, 2t)   & t \in \left[ 0, \dis \frac{1}{2} \right] \\
  \\
R(\zeta(x), 2t-1)   & t \in \left[ \dis \frac{1}{2}, 1 \right]
\end{array}
\right. .
$$

It is now not hard to check that $H$ has the desired properties.  $\Box$
\bigskip

%% file: examples.tex
\section{Examples}

In this section we present examples of pruning collections for
Smale's horseshoe map $f:{\Bbb R}^2\to{\Bbb R}^2$. We begin by choosing
a rigid model for $f$ and describing some well known results,
offered without proof. We also present some elementary concepts of
kneading theory which we will need (the reader is referred to the
books of Wiggins \cite{Wi}, Devaney \cite{De}, and de Melo and Van Strien
\cite{MS} for further details on the horseshoe and on 1-dimensional
dynamics.) We then get to the examples. In describing the dynamics
of the ``pruned" maps $f_P$ in each example, we will make several
assertions and only sketch the proofs. The reason for proceeding
thus is twofold. First, this is a section to give examples of
pruning collections and this aspect is presented fully.    Second, the
details we omit are part of a more general theory deserving of separate
treatment, which we intend to do in forthcoming papers.

\begin{figure}
\centerline{\psfig{file=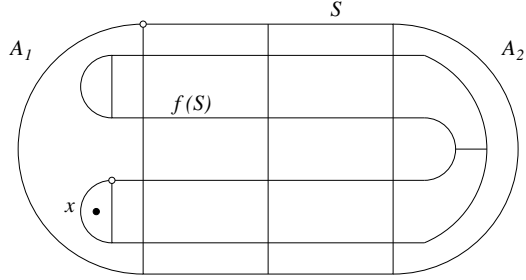,height=1.5in}}
\caption{A rigid model for the horseshoe.}\label{hsmodel}
\end{figure}

We now fix a rigid model of Smale's horseshoe map
$f:{\Bbb R}^2\to{\Bbb R}^2$.  Foliate the square
\[
S=\{ (x,y): |x|\leq \textstyle\frac 12, \
|y|\leq \textstyle\frac 12\}
\]
with horizontal {\it unstable} leaves and vertical
{\it stable} leaves, and begin by choosing the action  of $f$ on
$S$ as depicted in figure~\ref{hsmodel}. We require that $f$ should
stretch the unstable leaves uniformly, contract the stable leaves
uniformly, and map segments of unstable (respectively, stable)
leaf in $S\cap f^{-1}(S)$ onto segments of unstable (respectively, stable)
leaf in $S$. Morever, we choose $f$ to map the corner of $S$ marked
with a circle on figure~\ref{hsmodel}  onto the corner of $f(S)$ marked with a
circle. Extend $f$ to the half-disks $A_1$ and $A_2$ as depicted
in the diagram: let $f$ be a strict contraction of $A_1\cup A_2$,
so that there is a fixed point $x$ of $f$ lying in $A_1$ with the
property that $f^i(y)\to x$ \ as \ $i\to\infty$ for all
$y\in A_1\cup A_2$. Finally, extend $f$ over the rest of ${\Bbb R}^2$
without introducing any new nonwandering points.

The nonwandering set $\Omega(f)$ of $f$ consists of the fixed point
$\{x\}$ and an invariant Cantor set $\Lambda\subset S$. Morever,
there exists a homeomorphism $h:\Sigma\to\Lambda$, where
$\Sigma=\{0,1\}^{\Bbb Z}$ is the two-sided shift on two symbols,
which conjugates the shift map $\sigma:\Sigma\to\Sigma$ and
$f|_\Lambda:\Lambda\to\Lambda$.
$\Lambda$ is a hyperbolic invariant set and each point
$p\in\Lambda$ has one-dimensional stable and unstable manifolds
which intersect transversally. Notice that if two points $p_0$ and
$p_1$ lie on the stable (unstable) manifold of some point $q\in\Lambda$,
$p_0$ and $p_1$ are the endpoints of exactly one arc contained in
the stable (unstable) manifold of $p$.

We now describe the {\it unimodal order} on the one-sided shift space
$\Sigma_+=\{ 0,1\}^{\Bbb N}$ and, using it, define
{\it kneading sequences}.

\bigskip
\begin{defn} Let $s=s_0s_1\dots$  and
$t=t_0t_1\dots$   lie in  $\Sigma_+$   and suppose
$s_i=t_i$   for  $i < k$ {\it and} $s_k\neq t_k$.   We set
$s\lhd t$   if  $\sum^k_{i=1} s_i$   is even. We set
$s\unlhd t$   if either  $s=t$   or  $s\lhd t$.
\end{defn}

\begin{defn}
Let  $\sigma: \Sigma_+\to\Sigma_+$
be the shift map and  $\kappa\in\Sigma_+$.  We say
$\kappa$   is a  {\em kneading sequence}   if, for every
$n\in\Bbb N$, \ $\sigma^n(\kappa)\unlhd\kappa$.
\end{defn}

\bigskip
The unimodal order just defined is used in the study of
1-dimensional {\it unimodal} maps, i.e., piecewise monotone
endomorphisms of the interval   with exactly one critical
(turning) point. In this context, kneading sequences are defined as
the itinerary of the critical value. It is possible to check that
kneading sequences associated to unimodal maps satisfy the
definition above.

The unimodal order describes the horizontal and vertical ordering
of points in $\Lambda$ as follows: if $(x_1,y_1),(x_2,y_2)\in\Lambda$,
with $h(x_1,y_1)=s_{-2}s_{-1}\cdot s_0s_1\dots$ and
$h(x_2,y_2)=t_{-2}t_{-1}\cdot t_0t_1\dots$, then
\begin{eqnarray*}
&& x_1 < x_2 \Longleftrightarrow s_0 s_1 s_2\dots \lhd
   t_0 t_1 t_2 \dots \ \ {\mbox{and }} \\
&& y_1 < y_2 \Longleftrightarrow s_{-1} s_{-2}\dots \lhd
   t_{-1} t_{-2} \dots
\end{eqnarray*}

We shall often use the elements of $\Sigma$ to describe
points of $\Lambda$ without explicitly invoking the map $h$.
Thus, for example, we may talk about ``the fixed point $\bar 1$,"
``the periodic orbit $\overline{10011}$," or ``the point
$\overline 0 .0\overline{101}$."
Here, a bar over a group of symbols stands for infinite
repetition of the group. If the group is to the right (left) of the decimal
point, it should be repeated infinitely to the right (left) and if there is
no decimal point, the group should be repeated infinitely to both sides
(so $\overline 0 .0\overline{101}=\dots 000.0101101101\dots$ and
$\overline{10}=\dots 1010.1010\dots$ .) If the symbolic sequence is an element
of $\Sigma_+$, a bar over a group of symbols means infinite repetition of the
group to the right. Let $p\in\Lambda$ and $h(p)=
s_{-2}s_{-1}\cdot s_0s_1\dots$, we will sometimes refer
to $\dots s_{-2}s_{-1}\cdot s_0s_1\dots$
as the {\it symbolic representation} of $p$ and to
$\dots s_{-2}s_{-1}.$ and
$.s_0s_1\dots$ as the {\it symbolic vertical and horizontal coordinates 
of} $p$, respectively.

\begin{figure}
\centerline{\psfig{file=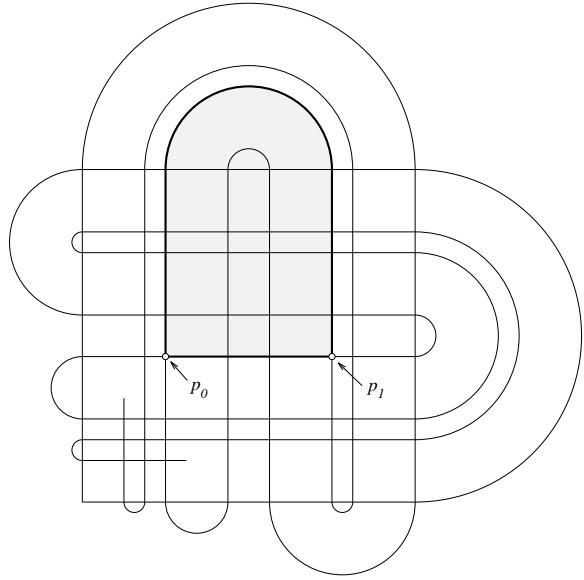,height=3in}}
\caption{\label{cedisk-ex}$(c,e)$-disk determined by $p_0=
\overline 010.011\overline 0$ and $p_1=\overline 010.111\overline 0$}
\end{figure}

If two points $p_0,p_1\in\Lambda$ have symbolic representation
$\dots t_{-2}t_{-1}\cdot t_0t_1\dots$ and
$\dots s_{-2}s_{-1}\cdot s_0s_1\dots$, respectively, $p_0$
and $p_1$ lie on the same stable (unstable) manifold if there
exists $N\in\Bbb Z$ such that $s_i=t_i$ for every $i\geq  N$
($i\leq N$). Consequently, $p_0$ and $p_1$ lie on the same stable {\it and}
unstable manifolds if their symbolic representations differ in at most
finitely many entries. If the symbolic representations of $p_0$ and
$p_1$ differ at exactly one entry, the stable and unstable arcs of
which they are the endpoints form a simple closed curve bounding a
closed disk which we denote by $D(p_0,p_1)$ (see figure~\ref{cedisk-ex}.) 
Because the boundary
of $D(p_0,p_1)$ is the union of a stable and an unstable arcs,
$D(p_0,p_1)$ is a $(c,e)$-disk for $f$ as defined in Section~\ref{cedisks} 
whose vertices are $p_0$ and $p_1$.

\bigskip
\noindent
{\sc Notation:}  If $p_0,p_1\in\Lambda$ lie on the same stable
(unstable) manifold, we denote the closed arc of stable (unstable)
manifold whose endpoints are $p_0$ and $p_1$ by
$[p_0,p_1]_s$ ($[p_0,p_1]_u$).
Let $s=s_0s_1\dots\in\Sigma_+$. We define the {\it vertical segment
with  horizontal coordinate} $s$ to be $[ \,\overline 0.s,
\overline 01.s]_s$ and denote it by ver$(s)$. A vertical
segment is thus an arc of stable manifold extending from the lowest to the
highest possible symbolic vertical coordinates (notice that 000$\dots$
and 100$\dots$ are, respectively, the smallest and largest elements of
$\Sigma_+$ in the unimodal order), having symbolic horizontal coordinate
$s$. Notice also that, if we use $\sigma$ to denote the shift map on 
$\Sigma_+$, $f({\mbox{ver}}(s))\subset{\mbox{ver}}(\sigma(s))$.

We are now ready to present examples of pruning collections for $f$.

\bigskip
\noindent
{\sc Remark:} So that the figures below are not hopelessly complicated and
unintelligible, we will represent the Cantor set $\Lambda$ as a solid
square. Formally, what we are depicting is the quotient of $\Lambda$ under
an equivalence relation which collapses the ``gaps" of the vertical and
horizontal Cantor sets, the product of which  is $\Lambda$.
The ambiguity thus created is  easily understood and the clarity gained
plentifully compensates it.

\begin{figure}
\centerline{\psfig{file=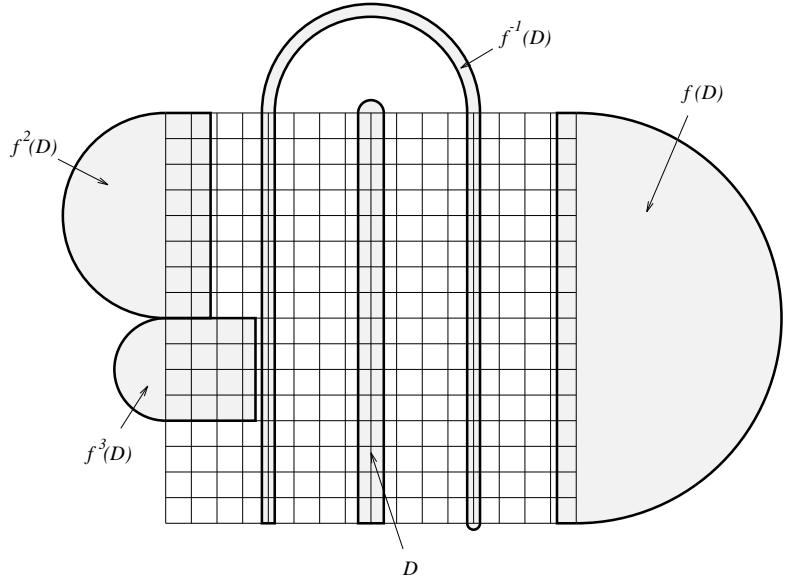,height=3in}}
\caption{A ``one-dimensional-like'' pruning front
 for the horseshoe.}\label{odl}
\end{figure}

\bigskip
{\bf Example 1}. Let $\kappa$ be a kneading sequence and
\[D=D(\overline 0.0\kappa,\overline 0.1\kappa)
\]
be the $(c,e)$-disk determined by $\overline 0.0\kappa$
and $\overline 0.1\kappa$ (that is, the disk bounded by
the union of $C=[ \,\overline 0.0\kappa,
\overline 0.1\kappa]_s$ and $E=[ \,\overline 0.0\kappa,
\overline 0.1\kappa]_u$). We claim that the collection
$\{D\}$, containing $D$ alone, is a pruning collection (see figure~\ref{odl}.)
In order to see this we have to show that, if
$f^k(I)\cap I\neq\emptyset$ (where $I$ is the interior of $D$),
then $f^k(D)\succ D$, if $k > 0$, and $f^k(D)\prec D$, if $k < 0$.
Notice that conditions (ii) and (iii) in Definition~\ref{longer}
are automatically satisfied since $C$ and $E$ are arcs of stable
and unstable manifolds which intersect transversally. All there is
left to check is that $f^n(C)\cap I=\emptyset$ for $n > 0$
and $f^m(E)\cap I=\emptyset$ for $m < 0$.

Notice that $E\subset[ \,\overline 0,\overline 0.1\overline 0 \,]_u$,
that
\[
f^{-1}([ \,\overline 0, \overline 0.1\overline 0 \,]_u) =
[ \,\overline 0, \overline 0.01\overline 0 \,]_u \subset
[ \,\overline 0, \overline 0.1\overline 0 \,]_u
\]
and that
\[
[ \,\overline 0,\overline 0.1\overline 0 \,]_u\cap I=\emptyset \ .
\]
Thus, if $m\leq -1$, \ $f^m(E)\cap I=\emptyset$.
On the other hand, observe that $f(C)=$ ver$(\kappa)$ and, therefore,
$f^n(C)\subset$ ver$(\sigma^{n-1}(\kappa))$ for $n\geq 1$.
If $f^n(C)\cap I\neq\emptyset$, then $f^n(C)\subset I$ and,
in fact, ver$(\sigma^{n-1}(\kappa))\subset I$. This implies that
\[
0\kappa\lhd \sigma^{n-1}(\kappa)\lhd 1\kappa
\]
and, applying $\sigma$ to this inequality, we get
$\kappa\lhd\sigma^n(\kappa)$, which contradicts the assumption
that $\kappa$ is a kneading sequence.

Let $f_\kappa$ denote the map obtained using Theorem~\ref{65} 
for $\overline P=D$. The family ${\cal F}=\{f_\kappa; \kappa$ is a kneading
sequence$\}$ mimics in dimension 2 a {\it full family} of
unimodal maps of the interval. In particular, ${\cal F}$ is an
uncountable family of 2-dimensional homeomorphisms passing from
trivial dynamics to a full horseshoe as $\kappa$ varies from
$\overline 0$ to $1\overline 0$.

\begin{figure}
\centerline{\psfig{file=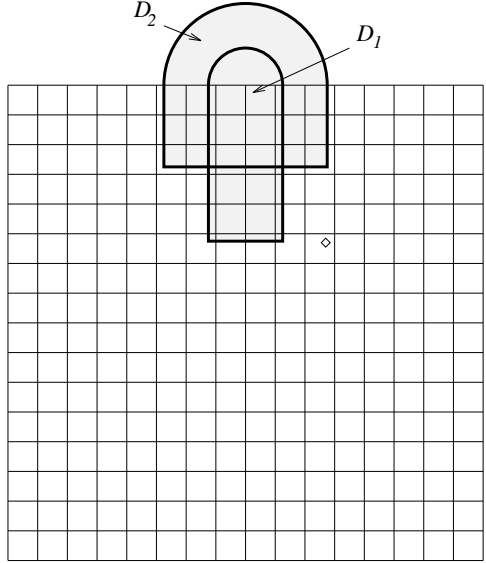,height=3in}}
\caption{The $(c,e)$-disks $D_1$ and $D_2$ of Example 2. $\diamond$ is the 
fixed point $\bar 1$.}\label{pf-ex2}
\end{figure}

\begin{figure}
\centerline{\psfig{file=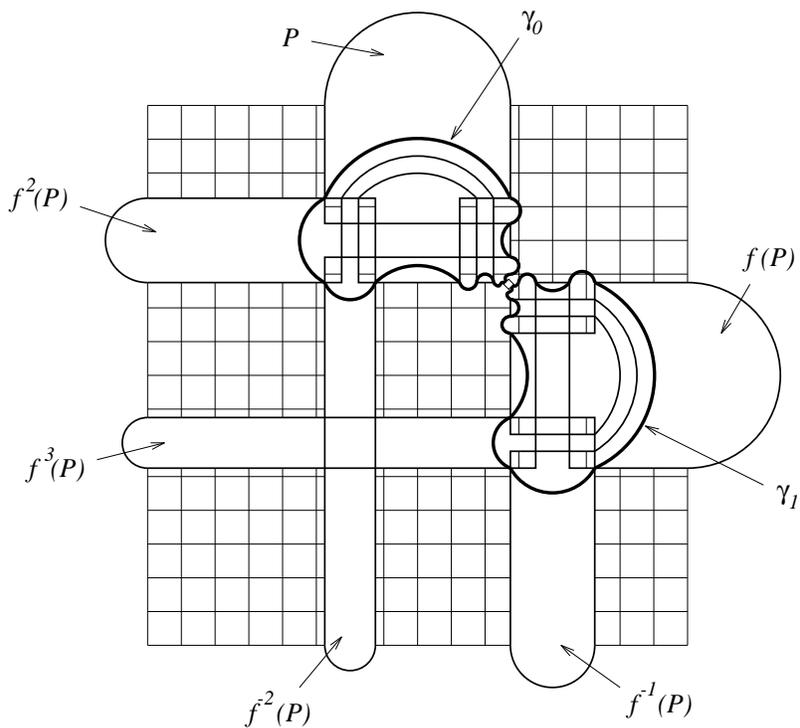,height=3.8in}}
\caption{The Jordan curves $\gamma_0$ and $\gamma_1$}\label{renorm}
\end{figure}

{\bf Example 2}. Consider the two $(c,e)$-disks
\[
D_1=D(\overline 1.010\overline 1, \ \overline 1.110\overline 1) \qquad
{\mbox{ and }} \qquad
D_2=D(\overline 101.0\overline 1, \ \overline 101.\overline 1) \ .
\]
Because of the periodicity in the coordinates of the vertices of
$D_1$ and $D_2$, it is easy to check that $\{D_1,D_2\}$ is a
pruning collection (see figure~\ref{pf-ex2}.) Let
 $\alpha_i(k)\subset f^k(D_i)$, for
$i=1,2$, $k\in\Bbb Z$, be closed cross-cuts as given by 
Proposition~\ref{28}. Then
\[
\gamma_0=\bigcup\{ \alpha_i(2k); \ i=1,2, \ k\in\Bbb Z\}
\]
and
\[
\gamma_1=\bigcup\{ \alpha_i(2k-1); \ i=1,2, \ k\in\Bbb Z\}
\]
are Jordan curves (see figure~\ref{renorm}) such that
$\gamma_0\cap\gamma_1=\{\,\overline 1\,\}$. If $f_P$ is the map given by 
Theorem~\ref{65} for the pruning front $\overline P=D_1\cup D_2$
and $U_0$ and $U_1$ are the closed disks bounded by $\gamma_0$ and
$\gamma_1$, $f_P$ interchanges $U_0$ and $U_1$, that is,
$f_P(U_0)=U_1$ and $f_P(U_1)=U_0$. Moreover, if $\Lambda_0=\Omega(f_P)\cap
U_0$ is the intersection of the nonwandering set of $f_P$ with
$U_0$, \  $f^2_P|_{\Lambda_0}:\Lambda_0\to\Lambda_0$ is topologically conjugated
to the full horseshoe $f|_\Lambda: \Lambda\to\Lambda$ restricted to the set $\Lambda$.
This is an example of a ``renormalizable" or ``reducible" map.

\begin{figure}
\centerline{\psfig{file=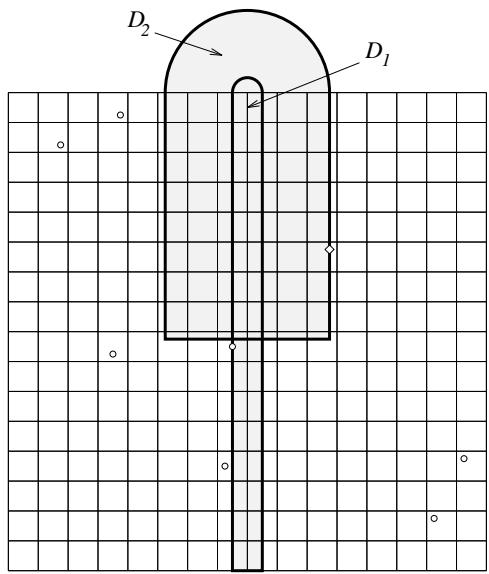,height=3in}}
\caption{The $(c,e)$-disks $D_1$ and $D_2$ of Example 3. The points marked by 
$\circ$ are the periodic orbit $s_{7}^{6}(0)$ and $\diamond$ is the fixed 
point $\bar 1$.}\label{pf-ex3}
\end{figure}

\bigskip
{\bf Example 3}. Consider the $(c,e)$-disks
\[
D_1=D(\overline 0.0\overline{1000100}, \
\overline 0.1\overline{1000100}) \qquad
{\mbox{ and }} \qquad
D_2=D(\overline 01100010.0\overline 1, \
\overline 01100010.\overline 1) \ .
\]
Notice that $D_1$ is of the kind $D(\overline 0.0\kappa,\overline
0.1\kappa)$,   since $\kappa = \overline{1000100}\in \Sigma_+$ 
is a kneading sequence, as may easily be verified.
 As in the previous example, the
periodicity in the coordinates of the vertices of $D_1$ and $D_2$
make it an easy computation to check that $\{D_1,D_2\}$ is a
pruning collection (see figure~\ref{pf-ex3}.)

\begin{figure}
\centerline{\psfig{file=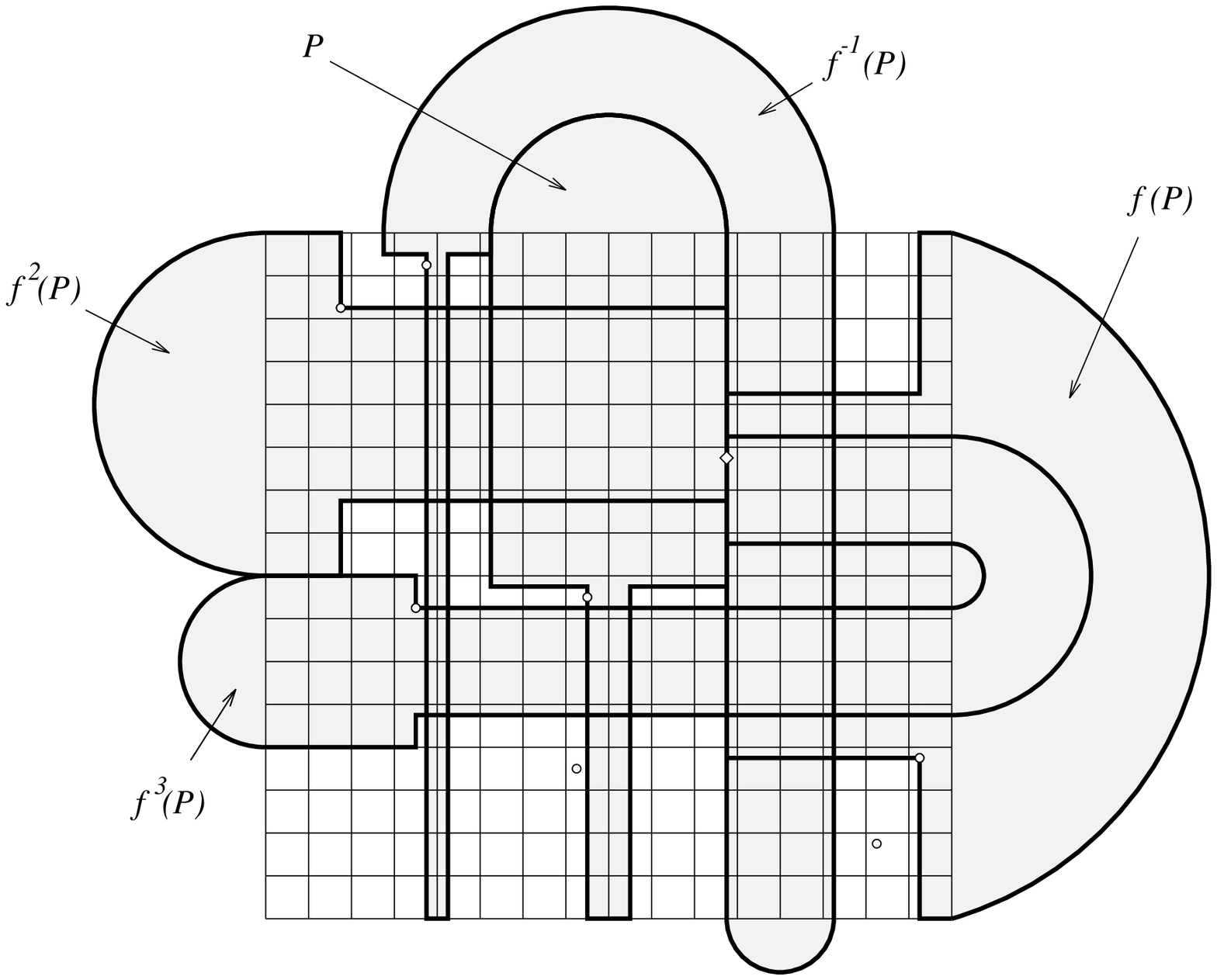,height=3in}}
\caption{A few images of $P$ under $f$.}\label{im-pf}
\end{figure}

Let $s^6_7(0)$ denote the periodic orbit containing the point
$\overline{1000100}$ (see \cite{HWh} for an explanation of this name).
Since none of its seven points  lies in $P=I_1\cup I_2$, none of them
lies in $\displaystyle{\bigcup_{n\in\Bbb Z}} f^n(P)$.
 In figure~\ref{im-pf}, $f^k(P)$, for $-1\leq k\leq 3$, are shown.
Notice that the nonwandering points of $f_P$ lie outside the shaded region.
(In fact, it is not hard to see that the points on the open $e$-side
$\stackrel{\circ}{E}_1$ of $D_1$ are also wandering under $f_P$.)

\begin{figure}
\centerline{\psfig{file=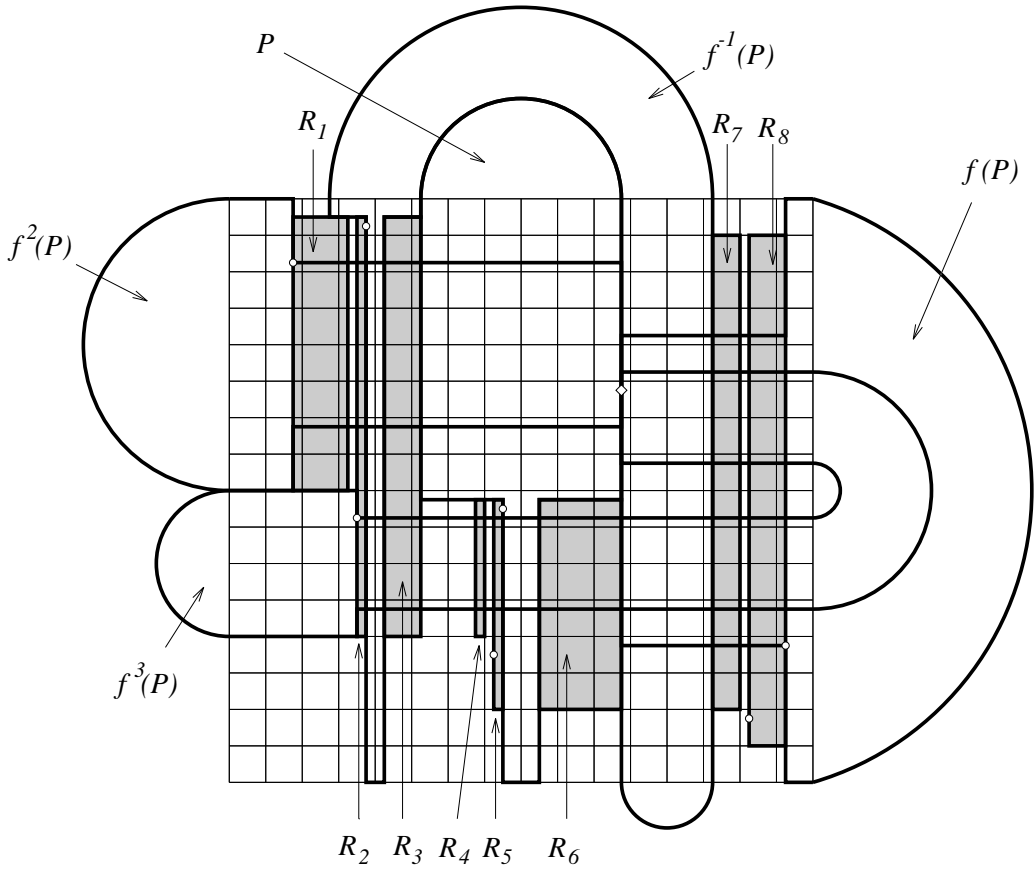,height=4.5in}}
\caption{The Markov Partition for $f_P$.}\label{mp}
\end{figure}

We claim that the map $f_P$ obtained using Theorem~\ref{65}
realizes the minimum topological entropy among all maps in the
isotopy class of $f$ relative to $s^6_7(0)$. In order to see this, we
construct a Markov Partition for $f_P$ like in figure~\ref{mp}.  
The horizontal and vertical sides of the rectangles $R_i$ are
contained in $\displaystyle{\bigcup_{n\in\Bbb Z}} f^n_P(E_1\cup E_2)$
and $\displaystyle{\bigcup_{n\in\Bbb Z}} f^n_P(C_1\cup C_2)$, respectively.
(In fact, it is enough to take the unions ranging from $n=-7$ to $n=7$,
say.) It is easy to see that, if we define $\Lambda_P=\Omega(f_P)\backslash
\{\overline 0,\overline 1\}$, where $\Omega(f_P)$ is the nonwandering set
of $f_P$ and $\overline 0$ and $\overline 1$ are the fixed of $f$ inside
$S$, which are also fixed points under $f_P$, then $\Lambda_P\subset
\displaystyle{\bigcup^8_{i=1}} R_i$. The vertices of each $R_i$ lie outside of
$\displaystyle{\bigcup_{n\in\Bbb Z}} f^n(P)$ and it therefore makes sense
to refer to them using their symbolic representation in $\Sigma$. In 
table~\ref{corners-mp} we give the symbolic horizontal and vertical
coordinates of the vertices of the rectangle $R_i$. The columns under $x_L$ 
and $x_R$ contain the left and right horizontal coordinates, respectively, 
whereas those under $y_L$ and $y_U$ contain the lower and upper vertical
coordinates, respectively.

\begin{table}
\begin{center}
\renewcommand{\arraystretch}{1.2}
\[  \begin{array}{|c|l|l|r|r|}  \hline
     &x_L     &x_R     &y_L     &y_U    \\  \hline
R_1  &.\bar{0001001} &.00101\bar{1000100} &\bar 011. &\bar 0110001.\\ \hline
R_2  &.\bar{0010010}&.\bar{0010001}&\bar 0110.&\bar 0110001. \\ \hline
R_3  &.01\bar{1000100}&.0\bar 1&\bar 0110.&\bar 0110001. \\ \hline
R_4  &.010\bar 1&.0101\bar{1000100}&\bar 0110.&\bar 01100010. \\ \hline
R_5  &.\bar{0100100}&.\bar{0100010}&\bar 01100.&\bar 01100010. \\ \hline
R_6  &.1\bar{1000100}&.\bar 1 &\bar 01100.&\bar 01100010. \\ \hline
R_7  &.10\bar 1&.101\bar{1000100}&\bar 01100.&\bar 011001.\\ \hline
R_8  &.\bar{1001000}&.\bar{1000100}&\bar 011000.&\bar 011001.\\ \hline
\end{array}  \]
\caption{The coordinates of the vertices of the rectangles $R_i$.}
\label{corners-mp}
\end{center}
\end{table}

\begin{figure}
\centering
\mbox{\subfigure[The rectangles $R_i$ and \ldots]
{\psfig{file=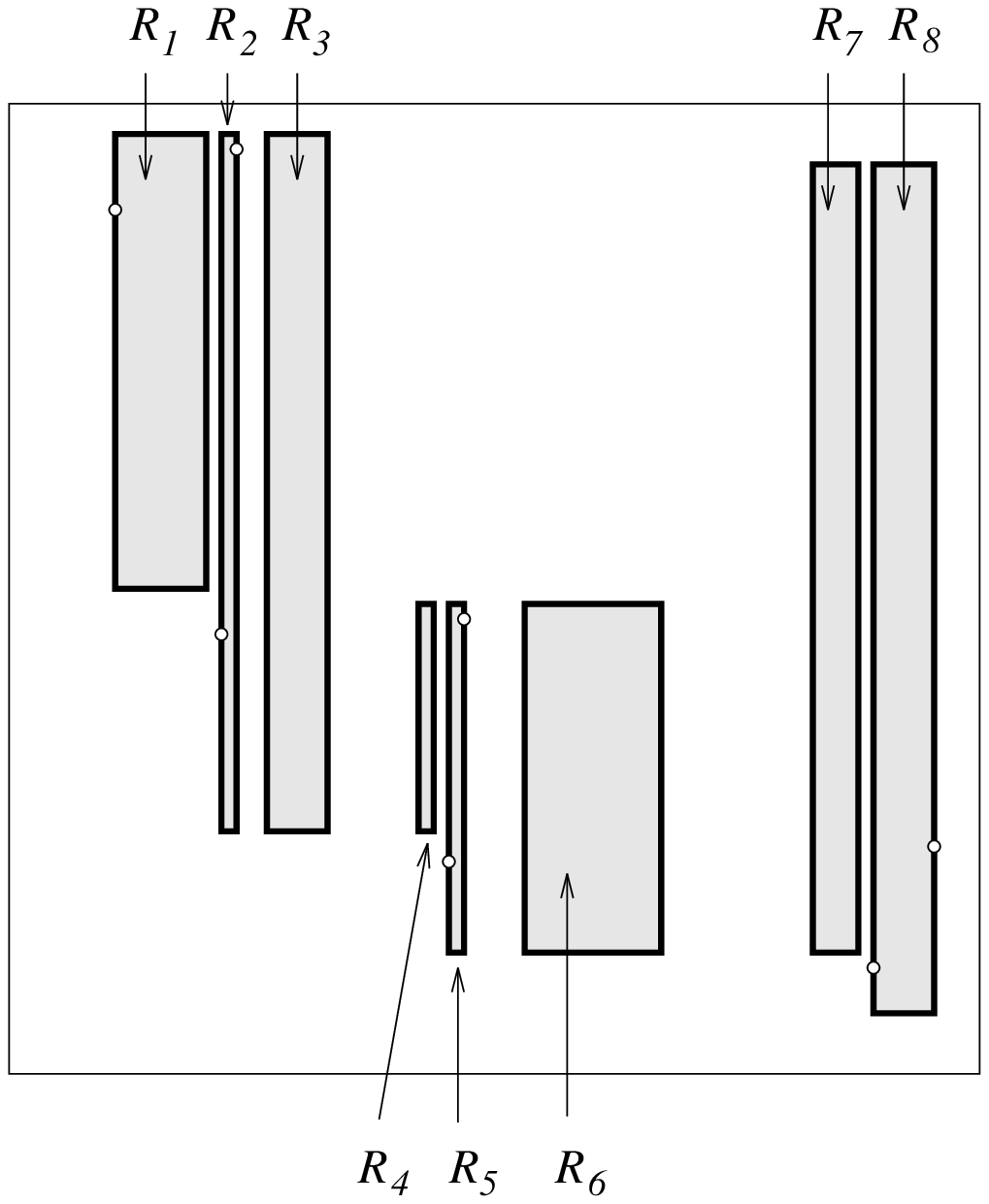,height=3in}}\qquad
      \subfigure[\ldots their images under $f_P$.]
{\psfig{figure=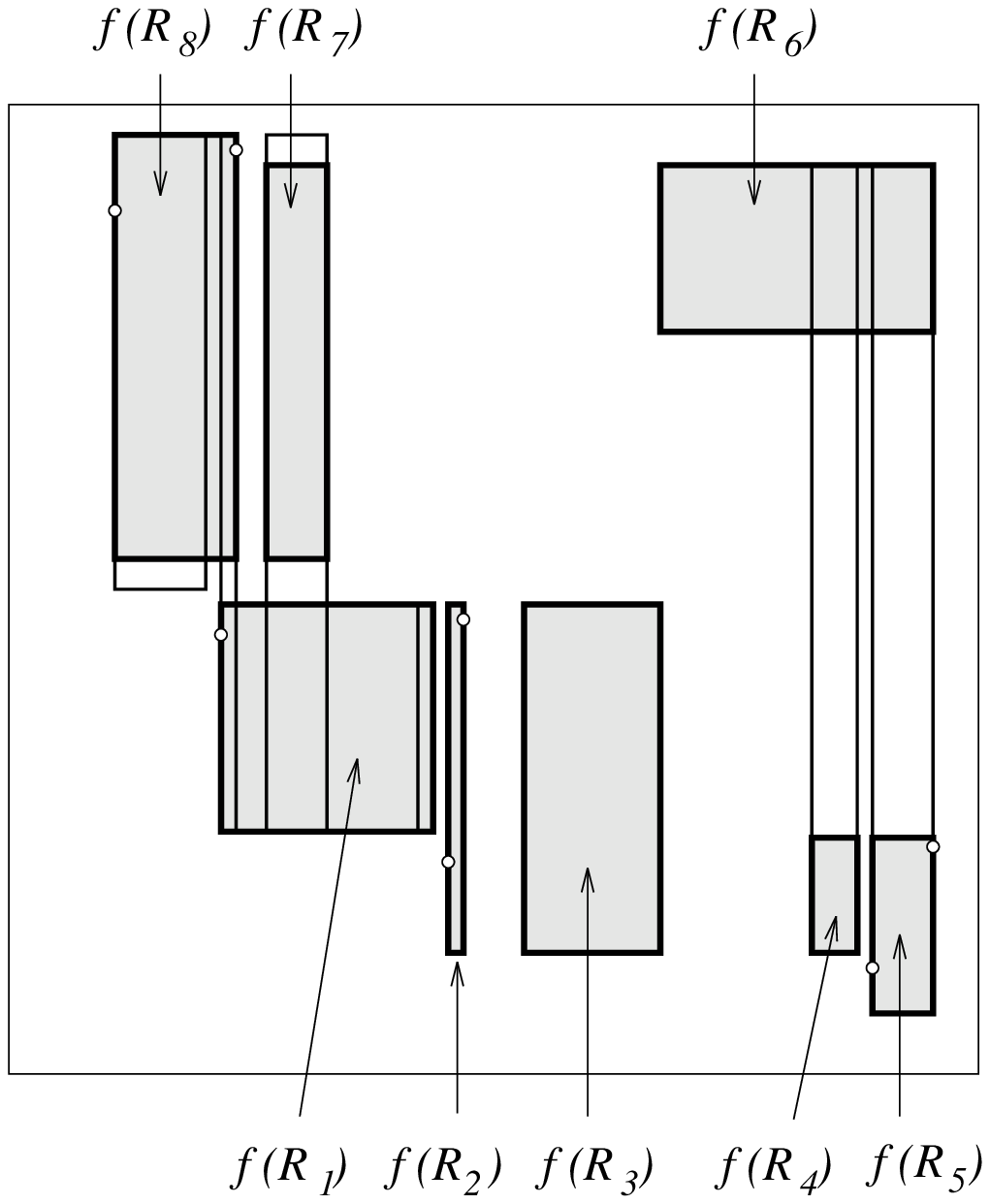,height=3in}}}
\caption{The Markov Partition and its image.}
\label{mapmp}
\end{figure}

In figure~\ref{mapmp} we show how the rectangles $R_i$ are mapped under $f_P$.
 The transition matrix $M=(m_{ij})$ 
associated with this partition is the $8\times 8$ matrix defined by
\[
m_{ij} =\cases 1 \ , & {\mbox{if }} \ I(f_P(R_j))\cap R_i\neq \emptyset \\
  0 \ , & {\mbox{otherwise}}
\endcases
\]
where $I(R_i)$ stands for the interior of the rectangle $R_i$.
Using the notation $R_j\to R_i$ for
$I(f_P(R_j))\cap R_i\neq\emptyset$, we have
$R_1\to R_2R_3R_4$, \ $R_2\to R_5$, \ $R_3\to R_6$, \ $R_4\to R_7$, \
$R_5\to R_8$, \ $R_6\to R_8R_7$,  \ $R_7\to R_3$,    \ $R_8\to R_2R_1$, \
so that
\[
M=\left[
\begin{array}{llllllll} 0 & 0 & 0 & 0 & 0 & 0 & 0 & 1 \\
  1 & 0 & 0 & 0 & 0 & 0 & 0 & 1 \\
  1 & 0 & 0 & 0 & 0 & 0 & 1 & 0 \\
  1 & 0 & 0 & 0 & 0 & 0 & 0 & 0 \\
  0 & 1 & 0 & 0 & 0 & 0 & 0 & 0 \\
  0 & 0 & 1 & 0 & 0 & 0 & 0 & 0 \\
  0 & 0 & 0 & 1 & 0 & 1 & 0 & 0 \\
  0 & 0 & 0 & 0 & 1 & 1 & 0 & 0
\end{array}   \right]
\]
Let
\[
\Sigma_M =\{ s=\dots s_{-2}s_{-1}\cdot s_0 s_1\dots\in
 \{ 0,1,\dots 8\}^{\Bbb Z}; m_{s_is_{i+1}}=1 \ \forall i\in\Bbb Z\}
\]
be the subshift of finite type associated to $M$. It is possible to show
that, if $s=\dots s_{-2}s_{-1}\cdot s_0 s_1\dots\in\Sigma_M$, then
$\displaystyle{\bigcap_{n\in\Bbb Z}}f^{-n}_P(R_{s_n})$ consists of a single
point in $\Lambda_P$ and that the map $k:\Sigma_M\to\Lambda_P$, given by
$k(s)=\displaystyle{\bigcap_{n\in\Bbb Z}}f^{-n}_P(R_{s_n})$, is a
topological conjugacy between the shift map $\sigma:\Sigma_M\to\Sigma_M$
and $f_P|_{\Lambda_P}:\Lambda_P\to\Lambda_P$. Under these circumstances,
$h(f_P)=\log\lambda$, where $\lambda$ is the spectral radius of $M$.
Using your favorite matrix computation program, you may check that
$\lambda=1.46557$ and that $\log\lambda=0.382244$, which agrees with table
1 of \cite{H2}. Although this is not a proof,
one may be obtained using the algorithm in \cite{BH} to find the
pseudo-Anosov homeomorphism in the isotopy class of $f$ rel.
$s^6_7(0)$. It is known that this map realizes the minimum topological
entropy in its isotopy class and it is possible to find a Markov Partition
for it with the same transition matrix $M$.

As was mentioned in the introduction, we intend to show in a forthcoming
paper that, as in Example 3, given a horseshoe periodic orbit collection
$\cal O$, there exists a pruning front $P=P(\cal O)$ such that
$f_P$ restricted to $\Omega(f_P)$ is semconjugated to the Thurston minimal
representative $\phi_{\cal O}$ in the isotopy class of $f$ rel. $\cal O$
and that $h(f_P)=h(\phi_{\cal O})$.

%% file: biblio.tex
\bibliographystyle{plain}